\documentclass{cornell}
\usepackage{latexsym}
\usepackage{epsfig}
\usepackage{graphicx}
\usepackage{amsfonts}
\usepackage{amsmath}
\usepackage{amsthm}
\usepackage{amscd}
\usepackage{amssymb,mathrsfs}

\def\R{\mathbb R}
\def\F{ F}
\def\G{G}
\def\P{P}
\def\C{\mathcal C}
\def\A{\mathcal A}
\def\B{\mathcal B}
\def\CC{\mathcal C}
\def\FF{\mathcal F}
\def\S{\mathbb S}
\def\H{\mathbb H}
\def\E{\mathbb E}

\def\e{\varepsilon}
\def\int{\mathrm{int}}
\def\cl{\mathrm{cl}}
\def\st{\mathrm{st}}
\def\dim{\mathrm{dim }}
\def\lk{\mathrm{lk}}
\def\vol{\mathrm{vol}}
\def\D{\Delta}
\def\glue{\oplus}

\def\spanr{\mathrm{span}}
\def\alphah{\widetilde\alpha}
\def\f{\bar f}
\def\compname{semi-constructible }

\newtheorem{theorem}{Theorem}[section]
\newtheorem{proposition}[theorem]{Proposition}

\newtheorem{corollary}[theorem]{Corollary}
\newtheorem{conjecture}[theorem]{Conjecture}
\newtheorem{lemma}[theorem]{Lemma}

\numberwithin{equation}{section}

\begin{document}
\title {Angle Sums on Polytopes and Polytopal Complexes}
\author{Kristin A. Camenga}
\conferraldate{August}{2006}
\maketitle

\makecopyright

\begin{abstract}

%The interior angle at a face of a convex polytope is the fraction of directions from the face that go into the interior of the polytope.  The $i$th angle sum of the polytope is the sum of the interior angle of all $i$-faces of the polytope. 
We will study the angle sums of polytopes, working to exploit the analogy between the $f$-vector of faces in each dimension and the $\alpha$-vector of angle sums. %by using the well-developed study of the $f$-vector. 
The Gram relation on the $\alpha$-vector is analogous to the Euler relation on the $f$-vector.  Similarly, the Perles relations on the angle sums of simplicial polytopes are analogous to the Dehn-Sommerville relations.  

First we describe the spaces spanned by the angle sums of certain classes of polytopes, as recorded in the $\alpha_{}$-vector and the $\alpha$-$f$-vector.  Families of polytopes are constructed whose angle sums span the spaces of polytopes defined by the Gram and Perles equations. This shows that the dimension of the affine span of the space of angle sums of simplices is $\left\lfloor \frac{d-1}{2}\right\rfloor$, and that of the combined angle sums and face numbers of simplicial polytopes and general polytopes are $d-1$ and $2d-3$, respectively.  %A tool used in proving these results is the $\gamma$-vector, an angle analog to the $h$-vector. 

Next we consider angle sums of polytopal complexes.  We define the angle characteristic on the $\alpha$-vector in analogy to the Euler characteristic.  Then we consider the effect of a gluing operation to construct new complexes on the angle and Euler characteristics.  We show that the changes in the two correspond and that, in the case of certain odd-dimensional polytopal complexes, the angle characteristic is half the Euler characteristic.  In particular, we show that many non-convex spheres satisfy the Gram relation and handle-bodies of genus $g$ constructed via gluings along disks have angle characteristic $1-g$.

Finally, we consider spherical and hyperbolic polytopes and polytopal complexes. Spherical and hyperbolic analogs of the Gram relation and a spherical analog of the Perles relation are known, and we show the hyperbolic analog of the Perles relations in a number of cases.  Proving this relation for simplices of dimension greater than 3 would finish the proof of this result.  Also, we show how constructions on spherical and hyperbolic polytopes lead to corresponding changes in the angle characteristic and Euler characteristic.  However, the angle characteristic and Euler characteristic do not have the 1:2 ratio that held for Euclidean polytopal complexes.
\end{abstract}

\begin{biosketch}
Kristin A. Camenga has been interested in mathematics and mathematics education from an early age.  Her mother, Althea Rood, is a secondary mathematics teacher and inspired in her both a love for mathematics and a passion for teaching it.  She was born November 30, 1976 in New London, Connecticut and spent her elementary years in Connecticut and Rhode Island.  A 1993 graduate of Triton High School in Dodge Center, Minnesota, Kristin planned to teach high school mathematics.  In 1997 she graduated from St. Olaf College in Northfield, Minnesota with a B.A. in mathematics and a concentration in computer science, completing her secondary mathematics teaching certification in the winter of 1997.  During 1998, she studied in Budapest, Hungary on a Fulbright Scholarship.  During this time, she studied graph theory and the Hungarian mathematics education system, observing a number of classrooms and studying the national curriculum.  Returning to Minnesota, she taught at Delano High School for two years, teaching Geometry, Algebra II, and Calculus II.  

In the summer of 2000, she married Andrew J. Camenga and moved to New York, where she began study for her doctorate in  Mathematics.  Kristin's combined interests in mathematics and education also led to the completion of an M.S. in Education in January 2006.  Her research focused on the variety of ways writing can be used in mathematics instruction. Beginning in Fall 2006, she will be teaching in the mathematics department at Houghton College in Houghton, New York, where she will combine her love for mathematics and education by teaching mathematics and working with pre-service teachers.  
\end{biosketch}

\begin{dedication}
To Andrew
\end{dedication}

\begin{acknowledgements}

I am very thankful to my advisor, Lou Billera, for guiding me and encouraging me in this research.  His wisdom in helping me to develop these problems and his knowledge of the general structure of work on polytopes and resources I might use have been invaluable.  I also wish to thank David  Henderson, Ed Swartz and Bob Connelly who have all advised me in this research, helping to locate references and pointing me toward ideas that have extended my research. John Hubbard also inspired some of this research by asking the right question at the right time.

I am thankful for my family, who have continually encouraged me as I pursued my Ph.D. and pushed me to try in the first place - believing in me much more than I believe in myself!  I am particularly grateful to my husband, Andrew, who has put up with me when I let other things go and encouraged me to dig back in when no progress was being made.  

Many friends also supported me in this process.  My best Friend created me with the necessary abilities in the first place and has provided so many opportunies along the way.  I am thankful for my church family who have believed in me and put up with regular progress reports. I am especially thankful to the other graduate students at Cornell who have encouraged me and spurred me on - still pushing me to do better!  In particular, I thank Josh Bowman, Lee Gibson, and Melanie Pivarski who frequently listened to me talk, shared ideas, and answered questions.  By the way Melanie - I think we both won the race!
\end{acknowledgements}

\contentspage

%\tablelistpage

\figurelistpage 

\normalspacing
\setcounter{page}{1}
\pagenumbering{arabic}
\pagestyle{cornell}

%%Thesis Introduction
\chapter{Introduction}
One of the motivating questions in combinatorics is whether a class of objects can be characterized by a set of combinatorial data. For polytopes and complexes, this has been studied in depth by considering the number of faces of each dimension and the inclusion structure on these faces.  We will seek to describe polytopes and complexes by studying their angle sums, which quantify some of the geometric aspects of these structures.  We will first consider Euclidean polytopes and then expand our consideration to complexes and polytopes in other geometries.

\section{Angle Sums and the Gram Relation}
Let $\P$ be a polytope, defined as the convex hull of vertices in $\R^d$, or, equivalently, as the bounded intersection of hyperplanes.  Polytopes considered here are assumed to be convex.  We will assume that $\P$ is a $d$-polytope, that is, the affine span of $\P$ is $d$-dimensional.  If $H$ is a hyperplane whose intersection with $\P$ is contained in $\partial\P$, we say that $\F = H \cap \P$ is \emph{face} of $\P$ and $H$ is a \emph{supporting hyperplane} of $\P$.  A face of dimension $i$ is called an $i$-face and $(d-1)$-faces are called \emph{facets}. We define $f_i(\P)$ as the number of $i$-faces of $\P$ for $i=0,\ldots, d$. Then $f_d(\P) = 1$ for all polytopes.  By convention, we define $f_{-1}(\P) = 1$ and think of this as counting the empty face. The \emph{$f$-vector} of a polytope $\P$ is $(f_{0}, f_{1}, f_{2}, \dots, f_{d})$.

The \emph{interior angle} at a face $\F \subset \P$ is defined by 
\[\alpha(\F, \P) = \frac{\vol(S_{\e}(x) \cap \P)}{\vol(S_{\e}(x))},\] 
where $x$ is in the relative interior of $\F$ and $S_\e(x)$ is the $(d-1)$-sphere of radius $\e$ centered at $x$ for $\e$ small enough to only intersect $\P$ in faces that contain $\F$. %%PICTURE?
Therefore, the interior angle measures the fraction of directions that one can move from the face $\F$ into $\P$.  The \emph{angle sums} of $\P$ are defined for $0\leq i\leq d$ as
\[\alpha_{i}(\P)= \sum_{i-\text{faces } \F\subseteq \P} \alpha(\F, P).\]
Since $\P$ is the only $d$-face of the polytope, $\alpha_d(\P) = 1$ for all polytopes.

%%EXAMPLE: CUBE

For example, we can consider a standard cube, $C$.  At any vertex $v$ of $C$, one-eighth of the directions from $v$ go into the cube.  Another way to think of this is that exactly eight cubes could be put together at $v$ to surround $v$.  Therefore, we say that the interior angle at $v$ is $\alpha(v, C) = \frac{1}{8}$.  In the same way, the angle at any edge $e$ is $\alpha(e, C) = \frac{1}{4}$ and the angle at any face $f$ is $\alpha(f, C) = \frac{1}{2}$. Since the cube is regular, the angle at each vertex, edge, and face is the same and we have the following angle sums for the cube: 
\[\alpha_{0}(C) = 8\Big(\frac{1}{8}\Big) = 1 \quad 
\alpha_{1}(C) = 12\Big(\frac{1}{4}\Big) = 3 \quad
\alpha_{2}(C) = 6\Big(\frac{1}{2}\Big) = 3 \quad
\alpha_{3}(C) = 1.\] 

Like the $f$-vector, we define the 
\emph{$\alpha{}$-vector},
\[\alpha(\P)=(\alpha_{0}(P), \alpha_{1}(P), \ldots, \alpha_{d}(P)),\] and the \emph{$\alpha{}$-$f$-vector}, 
\[\alpha\text{-}f(\P) = (\alpha_{0}(P), \alpha_{1}(P), \ldots ,
\alpha_{d}(P),f_{0}(P), f_{1}(P), \ldots , f_{d}(P)).\]
Occasionally, we will write 
\[\alpha\text{-}f(\P) = (\alpha_{0}(P), \alpha_{1}(P), \ldots ,
\alpha_{d}(P)|f_{0}(P), f_{1}(P), \ldots , f_{d}(P))\] 
to clarify where the angle sums end and face numbers begin.
The basic relation on the $\alpha$-vector is the following:
\begin{theorem}\label{Gram's Relation}
\textbf{Gram Relation}
For any $d$-polytope $\P$,
\[\sum_{i=0}^{d-1} (-1)^{i}\alpha_{i}(\P) = (-1)^{d-1}.\]
\end{theorem}
This relation is sometimes called the Gram-Sommerville relation \cite{Chen1} or the Brianchon-Gram relation \cite{McM2}.   

In the case of the cube, we can see that 
\[\alpha_{0}(C) - \alpha_{1}(C) + \alpha_{2}(C) = 1-3+3 = 1 = (-1)^{3-1}.\]
The 2-dimensional case of the Gram relation corresponds to the high school geometry theorem for the sum of angles in a polygon.  The 3-dimensional case was proved by Gram in 1874 \cite{G}.  Gr\"unbaum published the first accepted proof for dimension $d$ in 1967 \cite{Gr}, which we will follow below. Before beginning the proof, we note that the Gram relation is reminiscent of the Euler relation, which is the only linear relation on the face numbers of polytopes:
\begin{theorem}\label{Euler Relation}
\textbf{Euler Relation}
Let $\P$ be a $d$-polytope, and let $f_{i}(\P)$ be the number of faces of $\P$ of dimension $i$.  Then
\[\sum_{i=0}^{d-1} (-1)^{i}f_{i}(\P) = 1 + (-1)^{d-1}.\]
\end{theorem}
This will be needed in the proof of Gram's relation.  Gr\"unbaum also proves that the Gram relation is the only linear relation on angle sums using the fact that the Euler relation is the only linear relation on the $f$-vector \cite{Gr}. %p.304

More generally, for a collection of faces of maximal dimension $d-1$, we define the $f$-vector as we did for polytopes, $(f_{0}, f_1, \ldots, f_{d-1})$,  with $f_i$ counting the number of $i$-faces.  
Then the \emph{Euler characteristic} $\chi(\C)$ is: 
\[\chi(\C) = \sum_{i=0}^{d-1} (-1)^{i}f_i(\C).\] 
This will be needed in the proof of the Gram relation.  We also need to define the \emph{pyramid} on a $(d-1)$-polytope $Q$, the $d$-polytope formed by taking the convex hull of $Q$ and a point not in the affine span of Q. Any polytope formed in this manner is called a pyramid, and if a $d$-polytope can be constructed by taking the pyramid over an $(d-m)$-face $m$ times, we will call this an $m$-fold $d$-pyramid.

\begin{proof}[Proof of Theorem \ref{Gram's Relation} \cite{Gr}] 
We proceed in three steps: first we prove the relation for simplices; secondly, we decompose the polytope into pyramids over the facets from an interior point and show that if the pyramids satisfy the Gram relation, so does the whole polytope; thirdly, we show that any pyramid satisfies the Gram relation, and therefore it follows from the second step that all polytopes do.

\noindent \emph{PART I: Gram's Relation for Simplices}

Let $\D$ be a $d$-simplex.  For each of the facets $\F_1, \F_2,\dots \F_{d+1}$, let $C_i$ denote the half-space of $\R^d$ bounded by the hyperplane which is the affine span of $\F_i$ and containing $\D$. Since $\D$ is a simplex, for each $m$-tuple $(i_1, i_2, \ldots i_m)$ of supporting hyperplanes, there is some $(d-m)$-face $\F$ so that $\displaystyle{\bigcap_{j=1}^{m} \F_{i_j} = \F}$.  This set of facets is precisely the set of facets that contain the face $\F$.  If we consider the sphere, $S_{\F}$, used to determine the interior angle at a $(d-m)$-face $\F$ which is defined by an $m$-tuple of facets given by indices the $(i_1, i_2, \ldots i_m)$, then it is clear that 
\begin{equation}\label{facet intersections}
\bigcap_{j=1}^{m}\Big(C_{i_j} \cap S_{\F}\Big)=\Big(\bigcap_{j=1}^{m} C_{i_j}\Big) \cap S_{\F} = \D \cap S_{\F}.
\end{equation} 
That is, the part of the sphere determined by the angle at $\F$ is the intersections of the hemispheres determined by the facets that contain $\F$.  

To allow comparisons between the interior angles at different faces, we will define $V(\F)$ for each face $\F$ as the set of unit vectors in the unit sphere whose directions point from $\F$ into the polytope $\P$.  This could be considered a resizing of the polytope so that a unit sphere can be used to determine each interior angle.  Then $\displaystyle{\alpha(\F, \D) = \frac{\vol(V(\F))}{\vol(S^{d-1})}}$, where $S^{d-1}$ is the unit sphere.  However, the orientation of the unit sphere is fixed so that the angles at each face are represented by a corresponding portion of the unit sphere.  Then \eqref{facet intersections} can be rewritten as
\begin{equation}\label{interiorizingint}
\bigcap_{j=1}^{m}V(\F_{i_j}) = V\left(\bigcap_{j=1}^{m}\F_{i_j}\right).
\end{equation}

The argument is then carried out based on the principle of inclusion-exclusion, where we consider the unit $(d-1)$-sphere as the union of the hemispheres given by interior angles at the facets:  
\begin{align*}
\vol(S^{d-1}) & = \vol\left(\bigcup_{i=1}^{d+1} V(\F_i)\right)\\
& = \sum_{i=1}^{d+1} \vol\left(V(\F_i)\right) - \sum_{i, j = 1, i \neq j}^{d+1}\vol\left(V(\F_i) \cap V(\F_j)\right) + \cdots\\
& \qquad \qquad  + (-1)^{d}\vol\left(V(\F_1)\cap V(\F_2)\cap \cdots \cap V(\F_{d+1})\right)\\
& = \sum_{i=1}^{d+1} \vol\left(V(\F_i)\right) - \sum_{i, j = 1,i \neq j}^{d+1}\vol\left(V(\F_i \cap \F_j)\right) + \cdots \\
& \qquad \qquad + (-1)^{d}\vol\left(V(\F_1 \cap \F_2\cap \cdots \cap \F_{d+1})\right) \quad \text{by }\ref{interiorizingint}\\
& = \sum_{(d-1)-\text{faces } \F} \vol\left(V(\F)\right) - \sum_{(d-2)-\text{faces } \F}\vol\left(V(\F)\right) + \cdots\\
& \qquad \qquad+ (-1)^{d}\sum_{(-1)-\text{faces } \F}\vol\Big(V(\F)\Big). 
\end{align*}
Dividing both sides by $\vol(S^{d-1})$ and rewriting ratios of volumes as interior angles, we get that
\begin{align*}
1 &= \sum_{(d-1)-\text{faces } \F} \alpha(\F, \D) - \sum_{(d-2)-\text{faces } \F}\alpha(\F, \D) + \cdots \\
& \qquad \qquad \qquad + (-1)^{d-1} \sum_{0-\text{faces } \F} \alpha(\F, \D)+(-1)^{d}\sum_{(-1)-\text{faces } \F}\frac{\vol(\emptyset)}{\vol(S^{d-1})}\\
& = \alpha_{d-1}(\D) - \alpha_{d-2}(\D) + \cdots + (-1)^{d-1} \alpha_{0}(\D) + 0.
\end{align*}
Multiplying both sides by $(-1)^{d-1}$ gives the Gram relation for simplices.  This also gives some rationale for the convention $\alpha_{-1}(\P)=0$, since it corresponds to the normalized volume of the empty set.

\noindent \emph{PART II: Decomposition of a Polytope into Pyramids from an Interior Point}

Let $\P$ be a $d$-polytope with $f_{d-1}(\P) = e$ facets, $\F_1, \F_2, \ldots \F_e$, and let 0 be an interior point of $\P$.  Define $\P_i$ as the $d$-pyramid constructed as the convex hull of $\F_i$ and 0.  Then $\P$ is the union of the $\P_i$, with all intersection occurring on faces of the $\P_i$.  We claim that if each of the $\P_i$ satisfies the Gram relation, then $\P$ does.

This is dependent on the additive nature of interior angles.  First of all, we note that any $k$-face of some $\P_i$ that is in the interior of $\P$ will contribute a total of 1 to the sum $\displaystyle{\sum_{i=1}^{e}\alpha_k(\P_i)}$ since any direction from the face will be into the interior of $\P$ and therefore into one of the $\P_i$.  Also, for any $k$-face on the boundary of $\P$, the sum of the angles in the $\P_i$ including that face will equal the angle at that face in $\P$.  Since each of the $k$-faces in the interior is constructed as the convex hull of 0 and a $(k-1)$-face of $\P$, we have
\[\sum_{i=1}^{e}\alpha_k(\P_i) = \alpha_k(\P) + f_{k-1}(\P).\]
In particular, $\displaystyle{\sum_{i=1}^{e}\alpha_0(\P_i) = \alpha_0(\P) +1}$ since 0 is in the interior of $\P$.

Using the convention $f_{-1}(\P) = 1$, the following computation holds:
\begin{align*}
\sum_{k=0}^{d-1} (-1)^{k}\alpha_k(\P) & = \sum_{k=0}^{d-1}(-1)^{k}\sum_{i=1}^{e} \alpha_k(\P_i) - \sum_{k=0}^{d-1}(-1)^{k}f_{k-1}(\P)  \\ 
& =  \sum_{i=1}^{e}\sum_{k=0}^{d-1}(-1)^{k} \alpha_k(\P_i) + \sum_{k=-1}^{d-2}(-1)^{k}f_{k}(\P) \\
& =  \sum_{i=1}^{e}(-1)^{d-1}+ \sum_{k=-1}^{d-2}(-1)^{k}f_{k}(\P) \\
& = \sum_{k=-1}^{d-1}(-1)^{k}f_{k}(\P) = (-1)^{d-1},
\end{align*}
where the third equality follows from the Gram relation on the $\P_i$ and the last from the Euler relation.  Therefore, $\P$ satisfies the Gram relation.

\noindent \emph{PART III: Gram's Relation for all Pyramids}

To complete the proof, it is necessary to show that all of the pyramids over the faces used in the decomposition in Part II satisfy Gram's relation.  We will do this by proving that if Gram's relation holds for $m$-fold $d$-pyramids, it holds for $(m-1)$-fold $d$-pyramids, where $2 \leq m \leq d-1$.  Since Part I established the relation for simplices, or $(d-1)$-fold $d$-pyramids, this induction step will prove the relation for all pyramids and complete the proof. 

This induction step will be accomplished by decomposing a $(d-m)$-pyramid into $(d-m+1)$-pyramids, assuming that the Gram relation holds for the latter and showing that the relation holds for the former as a result. Suppose $\P$ is a $(d-m)$-fold $d$-pyramid where $m \geq 2$.  Let $\F^{*}$ be the $m$-face over which $\P$ is a pyramid.  Let 0 be an interior point of $\F^{*}$, and let $\P_i$ be the convex hull of 0 and a facet of $\P$ which does not contain $\F^{*}$ for $i = 1, 2, \ldots, s$.  Therefore $s$ is the number of facets of $\P$ which do not contain $\F^{*}$.  First of all, we note that $\P_1, \P_2, \ldots \P_s$ is a decomposition of $\P$.  If we take the ray from 0 through any point in $\P$, it must exit through some facet of $\P$ and therefore this point is contained in the pyramid built on that facet. Secondly, we note that each of the $\P_j$ is a $(d-m+1)$ pyramid since each facet had been a $(d-m)$-fold pyramid over some proper face of $\F^{*}$ and taking the convex hull with 0 is another iteration of the pyramid operation.  

For each $n > m$, denote the $n$-faces of $\P$ that contain $\F^{*}$ as $\F_i^n$ for $i \in I(n)$.  Let $\CC^n$  be the complex consisting of the $\F_i^n$ for $i \in I(n)$ and all their faces and let $\CC(\G)$ be the complex consisting of the face $\G$ and all its proper faces.  We also define the \emph{star} of a face $\G$ in a complex $\CC$ as $\st(\G; \CC) = \{\G' \in \CC: \G \subseteq \G'\}$.  Then for $0 \leq k \leq n-1$, 
\begin{equation}\label{stars and complexes}
f_k(\CC(\F_i^n) \setminus \CC^{n-1}) = f_k(\F_i^n) - f_k(\st(\F^{*}; \F_i^n));
\end{equation}
that is, the number of $k$-faces included in $\F_i^n$ that are not included in any $(n-1)$-dimensional face which includes $\F^{*}$ is the difference between the number of $k$-faces of $\F_i^n$ and the $k$-faces of the star of $\F^{*}$ in $\F_i^n$.  In particular, $s = f_{d-1}(\CC(\P) \setminus \CC^{d-1})$ where $\P = \F_1^d$. Using the Euler relation and the fact that the Euler characteristic of the star of a face in a polytope is 1 \cite{Gr}, %(GRUNBAUM p. 139 - cite and list theorem??) 
we have that 
\begin{equation}\label{Euler char of star}
\begin{split}
\sum_{k=1}^n (-1)^{k}f_{k-1}(\CC(\F_i^n) \setminus \CC^{n-1}) &= \sum_{k=0}^{n-1} (-1)^{k}f_k(\st(\F^{*}; \F_i^n)) - \sum_{k=0}^{n-1}(-1)^kf_k(\F_i^n) \\
& = 1 - (1+(-1)^{n-1}) = (-1)^n.
\end{split}
\end{equation}

Now we  consider the relationship between the angle sums of $\P$ and those of the $\P_j$, $1 \leq j \leq s$.  First, we note that $\displaystyle{\sum_{j=1}^s \alpha_0(\P_j)= \alpha_0(\P) + \alpha(\F^{*}, \P)}$.  This follows since the sum of the angles at vertices of the $\P_i$ not only count the angles at vertices of $\P$, but it also counts the angle at the new vertex, 0.  At this vertex, the total angle will sum to the interior angle at $\F^{*}$.  

For the quantity $\displaystyle{S_k=\sum_{j=1}^s \alpha_k(\P_j)}$ for $k > 0$, we can similarly divide this sum between the portion that counts angles at the $k$-faces that do not contain 0 and that which counts angles at $k$-faces that do contain 0.  If a face $\G$ does not contain 0, it is a $k$-face of $\P$ and $\displaystyle{\alpha(\G, \P) = \sum_{j=1}^s \alpha(\G, \P_j)}$ and these faces contribute exactly the same amount to $S_k$ as they do to $\alpha_k(\P)$.  Let $H(\G)$ be the uniquely determined face of $\P$ of smallest possible dimension which contains $\G$.  $H(\G)$ may be $\P$.  If a face $\G$ does contain 0, it was created as the convex hull of 0 and a $(k-1)$-face of $\P$ and $\F^{*} \subseteq H(\G)$ since 0 is in the relative interior of $\F^{*}$.  Then we have that 
\[ \alpha(H(\G), \P) = \sum_{\G \subseteq \P_j} \alpha(\G, \P_j)\] since the disjoint union of all the directions in the various $\P_j$ from $\G$ will equal the directions from $H(\G)$ into $\P$.  
Also, $\alpha(\G, \P) = \alpha(H(\G), \P)$ so we can write 
\begin{equation}\label{SK}
S_k = \sum_{0 \notin \G} \alpha(\G, \P) +  \sum_{0 \in \G} \alpha(H(\G), \P),
\end{equation}
where each $k$-face $\G$ of some $\P_j$ is considered once in the sums.
In the first sum of \eqref{SK}, the only faces of $\P$ we are not summing over are exactly those $k$-faces of $\P$ that contain $\F^{*}$.  Therefore, \[\sum_{0 \notin \G} \alpha(\G, \P) = \alpha_k(\P) - \sum_{i \in I(k)} \alpha(F_i^k, \P).\]

We can rewrite the second sum of \eqref{SK} by indexing according to $H(\G)$.  Thinking of $H(\G)$ as a fixed $F_i^n$, we consider all the $(k-1)$-faces of this $F_i^n$ for which $F_i^n$ is minimal according to inclusion among all faces of $\P$ that include $\F^{*}$.  This will then count all the $(k-1)$-faces of $F_i^n$ which do not belong to any $F_{i'}^{n-1}$, $i' \in I(n-1)$.  Then we can write
\begin{align*}
\sum_{0 \in \G} \alpha(H(\G), \P) & = \sum_{i,n} \sum_{\G: H(\G) = F_i^n} \alpha(H(\G), \P) \\
& = \sum_{n = \max(m,k)}^{d} \sum_{i\in I(n)} \alpha(F_i^n, \P) f_{k-1}(\CC(F_i^n) \setminus \CC^{n-1}).
\end{align*}
Putting these results together and writing  $l = \max(m,k)$, for $k > 0$ we have 
\begin{align*}
S_k = \sum_{j=1}^s \alpha_k(P_j) & = \sum_{0 \notin \G} \alpha(\G, \P) +  \sum_{0 \in \G} \alpha(H(\G), \P) \\
& = \Big(\alpha_k(\P) - \sum_{i \in I(k)} \alpha(F_i^k, \P)\Big) \\
& \qquad\qquad\qquad + \sum_{n = l}^{d} \sum_{i\in I(n)} \alpha(F_i^n, \P) f_{k-1}(\CC(F_i^n) \setminus \CC^{n-1}) \\
& = \alpha_k(\P) + \sum_{n = l}^{d} \sum_{i\in I(n)} \alpha(F_i^n, \P) \Big(f_{k-1}(\CC(F_i^n) \setminus \CC^{n-1}) - \delta_{nk}\Big),
\end{align*}
where $\delta_{nk}$ is the Kronecker delta.

Now we apply the assumption that the theorem holds for each $P_j$.
\begin{align*}
\sum_{k=0}^{d-1}(-1)^k \alpha_k (\P) & = 
\left(\sum_{j=1}^s \alpha_0(\P_s)- \alpha(\F^{*}, \P)\right)+ \sum_{k=1}^{d-1}(-1)^k\Bigg[\sum_{j=1}^s \alpha_k(P_j) \\
& \qquad \qquad- \Bigg( \sum_{n = l}^{d} \sum_{i\in I(n)} \alpha(\F_i^n, \P) \Big(f_{k-1}(\CC(\F_i^n) \setminus \CC^{n-1}) - \delta_{nk}\Big)\Bigg)\Bigg]\\
 & = \sum_{j=1}^s \sum_{k=0}^{d-1}(-1)^k \alpha_k(P_j) - \alpha(\F^{*}, \P) -  \sum_{k=1}^{d-1}(-1)^k\\
& \qquad \qquad \times \Bigg( \sum_{n = l}^{d} \sum_{i\in I(n)} \alpha(\F_i^n, \P) \Big(f_{k-1}(\CC(\F_i^n) \setminus \CC^{n-1}) - \delta_{nk}\Big)\Bigg)\\
& = s(-1)^{d-1} - \alpha(\F^{*}, \P) - \sum_{n = m}^{d-1} \sum_{i\in I(n)} \alpha(\F_i^n, \P)\\
& \qquad  \qquad \times \left[\sum_{k=1}^{d-1}(-1)^{k}f_{k-1}(\CC(\F_i^n) \setminus \CC^{n-1}) - (-1)^n\right] \\
& \qquad \qquad - \alpha(\P, \P)\sum_{k=1}^{d-1}(-1)^{k}f_{k-1}\left(\CC(\P) \setminus \CC^{d-1}\right),
\end{align*}
where in the last equality we simplify the first sum using the Gram relation, separate the $n=d$ term of the last sum, and switch summands of the last sum, remembering that if $k-1 > n$, then $f_{k-1}(\CC(\F_i^n) \setminus \CC^{n-1})=0$ so it adds nothing to the sum.  Then by the definition of $s$ and using that $|I(m)| =1$ and $\F_1^m = \F^{*}$,  
\begin{align*}
\sum_{k=0}^{d-1}(-1)^k& \alpha_k (\P) \\
&= (-1)^{d-1}f_{d-1}(\CC(\P) \setminus \CC^{d-1})- \alpha(\F^{*}, \P)\\
& \qquad  -\alpha(\F^{*}, \P)\left[\sum_{k=1}^{d-1}(-1)^{k}f_{k-1}\left(\CC(\F^{*})\right) - (-1)^m\right]\\
& \qquad -\sum_{n = m+1}^{d-1} \sum_{i\in I(n)} \alpha(\F_i^n, \P)\left[\sum_{k=1}^{d-1}(-1)^{k}f_{k-1}\left(\CC(\F_i^n) \setminus \CC^{n-1}\right) - (-1)^n\right] \\
& \qquad - \sum_{k=1}^{d-1}(-1)^{k}f_{k-1}\left(\CC(\P) \setminus \CC^{d-1}\right)\\
& = (-1)^{d-1}f_{d-1}\left(\CC(\P) \setminus \CC^{d-1}\right)- \alpha(\F^{*}, \P) \\
& \qquad - \alpha(\F, \P)\Big[-(1-(-1)^m) - (-1)^m\Big] \\
& \qquad - \sum_{n = m+1}^{d-1} \sum_{i\in I(n)} \alpha(\F_i^n, \P)\Big[(-1)^n - (-1)^n\Big] \\
& \qquad  + \sum_{k=0}^{d-2}(-1)^{k}f_{k-1}\left(\CC(\P) \setminus \CC^{d-1}\right)\\
& = \sum_{k=0}^{d-1}(-1)^{k}f_{k-1}\left(\CC(\P) \setminus \CC^{d-1}\right) \\
& = (-1)^{d-1} \quad \text{by } \eqref{Euler char of star}.
\end{align*} 
\end{proof}

In Gr\"unbaum's proof of Gram's theorem, the need for the Euler relation is clear.  Summing over polytopes in the decomposition leads to use of the Euler relation on the faces of the polytope.  In fact, Parts II and II of the proof are dependent only on the inclusion structure of the polytope and the Euler relation.  This fact will be used later to generalize this proof.  Shephard \cite{S} and, later, Welzl \cite{W} gave a more explicit connection between the Gram and Euler relations when they proved a connection between the two by considering a projection of a $d$-polytope $\P$ onto a $(d-1)$-polytope $\P'$.  Since the directions which lie along faces of the polytope account for a set of measure zero, it can be shown that for any $i$-face $\F$ of $\P$, $i\leq d-2$, \mbox{$\mathrm{Prob}(\F'$ is a face of $\P') = 1-2\alpha(\F,\P)$}.  Then
\begin{equation}\label{projection}
\alpha_{i} (\P) =\frac {f_{i}(\P) - E(f_{i}(\P'))}{2} \quad \text{for } 0 \leq i \leq d-2.
\end{equation}  
This allows for an easy translation between equations on $f$-vectors and those on $\alpha$-vectors.  Therefore, the following computation gives another proof of the Gram relation based on the Euler relation.
\begin{equation*}
\begin{split}
\sum_{i=0}^{d-1} (-1)^{i}\alpha_{i}(\P) & =  \sum_{i=0}^{d-2} \left[ (-1)^{i}\frac{f_{i}(\P) -
E(f_{i}(\P'))}{2} \right] +(-1)^{d-1} \frac {f_{d-1}(\P)}{2}
\\ & =  \frac {1}{2} \sum_{i=0}^{d-1} (-1)^{i} f_{i}(\P) - \frac {1}{2}
\sum_{i=0}^{d-2} E(f_{i}(\P')) \\ &
=  \frac {1}{2} (1+(-1)^{d-1}) - \frac {1}{2} (1 + (-1)^{d-2})\quad \text{ by the Euler relation}\\
& =  (-1)^{d-1}.
\end{split}
\end{equation*}

\section{Relations on Simplicial Polytopes}

A polytope $\P$ is \emph{simplicial} if all of its facets are simplices.  There are more relations on the $f$-vectors of simplicial polytopes than just the Euler relation.
\begin{theorem}\label{Dehn-Sommerville Relations}
\textbf{Dehn-Sommerville Relations}
For any simplicial polytope $\P$ and $-1 \leq k \leq d-2$, 
\[\sum_{j=k}^{d-1} (-1)^{j} {j+1\choose k+1} f_{j}(\P) =
(-1)^{d-1}f_{k}(\P).\]
\end{theorem}
The Dehn-Sommerville relation for $k=-1$ agrees with the Gram relation.  We can do a similar translation using \eqref{projection} for simplicial polytopes.  This depends on the fact
that if $\P$ is simplicial, its generic projection $\P'$ is also simplicial. Then we can determine the analogous relations for angle sums:
\begin{equation*}
\begin{split}
\sum_{j=k}^{d-1} (-1)^{j} {j+1\choose k+1} \alpha_{j}(\P) & =  \sum_{j=k}^{d-2}
(-1)^{j} {j+1\choose k+1} \left[ \frac {f_{i}(\P) - E(f_{i}(\P'))}{2} \right] \\ &
\qquad \qquad +(-1)^{d-1} {d\choose k+1} \frac {f_{d-1}}{2} \\ & =  \frac {1}{2}
(-1)^{d-1} f_{k}(\P) - \frac {1}{2} (-1)^{d-2} E(f_{k}(\P'))
\end{split}
\end{equation*}
by the Dehn-Sommerville relations.  Then 
\begin{equation*}
\begin{split}
\sum_{j=k}^{d-1} (-1)^{j} {j+1\choose k+1} \alpha_{j}(\P)& = \frac {1}{2}(-1)^{d-1} \left[ f_{k}(\P) + f_{k}(\P) - 2\alpha_{k}(\P))\right] \\ & =
(-1)^{d}(\alpha_{k}(\P) - f_{k}(\P)).
\end{split}
\end{equation*}
If we follow the conventions $f_{-1}(\P)=1$ and $\alpha_{-1}(\P) = 0$ (for $\dim(\P) \geq 1$) then we still have $E(f_{-1}(\P')) = 1 = f_{-1}(\P) - 2\alpha_{-1}(\P)$ and this proves the following theorem for $-1 \leq k \leq d-2$, earlier proved in \cite{Gr} in a manner analogous to Gr\"unbaum's proof of the Gram relation:
\begin{theorem}\label{Perles Relations}
\textbf{Perles Relations}
For any simplicial polytope $\P$ and $-1 \leq k \leq d-2$, 
$$\sum_{j=k}^{d-1} (-1)^{j} {j+1\choose k+1} \alpha_{j}(\P) =
(-1)^{d}(\alpha_{k}(\P) - f_{k}(\P)).$$ \hfill $\Box$
\end{theorem}
The Gram relation is the $k=-1$ case of the Perles relations. For $k=d-1$, the Perles relations give $\alpha_{d-1}(\P) = \frac{1}{2} f_{d-1}(\P)$, which is true since every facet contributes $\frac{1}{2}$ to $\alpha_{d-1}(\P)$. The Dehn-Sommerville relations are the only linear relations on the $f$-vectors of simplicial polytopes \cite{Gr}. %p.146 - dim \floor(d/2)
We will show in Chapter 2 that the Perles relations and the Dehn-Sommerville relations are the only linear relations on the $\alpha$-$f$-vectors of simplicial polytopes.
 
With such close parallels between relations on the $f$-vector and those on the $\alpha$-vector, it is reasonable to consider other combinatorial results and possible analogs on angle sums.  In the next two sections we introduce two areas of study of the $f$-vector that we will explore the angle sum analogs of in later chapters.
\section{Affine Spans of $f$-Vectors}

    It is known that the Euler relation and the Dehn-Sommerville relations are the only linear relations on the $f$-vectors of general and simplicial polytopes, respectively.  This means that the affine span of the $f$-vectors of $d$-polytopes is $(d-1)$-dimensional and the affine span of the $f$-vectors of simplicial $d$-polytopes is $\left\lfloor \frac{d}{2}\right\rfloor$-dimensional.

    This can be shown either by demonstrating that any other linear relations are implied by these relations or by providing a set of polytopes that span the affine space defined by the relations. Gr\"unbaum \cite{Gr} shows that a set of cyclic polytopes can be used to span the affine hyperplane determined by Euler's relation.  %p.147
 Bayer and Billera \cite{BB} use pyramid and stellar subdivision constructions to construct two sets of polytopes, one of which spans the affine space defined by the Euler relation and the other of which spans the affine space defined by the Dehn-Sommerville relations.

The pyramid over a polytope $Q$, which we will write as $PQ$, was defined earlier.  We define a point $x$  to be \emph{beyond} a facet of a polytope $Q$ if $x$ and $Q$ lie on opposite sides of the hyperplane which is the affine span of the facet. Let $Q$ be a simplicial $d$-polytope with proper face $\F$.  The \emph{stellar subdivision} of a face $\F$ in $Q$, $\st^{*}(\F,Q)$, is the simplicial $d$-polytope which is the convex hull of $Q \cup \{x\}$, where $x$ is beyond exactly those facets which contain $\F$. 
%In the case that $\F$ is a facet, the stellar subdivision of $\F$ in $Q$ is equivalent to attaching a shallow pyramid to $\F$ and deleting $\F$. 
%PICTURE?

If $Q$ is a $d$-simplex and $\F$ is a $(d-k)$-face of $Q$, $\st^{*}(\F,Q)$ will be denoted $T_{k}^{d}$.  In this case, the new vertex $x$ is beyond exactly $k$ of the hyperplanes defined by facets. For convenience, the $d$-simplex is denoted as $T_0^{d}$.  Then the $f$-vectors of the simplicial polytopes $T_{k}^{d}$, for $0 \leq k \leq \lfloor \frac{d}{2}\rfloor$ span the affine space defined by the Dehn-Sommerville relations \cite{BB}. %p.218

Define $T_{k}^{d,r}$ as the $r$-fold pyramid over the $(d-r)$-polytope $T_k^{d-r}$, where $0 \leq r \leq d-2$ and $1 \leq k \leq \lfloor \frac{d-r}{2}\rfloor$.   Then the $f$-vectors of the $d$-polytopes $T_1^{d,r}$, $0 \leq r \leq d-2$, and that of the $d$-simplex $T_0^{d}$ span the affine space determined by the Euler relation \cite{BB}. % p. 219

For later reference, we also define a third construction operation used by Bayer and Billera.  The \emph{bipyramid} over a $(d-1)$-polytope $Q$, denoted $BQ$, is the $d$-polytope formed by taking the convex hull of $Q$ and a line segment which meets $Q$ in a relative interior point of each.  Equivalently, we could think of this as taking two copies of the pyramid over $Q$ and identifying their bases so that the new vertices in the pyramid are on opposite sides of the hyperplane %of $\R^{d+1}$ 
which is the affine span of $Q$. Then we have the following values of $f_i$ for $PQ$ and $BQ$ %of pyramids and bipyramids 
\cite{BB,Gr}:
\begin{align*}
f_{i}(PQ) &= f_{i}(Q) + f_{i-1}(Q) \text{ for }i\leq d-1, \\
f_{d}(PQ) &= 1,
\end{align*}
\begin{align*}
f_{i}(BQ) &= f_{i}(Q) + 2f_{i-1}(Q) \text{ for }i\leq d-2, \\
f_{d-1}(BQ) &= 2f_{d-2}(Q)\\
f_{d}(BQ) & = 1.
\end{align*}

Bayer and Billera work with the \emph{$h$-vector}, a linear transformation of the $f$-vector.  The $h$-vector is defined on a simplicial polytope $\P$ as
\[h(\P) = \left(h_0(\P), h_1(\P), \ldots, h_d(\P)\right),\] where
\begin{equation}\label{h-vector}
h_i(\P) = \sum_{j=0}^{i}(-1)^{i-j}{d-j \choose d-i}f_{j-1}(\P).
\end{equation}
This transformation can be inverted in the following way:
\begin{equation}\label{h to f}
f_j(\P) = \sum_{i=0}^{j+1}{d-i \choose d-j-i}h_i(\P).
\end{equation}
Therefore, the linear independence of a set of $f$-vectors is equivalent to the linear independence of the corresponding set of $h$-vectors.  We can also consider the effect of different constructions on the $h$-vector as we have on the $f$-vector.  This will be done in the next chapter.
One of the strengths of the $h$-vector is that it allows a more symmetric reformulation of the Dehn-Sommerville relations, proved in \cite{BB}:
\begin{theorem}\textbf{Dehn-Sommerville Relations}
If $\P$ is a simplicial polytope,  
\[h_i(\P) = h_{d-i}(\P) \quad \text{for } i = 0, \ldots, \left\lfloor \frac{d}{2} \right\rfloor.\]\hfill $\Box$
\end{theorem}

We consider one more construction on polytopes.  Each polytope $\P$ has an associated \emph{dual} polytope $\P^{*}$.  Two polytopes, $\P$ and $\P^{*}$, are dual if there is a bijection between the faces of the two that is inclusion-reversing. That is, every vertex of $\P$ corresponds to a facet of $\P^{*}$, every edge of $\P$ corresponds to a codimension 2 face of $\P^{*}$, etc.  Therefore, the $f$-vector of $\P^{*}$ is the reverse of the $f$-vector of $\P$.  In $\R^{3}$ we can see that the cube is dual to the octahedron. %%PICTURE?
Geometrically, we can construct $\P^{*}$ by placing $\P$ so that the origin is in its interior and taking the \emph{polar} of the set of points in $\P$:
\[\P^{*} = \{y \in \R^{d} : \left<x,y\right> \leq 1 \text{  for all } x \in \P\}.\]
Clearly, the dual is an involution on polytopes, so $(\P^{*})^{*}= \P$ \cite{Gr}.

In Chapter 2, we will construct polytopes whose $\alpha$-vectors and $\alpha$-$f$-vectors span the spaces of determined by the Gram and Perles relations.
%Complexes
\section{Euler-type relations for complexes}

The $f$-vectors of complexes have also been widely studied.  A \emph{simplicial complex} is a set of simplices $\C$ such that 
\begin{itemize} 
\item If $\F \in \C$ and $\G$ is a face of $\F$, then $\G \in \C$.

\item If $\F, \G \in \C$, then $\F \cap \G$ is a face of each.
\end{itemize}

Each face will be identified with its set of vertices.  If $\C$ has vertex set $\{v_1, \ldots, v_n\}$, the \emph{geometric realization} of $\C$, $|\C|$, is the union over all faces $\{v_{i_1},\ldots, v_{i_j}\}$ of $\C$ of the convex hull of $\{e_{i_1}, \ldots, e_{i_j}\}$ where $\{e_1, \ldots, e_n\}$ is the standard basis in $R^n$.  Two complexes are said to be homeomorphic if their geometric realizations are.  A simplicial complex is \emph{pure} if the facets (maximal faces) all have the same dimension.

The \emph{link} of a face $\F \in \C$ is the complex
\[\lk(\F, \C) = \{\G: \G \cup \F \in \C, \G \cap \F= \emptyset \}\]

If two complexes are homeomorphic, they have the same Euler characteristic.  Since the Euler relation says that the boundary of a $d$-polytope has Euler characteristic $1+(-1)^{d-1}$, this tells us that any complex $\C$ homeomorphic to a $(d-1)$-sphere has $\chi(\C) = 1+(-1)^{d-1}$.  Therefore, by fixing the homeomorphism type of the complex, the Euler characteristic allows us to determine $f$-vectors which cannot be $f$-vectors of spheres.  On the other hand, when we do not fix the topology of a complex, the Euler characteristic is most frequently used to show that two complexes are not homeomorphic by showing that they have distinct Euler characteristics.  

We define two subclasses of complexes, following the terminology in Swartz \cite{Sw}: semi-Eulerian complexes and homology manifolds.  A complex $\C$ is \emph{semi-Eulerian} if $\chi(\lk (\F, \C)) = \chi(S^{d-\dim (\F)-1})$ for all faces $\F \in \C$. That is, the link of a face has the Euler characteristic of a sphere of appropriate dimension.  If, in addition, $\chi(\C) = \chi(S^{d-1})$, $\C$ is called an \emph{Eulerian complex}.  Semi-Eulerian complexes were called Eulerian manifolds in \cite{Klee}.  If we fix a field $k$ then a complex is a \emph{$k$-homology manifold} if, for all $x \in |\C|$, $\widetilde{H}_i(|\C|, |\C| - x; k) = 0$ when $i < d-1$ and equals either $k$ or 0 when $i=d-1$.  This is equivalent to saying that every non-empty face $\F$ of $\C$ has $k$-homology isomorphic to a sphere or ball of dimension $d-\dim(\F) -1$.  A $k$-homology manifold without boundary is therefore a semi-Eulerian complex.

Since the Euler characteristic is based on the Euler relation, we can ask if there are analogs of the Dehn-Sommerville relations that give more information about simplicial complexes.  In fact, these relations apply to a wide class of complexes, rather than differentiating between them.

>From Klee \cite{Klee} we know that a variant of the Dehn-Sommerville relations holds on all semi-Eulerian complexes:
\begin{theorem}[Klee]\label{Klee DS}
If $\C$ is a simplicial $(d-1)$-semi-Eulerian complex, then $\C$ satisfies the Dehn-Sommerville relations for $0 \leq k \leq d-2$, that is, 
\[\sum_{j=k}^{d-1} (-1)^j {j+1 \choose k+1} f_{j}(\C) = (-1)^{d-1}f_{k} (\C).\]\hfill $\Box$
\end{theorem}
These relations are sometimes denoted $E_k^d$ for $k=0, \ldots, d-2$ \cite{BB}. Taking $k=d-1$, the relation is simply an identity.  Taking $k=-1$ in the left hand side gives $\displaystyle{\sum_{i=-1}^{d-1} (-1)^j f_{j}(\C)}$, which equals $\chi(\C) - 1.$  If $\C$ is a simplicial polytope, the Euler relation tells us that $\chi(\C) = 1+(-1)^{d-1}$.  In the case of a complex, $\chi(\C)$ varies, so the $k=-1$ case is not included in the theorem.  If $\C$ is a $(d-1)$-dimensional semi-Eulerian complex, where $d$ is even, a linear combination of these relations shows that $\chi(\C) = 0$, so the complex is Eulerian. %SHOW!! Rewrite in h-form??

Let $\C$ and $\C'$ be two simplicial complexes with facets $\F$ and $\F'$, respectively.  Choose a bijection between the vertices of $\F$ and $\F'$.  The \emph{connected sum} of $\C$ and $\C'$, $C \# \C'$, is the complex constructed by identifying the vertices and corresponding faces of $\C$ and $\C'$ according to the chosen bijection and then removing the facet $\F=\F'$.  If both complexes are homology manifolds without boundary, then the connected sum is as well.  One special case of this construction is a \emph{stacked} polytope.  A complex is called a stacked polytope if it is a simplex or the connected sum of a simplex and a stacked polytope. This can also be thought of as iteratively taking pyramids over a facet.  For example, we can take the pyramid over one facet of a tetrahedron and get the bipyramid over a triangle.  The boundaries of these polytopes are called \emph{stacked spheres}.
%PICTURE?
Although any stacked polytope can be made while maintaining the convexity of the polytope, poorly choosing the apex for a pyramid may result in a non-convex set.  This does not change the combinatorics of the complex, so we will allow the construction to create non-convex sets, referring to the the boundary complex rather than the polytope at these times to avoid confusion.  

Another method to create new complexes is via \emph{handle addition}.  If $\F$ and $\F'$ are disjoint facets of $\C$ and a bijection between the two facets is chosen, then we can identify vertices and corresponding faces according to the bijection and remove $\F = \F'$.  As long as identified vertices are not neighbors of each other or both neighbors of the same vertex, the resulting complex will be a simplicial complex obtained by handle addition.  If the original complex is a homology manifold without boundary, the new complex is as well.  When we begin with a complex homeomorphic to $S^2$ and do a sequence of $g$ handle additions, the complex is a surface of genus $g$.  It is known a surface of genus $g$ has Euler characteristic $\chi(\C) = 2-2g$.

In Chapters 3 and 4, we will make more general complexes by starting with a set of polytopes rather than simplices.  We will consider the Euler characteristic and Gram-like relations on these complexes, considering complexes made by generalizations of the connected sum and handle addition constructions.  

%Walkup's \H_1
%CONJECTURES? Kalai & Kuhnel

%%Connections to angle-polynomial?

\chapter{Affine Spans of Angle Sums}

In this chapter we will consider the spaces spanned by the $\alpha$-vectors of simplices and the $\alpha$-$f$-vectors of simplicial polytopes and general polytopes.  We will construct families of polytopes whose $\alpha$-vectors and $\alpha$-$f$-vectors span the spaces defined by the Gram and Perles equations. In the first section, we will construct these polytopes and, in the second section, we will define the $\gamma$-vector, an analog to the $h$-vector, and consider the effect of the constructions on both the $h$-vectors and the $\gamma$-vectors  of polytopes. Then in the third section we use the $\gamma$-vector to show that the $\alpha$-vectors and $\alpha$-$f$-vectors of the constructed polytopes span the spaces defined by the relations on angle sums and face numbers mentioned in Chapter 1. That is, we show that the dimension of the affine span of the space of $\alpha{}$-vectors of simplices is $\left\lfloor \frac{d-1}{2}\right\rfloor$, the dimension of the affine span of $\alpha{}$-$f$-vectors of simplicial polytopes is $d-1$, and the dimension of the affine span of $\alpha{}$-$f$-vectors of general polytopes is $2d-3$.  
\section{Construction of Polytopes}

We will define two construction operations, the pyramid and prism operations, that create polytopes with varying angle sums. Each polytope will be constructed from a polytope of dimension one lower. These are similar to the pyramid and bipyramid constructions done by Bayer and Billera \cite{BB}, although, rather than bipyramids, we will build the dual polytope, prisms.  For a $(d-1)$-polytope $Q$, we will denote a $d$-pyramid over it as $PQ$ and the d-prism over it as $B^{*}Q$, following Bayer and Billera's notation for pyramids and bipyramids. However,  we will fix the geometry of the polytopes and not just the combinatorial structure.  

The \emph{prism} over the $(d-1)$-polytope $Q$, $B^{*}Q$, is $Q \times I$, 
where $I=[0,k]$ for some $k$. Then any $i$-face $\F$ of $B^{*}Q$ is either 
an $i$-face of one of $Q \times \{0\}$ or $Q \times \{k\}$, or, for some 
$(i-1)$-face $\G \subseteq Q$, $\F=\G \times I$, which is perpendicular to 
both $Q \times \{0\}$ and $Q \times \{k\}$. If $\F$ is a face of this 
latter type, then $\alpha(\F, B^{*}Q) = \alpha(\G, Q)$.  No angles change 
as the distance between the two copies of $Q$ varies, so the angle sums do 
not depend on $k$, but only on the prism construction.  Then, using the 
convention $\alpha_{-1}=0$, 
we have the following relationships on the $f$-vector and angle sums: 
\begin{equation}\label{alpha and B*} \begin{split} f_{0}(B^{*}Q) & = 
2f_{0}(Q), \\ f_{i}(B^{*}Q) & = 2f_{i}(Q) + f_{i-1}(Q) \quad \text{for } 1 
\leq i \leq d,\\ \alpha_{i}(B^{*}Q)& = \alpha_{i}(Q)+\alpha_{i-1}(Q) \quad 
\text{for } 0 \leq i \leq d. \end{split} \end{equation}

\singlespacing
\begin{figure}
\begin{center}
\includegraphics[width=5.0in]{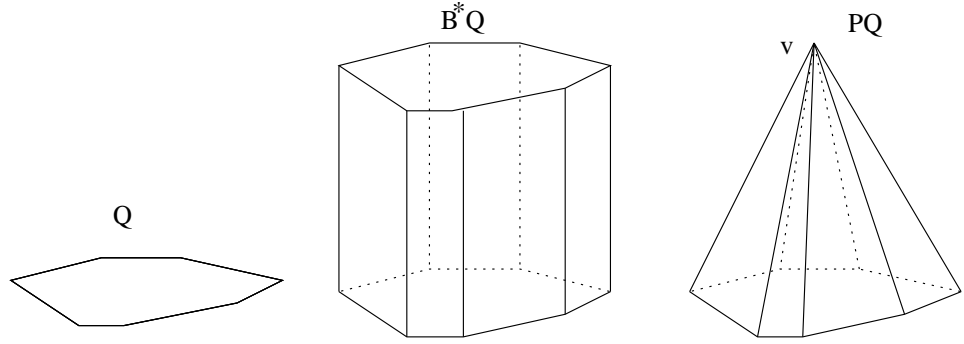}
\caption[The prism and pyramid over a polytope $Q$.]{The polytope $Q$; $B^{*}Q$, the prism over $Q$; and $PQ$, the pyramid over $Q$.}
\end{center}
\end{figure}
\normalspacing

Now we fix the geometry of $PQ$, the \emph{pyramid} over the polytope $Q$.  We start by placing a $(d-1)$-dimensional polytope $Q$ in the hyperplane $x_{d}=0$ in $\R^{d}$. We then place a vertex, $v$, along the line through the centroid of $Q$ and perpendicular to $Q$, so that it has $d$th coordinate $k>0$.  $PQ$ is then the convex hull of $v$ and $Q$, agreeing with our earlier definition of a pyramid.  An $i$-face of $PQ$ is either an $i$-face of $Q$ and therefore part of the base of the pyramid, or the convex hull of $v$ and an $(i-1)$-face of $Q$.  We will refer to the latter as \emph{sides}. The angles formed between the sides and faces in the base increase as $k$ does. For this reason, we will denote the pyramid by $P_{k}Q$ to specify the height of $v$ and fix the geometry of the construction.  Taking a pyramid over a point $d$ times results in a $d$-simplex. Therefore, we will denote d-simplices as $P^d$, assuming a starting polytope of a point when one is not explicitly given.  For any $k$, the pyramid operation has the following effect on the $f$-vector \cite{BB}:
\begin{equation}\label{f and P}
\begin{split}
f_{i}(PQ) & = f_{i}(Q) + f_{i-1}(Q) \quad \text{for } 0 \leq i \leq d-1,\\
f_{d}(PQ) & = 1.
\end{split}
\end{equation}

We note two limiting cases of the pyramid operation: the case as $k$ tends toward $0$ and the case where $k$ tends toward infinity.  We will denote these constructions by $P_{0}Q$ and $P_{\infty}Q$, respectively. Although neither is actually a $d$-pyramid, since $P_{0}Q$ is $(d-1)$-dimensional and $P_{\infty}Q$ is not bounded, one can easily find the limits of the angle sums as $k$ tends to $0$ or infinity, and we will define these values as the angle sums for $P_{0}Q$ and $P_{\infty}Q$.  Since the values of the angle sums vary continuously as $k$ does, we can find pyramids with angle sums that are arbitrarily close to those of $P_{0}Q$ and $P_{\infty}Q$.We can picture $P_{0}Q$ as a `flat' pyramid, with two copies of $Q$ glued together, one of which has an extra vertex joined to each proper face.  $P_{\infty}Q$ can be pictured as a prism with infinite height or no top.

For $P_{0}Q$, all angles made between the base and sides tend to $0$, so any interior angles at proper faces of
 the base are $0$, and the interior angle at the base and at faces including the apex $v$ are all $\frac{1}{2}$.
 Therefore, all the angles sums are dependent on the $f$-vector of the 
base $Q$. Then, using the conventions $f_{-1}(Q) = 1$ and $\alpha_{-1}(Q) 
= 0$, we have \begin{equation}\label{alpha P0} \begin{split} 
\alpha_i(P_{0}Q) & = \frac{1}{2} f_{i-1}(Q) \quad \text{for }0 \leq i \leq 
d-2,\\ \alpha_{d-1}(P_{0}Q) & = \frac{1}{2}f_{d-2}(Q) + 
\frac{1}{2}. \end{split} \end{equation}

For $P_{\infty}Q$, angles between the sides and base tend to right angles, 
so for any face $\G \subseteq Q$, $\alpha(\G, P_{\infty}Q)= \frac{1}{2} 
\alpha(\G, Q)$.  For faces $\F \subsetneq P_{\infty}Q$ that are the convex 
hull of a face $\G \subseteq Q$ and $v$, the interior angle at $F$ is the 
same as it was at $\G$, that is, $\alpha(\F, P_{\infty}Q) = \alpha(\G, 
Q)$. Using the convention that $\alpha(\emptyset, Q) = 0$, this also 
applies to $v$ itself: $\alpha(v, P_{\infty}Q) = 0$. Therefore: 
\begin{equation}\label{alpha and Pinfty} \begin{split} 
\alpha_{i}(P_{\infty}Q) & = \frac{1}{2} \alpha_{i}(Q)+\alpha_{i-1}(Q), 
\quad \text{for }0 \leq i \leq d-1 \\ \alpha_{d}(P_{\infty}Q) 
&= 1. \end{split} \end{equation}

We will sometimes want to iterate these constructions; we will write $C^{k}Q$ when we wish to apply a construction $C$ $k$ times in succession to $Q$. Taking a pyramid over a point $d$ times results in a $d$-simplex. Therefore, we will denote $d$-simplices as $P^d$, assuming a starting polytope of a point when one is not explicitly given.
\section{The $\gamma$-vector}

In analogy to the $h$-vector, we define the \emph{$\gamma$-vector} as  \[\gamma(\P) = (\gamma_0(\P), \gamma_1(\P), \ldots, \gamma_d(\P)),\] where
\begin{equation}\label{gamma-vector}
\gamma_i(\P) = \sum_{j=0}^{i}(-1)^{i-j}{d-j \choose d-i}\alpha_{j-1}(\P).
\end{equation}
We note that $\gamma_0(\P) = 0$, $\gamma_1(\P) = \alpha_0(\P)$, and $\gamma_d(\P) = 1$ for all polytopes $\P$.
For convenience in considering the $\gamma$-vector of polytopes of increasing dimension made by the pyramid and bipyramid constructions,  we will define $\gamma_i(\P)$ for $i$ beyond the $\gamma$-vector.  We define %$h_i(\P) = 0$  and 
$\gamma_i(\P) = 0$ for $i < 0$  and $\gamma_i(\P) = 1$ for $i>\dim (\P)$ for all $\P$.

Kleinschmidt and Smilansky \cite{KS} defined a vector that agrees with the $\gamma$-vector on spherical simplices, calling the entries $\sigma_i(\D)$.  In this case, the sphere was decomposed into regions that were all simplices by the great spheres that defined $\D$, and $\sigma_i(\D)$ measured the area of all the regions that were reached from $\D$ by crossing $i$ great spheres.  We have chosen a different name for our vector to avoid confusion in the definition. 
 
The matrix which transforms $\left(\alpha_{-1}(\P), \alpha_{0}(\P), \ldots \alpha_{d-1}(\P)\right)$ to $\gamma(\P)$ is lower triangular with entries of 1 along the diagonal.  Therefore this transformation is invertible.  We define the $\gamma$-$h$-vector,
$(\gamma_0(Q), \ldots, \gamma_d(Q)| h_0(Q) \ldots, h_d(Q))$, and note that it is an invertible linear transformation of the $\alpha$-$f$-vector.

As with the $h$-vector formulation of the Dehn-Sommerville relations, we can rewrite the Perles relations in terms of the $\gamma$-vector:
\begin{theorem}\label{h-Perles}
For a simplicial $d$-polytope $\P$, 
\[\gamma_i(\P) + \gamma_{d-i}(\P) = h_i(\P) \quad \text{for } 0 \leq i \leq d.\]
\end{theorem}
\begin{proof}
The proof follows the one given for Corollary 2.2 in \cite{BB}. 
As in the transformation of the Dehn-Sommerville relations on the $f$-vector to their $h$-vector form, we take the linear combination
\[\sum_{i=0}^{r}(-1)^i{d-i \choose d-r}S_{i-1}^d,\]
where $S_k^{d}$ is
\[\sum_{j=k}^{d-1} (-1)^{j} {j+1\choose k+1} \alpha_{j}(\P) =
(-1)^{d}(\alpha_{k}(\P) - f_{k}(\P)),\]
the $k$th Perles relation on simplicial $d$-polytopes.  On the right hand side, the sum becomes
\begin{equation*}
\begin{split}
\sum_{i=0}^r(-1)^i&{d-i \choose d-r}(-1)^{d}(\alpha_{i-1}(\P) - f_{i-1}(\P))\\
& = (-1)^d\left[\sum_{i=0}^r(-1)^i{d-i \choose d-r}\alpha_{i-1}(\P) - \sum_{i=0}^r(-1)^i{d-i \choose d-r}f_{i-1}(\P) \right]\\
& = (-1)^{d}\left((-1)^{r}\gamma_r(\P) - (-1)^{r}h_r(\P)\right)\\
& = (-1)^{d-r}\left(\gamma_r(\P) - h_r(\P)\right).
\end{split}
\end{equation*}
Regarding the left hand side, we see that 
\begin{equation*}
\begin{split}
\sum_{i=0}^r(-1)^{i}{d-i \choose d-r}&\sum_{j=i-1}^{d-1} (-1)^{j} {j+1\choose i} \alpha_{j}(\P)  \\
&= \sum_{i=0}^r(-1)^{i}{d-i \choose d-r}\sum_{j=i}^{d} (-1)^{j-1} {j\choose i} \alpha_{j-1}(\P)\\
&=\sum_{j=0}^{d} (-1)^{j-1}\alpha_{j-1}(\P)\sum_{i=0}^{j}(-1)^{i}{d-i \choose d-r} {j\choose i}.
\end{split}
\end{equation*}
Then we apply the identity 
\[\sum_{s=0}^n (-1)^{s}{s+m \choose t}{n \choose s} = (-1)^{n}{m \choose t-n}\] 
to simplify the interior sum
\begin{equation*}
\begin{split}
\sum_{i=0}^{j}(-1)^{i}{d-i \choose d-r}{j\choose i}  & = \sum_{s=0}^{j}(-1)^{j-s}{d-j+s \choose d-r} {j \choose s}\\
&= {d-j \choose d-r-j} = {d-j \choose r}.
\end{split}
\end{equation*}
Therefore the left hand side simplifies to 
\begin{equation*}
\begin{split}
\sum_{i=0}^r(-1)^{i}{d-i \choose d-r}\sum_{j=i-1}^{d-1} (-1)^{j} {j+1\choose i} \alpha_{j}(\P) & = \sum_{j=0}^{d} (-1)^{j-1}\alpha_{j-1}(\P){d-j \choose r}\\
& = \sum_{j=0}^{d-r} (-1)^{j-1}{d-j \choose r}\alpha_{j-1}(\P)\\
& = (-1)^{d-r+1}\gamma_{d-r}(\P).
\end{split}
\end{equation*}
Putting these results together we see that
\[(-1)^{d-r+1}\gamma_{d-r}(\P) = (-1)^{d-r}\left(\gamma_r(\P) - h_r(\P)\right)\]
or
\[\gamma_{d-r}(\P) =  h_r(\P)-\gamma_r(\P).\]
\end{proof}

In preparation for using the constructions to create affinely independent $\gamma$-vectors, we consider the effect of the pyramid and prism constructions on the $\gamma$-vector.  In the following propositions, we consider the $h$-vector entries strictly as a linear combination of the $f$-vector entries and do not assume that the polytope is simplicial.  
\begin{proposition}\label{gamma and P0}
If $Q$ is a $(d-1)$-polytope, 
\[h(PQ) =(h(Q), 1),  \]
or, equivalently, \[h_i(PQ) = h_i(Q) \quad \text{for }0 \leq i \leq d-1 \text{ and } h_{d}(PQ) = 1.\]
Also,
\[\gamma_i(P_{0}Q) = \frac{1}{2}h_{i-1}(Q) \quad \text{for }0 \leq i \leq d-1 \text{ and } \gamma_{d}(P_{0}Q) = 1.\]
\end{proposition}
\begin{proof}
The first relation is proved in Proposition 3.1 of \cite{BB}.

>From \eqref{alpha P0}, we know that if Q is a $(d-1)$-polytope, $\alpha_j(P_{0}Q) = \frac{1}{2}f_{j-1}(Q)$ for $0 \leq j \leq d-2$ and $\alpha_{-1}(P_{0}Q) = 0$.  Therefore, for $0 \leq i \leq d-1$,
\begin{equation*}
\begin{split}
\gamma_i(P_{0}Q) &= \sum_{j=0}^{i}(-1)^{i-j}{d-j \choose d-i}\alpha_{j-1}(P_{0}Q)\\%alpha_-1 = 0
& = \frac{1}{2}\sum_{j=1}^{i}(-1)^{i-j}{d-j \choose d-i}f_{j-2}(Q)\\
& = \frac{1}{2}\sum_{j=0}^{i-1}(-1)^{i-j-1}{d-j-1 \choose d-i}f_{j-1}(Q)\\
& = \frac{1}{2}\sum_{j=0}^{i-1}(-1)^{(i-1) - j} {(d-1) - j \choose (d-1) - (i-1)} f_{j-1}(Q)\\
& = \frac{1}{2}h_{i-1}(Q).
\end{split}
\end{equation*}
\end{proof}
The polytope $PQ$ is simplicial only if $Q$ is a simplex since $Q$ is a facet of $PQ$.  Since the $h$-vector of a line segment is $(1,1)$, induction using the proposition shows that that the $h$-vector of a $d$-simplex $\D$ is $(1,1,\ldots,1)$, so that $\gamma_i(\D) + \gamma_{d-i}(\D) = 1$ for $0 \leq i \leq d$ by Theorem \ref{h-Perles}.  In particular, the proposition shows that $\gamma(P_0^{d-1}P) = (0,\frac{1}{2}, \ldots \frac{1}{2}, 1)$.

\begin{proposition}\label{gamma and Pinfinity}
If $Q$ is a $(d-1)$-polytope, 
\[\gamma(P_{\infty}Q) = \frac{1}{2}\left[(0, \gamma(Q))+ (\gamma(Q),1)\right], \]
or, equivalently,
\[\gamma_i(P_{\infty}Q) = \frac{1}{2}\gamma_i(Q) + \frac{1}{2}\gamma_{i-1}(Q) \quad \text{for } 0 \leq i \leq d.\]
More generally, and using the extended $\gamma$-entries, 
\[\gamma_i((P_{\infty})^{k}Q) = \frac{1}{2^k}\sum_{j=0}^{k}{k \choose j}\gamma_{i-j}(Q) \quad \text{for } 0 \leq i \leq d.\]
\end{proposition}
\begin{proof}
By \eqref{alpha and Pinfty}, $\alpha_{i}(P_{\infty}Q) =  \frac{1}{2} 
\alpha_{i}(Q)+\alpha_{i-1}(Q)$ for $0 \leq i \leq d-1$.
Then we calculate
\begin{equation*}
\begin{split}
\gamma_i(P_{\infty}Q) &= \sum_{j=0}^{i}(-1)^{i-j}{d-j \choose d-i}\alpha_{j-1}(P_{\infty}Q)\\ 
& = \sum_{j=0}^{i}(-1)^{i-j}{d-j \choose d-i}\left(\frac{1}{2}\alpha_{j-1}(Q)+\alpha_{j-2}(Q)\right)\\%%alpha_{-1} = 0
& = \frac{1}{2}\sum_{j=0}^{i}(-1)^{i-j}{d-j \choose d-i}\alpha_{j-1}(Q) \\
& \qquad \qquad \qquad + \sum_{j=0}^{i-1}(-1)^{(i-1)-j}{(d-1)-j \choose (d-1)-(i-1)}\alpha_{j-1}(Q)\\
& = \frac{1}{2}\left[\sum_{j=0}^{i}(-1)^{i-j}{d-j-1 \choose d-i}\alpha_{j-1}(Q) \right.\\
& \qquad \qquad \qquad \left.+ \sum_{j=0}^{i}(-1)^{i-j}{d-j-1 \choose d-i-1}\alpha_{j-1}(Q)\right] + \gamma_{i-1}(Q)\\ %{d-i-1 \choose d-i}=0
& = -\frac{1}{2}\gamma_{i-1}(Q)+ \frac{1}{2}\gamma_{i}(Q) + \gamma_{i-1}(Q)\\
& = \frac{1}{2}\gamma_{i}(Q) + \frac{1}{2}\gamma_{i-1}(Q).
\end{split}
\end{equation*}

As we iterate the $P_{\infty}$ construction, we see that 
\[\gamma_i((P_{\infty})^2Q) = \frac{1}{4}\gamma_{i}(Q) + 
\frac{1}{2}\gamma_{i-1}(Q) + \frac{1}{4}\gamma_{i-2}(Q)\] and 
\[\gamma_i((P_{\infty})^2Q) = \frac{1}{8}\gamma_{i}(Q) + 
\frac{3}{8}\gamma_{i-1}(Q) + \frac{3}{8}\gamma_{i-2}(Q) + 
\frac{1}{8}\gamma_{i-3}(Q).\] In each iteration we see that 
$\gamma_i((P_{\infty})^kQ)$ is a linear combination of $\gamma_j(Q)$ for 
$i-k \leq j \leq i$ and the coefficient of $\gamma_j(Q)$ is half the 
sum of the coefficients of $\gamma_j(Q)$ and $\gamma_{j-1}(Q)$ in the 
linear combination for $\gamma_i((P_{\infty})^{k-1}Q)$.  Therefore, the 
coefficient of $\gamma_j(Q)$ in $\gamma_i((P_{\infty})^kQ)$ is 
$\frac{1}{2^k}{k \choose j}$. \end{proof}

\begin{proposition}\label{gamma and B*}
If $Q$ is a $(d-1)$-polytope, 
\[\gamma(B^{*}Q) = (\gamma(Q), 1),\] or, equivalently, \[\gamma_i(B^{*}Q) = \gamma_i(Q)  \quad \text{for all } i.\]
\end{proposition}
\begin{proof}
By \eqref{alpha and B*}, $\alpha_{i}(B^{*}Q) =  
\alpha_{i}(Q)+\alpha_{i-1}(Q)$ for $0 \leq i \leq d$ .

Then we can calculate:
\begin{equation*}
\begin{split}
\gamma_i(B^{*}Q) &= \sum_{j=0}^{i}(-1)^{i-j}{d-j \choose d-i}\alpha_{j-1}(B^{*}Q)\\ 
%& = \sum_{j=0}^{i}(-1)^{i-j}{d-j \choose d-i}\left(\alpha_{j-1}(Q)+\alpha_{j-2}(Q)\right)\\%%alpha_{-1} = 0
& = \sum_{j=0}^{i}(-1)^{i-j}{d-j \choose d-i}\alpha_{j-1}(Q) + \sum_{j=0}^{i-1}(-1)^{i-j-1}{d-j-1 \choose d-i}\alpha_{j-1}(Q)\\
& = \sum_{j=0}^{i}(-1)^{i-j}\left({d-j \choose d-i} - {d-j-1 \choose d-i}\right)\alpha_{j-1}(Q)\\
& = \sum_{j=0}^{i}(-1)^{i-j}{(d-1)-j \choose (d-1)-i}\alpha_{j-1}(Q)\\
& = \gamma_i(Q).
\end{split}
\end{equation*}
This calculation applies for $i = 0, 
\ldots, d-1$ and, since $\gamma_d(\P) = 1$ for all $d$-polytopes $\P$ and we have defined $\gamma_d(Q) = 1$ in the extended $\gamma$-vector, $\gamma_d(B^{*}Q) = \gamma_d(Q)$.  Therefore we can write $\gamma_i(B^{*}Q) = \gamma_i(Q)$ for all $i$ or $\gamma(B^{*}Q) = (\gamma(Q), 1)$.
\end{proof}

\section{Spans of $\alpha$ and $\alpha$-$f$-vectors}
Using the prism and pyramid constructions, we can now build families of polytopes with affinely independent
$\alpha$-vectors or $\alpha$-$f$-vectors.  We will use these families to span the spaces of $\alpha$-vectors and $\alpha$-$f$-vectors
defined by the Gram and Perles relations.

In order to prove results about the affine span of $\alpha$-vectors and $\alpha$-$f$-vectors, we will want to work with the constructions $B^{*}$, $P_0$ and $P_{\infty}$.  However, as mentioned before, $P_0$ and $P_{\infty}$ are limiting cases of the pyramid construction and do not create $d$-polytopes.  Therefore, we need the following lemma to tell us that, when we create $\alpha$-vectors using these constructions, we can find a set of $d$-polytopes that maintain the independence properties of the $\alpha$-vectors and the $\alpha$-$f$-vectors of the polytopes. 
 
\begin{lemma}\label{backing off limiting}
Let $\e>0$ be given, and let $Q_i$, $i=0,\ldots, k$, be a set of $d$-'polytopes' where each $Q_i$ is a $d$-polytope or has form  $(P_0)^{k}Q$ or $(P_{\infty})^{k}(P_0)^{l}Q$ for some nonnegative integers $k$ and $l$ and a polytope $Q$.  

Then there is a set of $d$-polytopes $Q'_i$, $i=0,\ldots, k$, where $Q'_i$ and $Q_i$ have the same $f$-vector and $\left|\alpha_j(Q'_i)- \alpha_j(Q_i)\right|< \e$ for all $i$ and $j$.  Further, if the $\alpha$-vectors or $\alpha$-$f$-vectors of the $Q_i$ are affinely independent, $\e$ can be chosen so that the $\alpha$-vectors or $\alpha$-$f$-vectors of the $Q'_i$ are also affinely independent.
\end{lemma}
\begin{proof}
For each $i$ we will define constants $M_i$ and $\delta_i$.  Suppose $Q_i$ has form $P_{\infty}Q$.  Since the angle sums are continuous, we can choose $M_{i}^{j}$ so that $\alpha_j(P_{N}Q)$ is within $\e$ of $\alpha_j(P_{\infty}Q)$, for any $N \geq M_i^{j}$.  Let $M_{i}=\max_{j}M^{j}_{i}$.  Then for any $N \geq M_{i}$, $\left|\alpha_{j}(P_{N}Q) - \alpha_{j}(P_{\infty}Q)\right| < \e$ for all $j$. 

If $Q_i = (P_{\infty})^kQ$, we can iterate this process with difference $\e/k$.  Starting with $M_{i_0} = 1$, we iteratively choose $M_{i_{m}}\geq M_{i_{m-1}}$ by the same process as above so that  $\left|\alpha_{j}((P_{\infty})^{k-m}(P_{N})^m Q) - \alpha_{j}((P_{\infty})^{k-m+1}(P_{N})^{m-1} Q)\right| < \e/k$ for $N \geq M_{i_{m}}$ and all $j$.  Then we let $M_{i} = M_{i_k}$ so that if $N \geq M_i$,  $\left|\alpha_j((P_{N})^{k} Q) - \alpha_j((P_{\infty})^{k} Q)\right|<\e$.   

If $Q_i$ has form $P_{0}Q$, an analogous argument finds $\delta_i$ such that for all $\delta\leq \delta_i$ and $j$, $\left|\alpha_j(P_{\delta} Q) - \alpha_j(P_{0} Q)\right|<\e$.  Since the $\alpha$-vector of $P_{0}Q$ is entirely determined by the combinatorics of $Q$, this one step is also sufficient to choose $\delta_i$ for $Q_i = (P_{0})^{k}Q$.
\newpage
Suppose $Q_i$ has form $(P_{\infty})^{k}(P_0)^{l}Q$ with $l \geq 1$.  First choose $\delta_i$ so that for all $\delta\leq \delta_i$, $\left|\alpha_j((P_{\delta})^{l} Q) - \alpha_j((P_{0})^{l} Q)\right|<\frac{\e}{2^{k+1}}$ for all $j$.  Then by \eqref{alpha and Pinfty} the $P_{\infty}$ construction will less than double any differences in angle sum values, so \[\left|\alpha_j((P_{\infty})^k(P_{\delta})^{l} Q) - \alpha_j((P_{\infty})^k(P_{0})^{l} Q)\right|<\e/2\]
for $\delta\leq \delta_i$ and all $i$ and $j$. Next choose $M_i$ so that  
\[\left|\alpha_{j}((P_{N})^k(P_{\delta})^{l}Q) - \alpha_{j}((P_{\infty})^k(P_{\delta})^{l}Q)\right| < \e/2\] for $N \geq M_{i}$ and all $i$ and $j$.  Then \[\left|\alpha_{j}((P_{N})^k(P_{\delta})^{l}Q) - \alpha_{j}((P_{\infty})^k(P_{0})^{l}Q)\right| < \e\]
 for $N \geq M_{i}$, $\delta\leq \delta_i$,  and all $j$.  

Now we choose
\[Q'_i = \begin{cases}
(P_{\delta_i})^{k}Q & \text{ if } Q_i = (P_{0})^{k}Q\\
(P_{M_i})^k(P_{\delta_i})^{l}Q & \text{ if } Q_i = (P_{\infty})^{k}(P_{0})^{l}Q\\
Q_i & \text{ if } Q \text{ is a }d\text{-polytope}.
\end{cases}\]
The $f$-vectors of $Q_i$ and $Q'_i$ are the same since they are pyramids of the same degree over the same polytope and $\left|\alpha_j(Q'_i)- \alpha_k(Q_i)\right|< \e$ for all $i$ and $j$. 

Since affine independence is an open condition, if the $Q_i$ have affinely independent $\alpha$-vectors or $\alpha$-$f$-vectors, we can choose $\e$ small enough that the $Q'_i$ given above have affinely independent $\alpha$-vectors or $\alpha$-$f$-vectors, respectively.  
\end{proof}

The $P_{\infty}$ construction will be particularly useful for increasing the dimension of a set of polytopes and maintaining the affine independence of their $\alpha$-vectors, as shown in the next lemma.
\begin{lemma}\label{P_{infty} independence}
If a set of $(d-1)$-polytopes $Q_i$, $i = 0,\ldots, k$, has affinely independent $\alpha$-vectors, then the set of $d$-polytopes $P_{\infty}Q_i$, $i = 0,\ldots, k$, also has affinely independent $\alpha$-vectors.
\end{lemma}
\begin{proof}
Since the last entry of $\alpha(Q)$ is 1 for every polytope, the affine independence of a set of $\alpha$-vectors is equivalent to their linear independence.  Also, the $\gamma$-vector is an invertible linear transformation of the $\alpha$-vector, so the linear independence of a set of $\alpha$-vectors is equivalent to the linear independence of the corresponding set of $\gamma$-vectors.  We will work with the $\gamma$-vectors for ease of computation.

Based on Proposition \ref{gamma and Pinfinity}, we can write
\begin{equation}\label{Pinfty trans on gamma}
\gamma(P_{\infty}Q) =\textbf{A}\begin{bmatrix} \gamma(Q)\\1\end{bmatrix}
\end{equation}
where
\begin{equation*}
\textbf{A} = \frac{1}{2}\left[\begin{array}{ccccc}
 1& & & &  \\
 1&1& & &  \\
  &1&1& &   \\
 &&\ddots&\ddots& \\ 
  & & &1&1\\
\end{array}\right],
\end{equation*}
a $(d+1)\times(d+1)$ matrix where all other entries are 0.  Clearly $\textbf{A}$ is invertible.  Its inverse is
\begin{equation}\label{Ainverse}
\textbf{A}^{-1} = \left[\begin{array}{rrrrr}
 2 & & & &\\ 
 -2 & 2 &  & &   \\
 2 & -2 & 2 &&\\
\vdots& & &\ddots\\ 
&\ldots&2&-2&2\\
\end{array}\right],
\end{equation}
where all entries on and below the diagonal alternate between $2$ and $-2$ and entries above the diagonal are 0. 
Since the matrix is invertible, we know that the $\gamma$-vectors of $P_{\infty}Q_{i}$ for $i=1,\ldots,k$ are linearly independent since the $\gamma$-vectors of $Q_{i}$ for $i=1,\ldots,k$ were. 
\end{proof}

\begin{theorem}\label{Angle sums of simplices} The affine span of the $\alpha$-vectors of $d$-simplices has dimension
$\left\lfloor\frac{d-1}{2}\right\rfloor$.
\end{theorem}
\begin{proof} Let $A$ be the affine space spanned by the $\alpha$-vectors of $d$-simplices.  We want to show that $\dim (A) = \left\lfloor\frac{d-1}{2}\right\rfloor$.  As in Lemma \ref{P_{infty} independence}, we will show the equivalent fact that the space of $\gamma$-vectors has dimension $\left\lfloor\frac{d-1}{2}\right\rfloor$ in order to simplify calculation.

We first prove that $\left\lfloor\frac{d-1}{2}\right\rfloor$ is an upper bound on the dimension of $A$ and then construct a family of polytopes to show this bound is achieved. If $\D$ is a $d$-simplex, we know that $f_{i}(\D) = {d+1 \choose i+1}$. So by Theorem \ref{h-Perles} and Lemma \ref{gamma and P0}, the Perles equations on the $\gamma$-vector become:
\begin{equation*}
\gamma_k(\D) + \gamma_{d-k}(\D) = 1,
\end{equation*}
which we will call $S^d_{k}(\D)$.  Then the relations 
$S^d_{0}(\D),S^d_{1}(\D), \ldots, S^d_{\left\lfloor\frac{d}{2}\right\rfloor}(\D)$ are clearly independent.  Since the $\gamma$-vector is $(d+1)$-dimensional and all $\alpha$-vectors lie in the plane $\gamma_{d} = 1$, we get:
\begin{equation*}
\dim(A) \leq d+1-\left(\left\lfloor\frac{d}{2}\right\rfloor+1\right)- 1 =\left\lfloor\frac{d-1}{2}\right\rfloor.
\end{equation*}

We will prove that $\dim(A) \geq \left\lfloor\frac{d-1}{2}\right\rfloor$ by constructing a set of $\left\lfloor\frac{d-1}{2}\right\rfloor + 1$ simplices whose $\gamma$-vectors are affinely independent. 
  The proof will proceed by induction on $d$, first using the limiting constructions $P_{\infty}$ and $P_0$ and then finding $d$-polytopes whose angle sums are arbitrarily close to these polytopes.  For $d=1$ and $d=2$, a line segment and a triangle (denoted $P$ and $P^{2}$, respectively) provide the one element needed for the basis. 

Suppose we have a set of $\left\lfloor \frac{d-3}{2} \right\rfloor +1 = \left\lfloor \frac{d-1}{2} \right\rfloor$ simplices $Q_i$ for $i = 1,\ldots,\left\lfloor \frac{d-1}{2} \right\rfloor$ in dimension $d-2$, $d \geq 3$, with linearly independent $\gamma$-vectors. We claim that the vectors $\gamma\Big(P_0^{d-1}P\Big), \gamma\Big(P_{\infty}^2Q_1\Big), \ldots, \gamma\Big(P_{\infty}^2Q_{\left\lfloor \frac{d-1}{2} \right\rfloor}\Big)$ are linearly independent. If this claim is true, then we have $\left\lfloor \frac{d-1}{2} \right\rfloor+1$ simplices with linearly independent $\gamma$-vectors in dimension $d$ as needed for the theorem. Since the theorem is true for $d=1$ and $d=2$, the proof of this claim will finish the proof of the theorem.  

We know the vectors $\gamma\Big(P_{\infty}^2Q_1\Big), \ldots, \gamma\Big(P_{\infty}^2Q_{\left\lfloor \frac{d-1}{2} \right\rfloor}\Big)$ are linearly independent by applying Lemma \ref{P_{infty} independence} twice.  We will show that adding the vector $\gamma(P_0^{d+1}P)$ increases the linear span by showing that  the linear span of the inverse images of the $\left\lfloor \frac{d-1}{2} \right\rfloor+1$ $\gamma$-vectors of $d$-polytopes under the $P_{\infty}^2$ transformation must be greater than the linear span of the $\gamma(Q_i)$. More specifically we will show that  $\left(\textbf{A}^{-1}\right)^2\left(\gamma(P_0^{d-1}P)\right)$ must be outside the linear span of $\gamma(Q_i)$ for $i=0, \ldots, \left\lfloor \frac{d-1}{2} \right\rfloor$, where $\textbf{A}^{-1}$ is as in \eqref{Ainverse}.

Now by Proposition \ref{gamma and P0}, $\gamma(P_0^{d-1}P) = \left(0, \frac{1}{2},\frac{1}{2}, \ldots, \frac{1}{2}, 1\right)$.  Then if $d\geq3$, 
\[\left(\textbf{A}^{-1}\right)^2\left(\gamma(P_0^{d-1}P)\right) = \begin{bmatrix} 0\\2\\ \vdots\end{bmatrix},\]
where the entries alternate in sign. Therefore, the last three entries do not have the same value. But each of the vectors \[\left(\textbf{A}^{-1}\right)^2\left(\gamma\left(P_{\infty}^2Q_i\right)\right) = \begin{bmatrix} \gamma(Q_i)\\1\\1\end{bmatrix}\]
has last three entries 1, as do the extended $\gamma$-vectors of all $(d-1)$-polytopes, and therefore all the vectors in the span of the $\gamma(Q_i)$ must have the same value on last three entries.  Therefore $\left(\textbf{A}^{-1}\right)^2\left(\gamma(P_0^{d-1}P)\right)$ is outside the linear span of the $\gamma(Q_i)$ for $i=1, \ldots, \left\lfloor \frac{d-1}{2} \right\rfloor)$ and
\[\dim\left(\spanr\{\gamma(P_0^{d-1}P), \gamma(P_{\infty}^2Q_i)\}\right) = \dim\left(\spanr\{\gamma(P_{\infty}^2Q_i)\}\right) + 1.\]

Therefore, the following sets of simplices inductively constructed above for dimension $d$ have affinely independent $\alpha$-vectors:
\[P_{\infty}^{d-1}P, P_{\infty}^{d-3}P_{0}^{2}P, P_{\infty}^{d-5}P_{0}^{4}P, \ldots, P_{0}^{d-1}P \qquad \text{if }d\text{ is odd}\]
and 
\[P_{\infty}^{d-2}P^2, P_{\infty}^{d-4}P_{0}^{2}P^2, P_{\infty}^{d-6}P_{0}^{4}P^2, \ldots, P_{0}^{d-2}P^2 \qquad \text{if }d\text{ is even}.\]
(The latter uses that $\alpha(P_{0}P) = \alpha(P^2)$ since there is only one $\alpha$-vector for triangles.)

Then by Lemma \ref{backing off limiting}, we know that we have a set of $\left\lfloor \frac{d-1}{2} \right\rfloor + 1$ $d$-simplices with affinely independent $\alpha$-vectors. 
\end{proof}

In $\R^{3}$, this theorem tells us that the $\alpha$-vectors of two polytopes span the affine space of $\alpha$-vectors of simplices: those for $P_{\delta}P^{2}$ and for $P_{N}P^{2}$, a very short tetrahedron and a very tall one.  Therefore, even though there is only one $f$-vector, we have a one-parameter family of $\alpha$-vectors of tetrahedra and this one parameter can describe the height of the simplex.  Similar tetrahedra have the same $\alpha$-vector and changing the base of the tetrahedron reparametrizes the family of $\alpha$-vectors rather than giving any new descriptions.
%Question: as far as the alphavectors are concerned do we get a segment or a line?
 
Using what is known about the affine span of the $f$-vectors of simplicial polytopes together with the results of the preceding theorem, we can determine the affine span of the $\alpha$-$f$-vectors of simplicial polytopes. It is appropriate to consider this vector rather than the $\alpha$-vector in describing the angles of simplicial polytopes since the Perles relations refer to face numbers as well as angle sums.
\begin{theorem}\label{Angle sums of simplicial polytopes} The affine span of the $\alpha$-$f$-vectors of simplicial $d$-polytopes has dimension $d-1$.  The space is spanned by $\left\lfloor\frac{d+1}{2}\right\rfloor$ simplices, as in Theorem \ref{Angle sums of simplices}, and $\left\lfloor\frac{d}{2}\right\rfloor$ non-simplices which are combinatorially independent simplicial polytopes. 
\end{theorem}
\begin{proof} Let $A_S$ be the affine space spanned by the $\alpha$-$f$-vectors of simplicial polytopes.  Since the $h$-vector is computed independently of the $\gamma$-vector, and each is an invertible linear transformation, the affine independence of the $\alpha$-$f$-vectors of  is equivalent to the affine independence of the corresponding $\gamma$-$h$-vectors.  

By the same argument as in Theorem \ref{Angle sums of simplices}, there are $\lfloor \frac{d}{2}\rfloor+1$ Perles relations that are independent with regard to angle sums, and in this case each includes a different element of the $h$-vector.  The other relation on the angle sums is that $\gamma_d(P)=1$ for all polytopes.  Similarly, the Dehn-Sommerville relations on the $h$-vector show that there are $\lfloor \frac{d+1}{2}\rfloor$ independent Dehn-Sommerville relations, since if $d$ is even, $h_{\lfloor \frac{d}{2}\rfloor} = h_{\lfloor \frac{d}{2}\rfloor}$ does not provide a new relation. We also have the relation $f_d(P)=1$.  Since these relations include no angle sums, they are in turn independent of the Perles relations.  Therefore, 
\begin{equation*}
\text{dim}(A_S) \leq 2d+2-\left(\left\lfloor\frac{d}{2}\right\rfloor+1\right)- 1 - 
\left(\left\lfloor \frac{d+1}{2}\right\rfloor\right) - 1=d-1.
\end{equation*}

The affine span of the $h$-vectors of simplicial $d$-polytopes is 
of dimension $\left\lfloor\frac{d}{2}\right\rfloor$.  In Bayer and Billera 
\cite{BB}, a set of $\left\lfloor\frac{d}{2}\right\rfloor + 1$ simplicial 
polytopes with affinely independent $h$-vectors is given, spanning the 
space defined by the Dehn-Sommerville equations.  This basis includes one 
simplex.  We can combine the $\left\lfloor\frac{d}{2}\right\rfloor$ 
non-simplices of this basis with the 
$\left\lfloor\frac{d+1}{2}\right\rfloor$ simplices given in Theorem 
\ref{Angle sums of simplices}.

If the $\gamma$-$h$-vectors of this set, $Q_i$ for $i=1,\ldots, d$ are affinely dependent, then there exist $\lambda_i$, $i=1,\ldots,d$ such that
\[\sum_{i=1}^{d} \lambda_i\left(\gamma\text{-}h(Q_i)\right)= 0,\]
where $\displaystyle{\sum_{i=0}^{d}\lambda_i = 0}.$
But if the $Q_i$, $i = 1,\ldots,\left\lfloor\frac{d+1}{2}\right\rfloor$, are simplices and we consider just the $h$-vector entries, we can rewrite the above as 
\[0=\sum_{i=1}^{\left\lfloor\frac{d+1}{2}\right\rfloor} \lambda_ih(Q_i) +\sum_{i=\left\lfloor\frac{d+1}{2}\right\rfloor+1}^{d-1} \lambda_ih(Q_i)=\lambda'h(Q_1)+ \sum_{i=\left\lfloor\frac{d+1}{2}\right\rfloor+1}^{d-1} \lambda_ih(Q_i)\]
where $\lambda' = \displaystyle{\sum_{i=1}^{\left\lfloor\frac{d+1}{2}\right\rfloor} \lambda_i}$.   Since these $h$-vectors are affinely independent \cite{BB}, $\lambda' = 0$ and $\lambda_i=0$ for $i =\left\lfloor\frac{d+1}{2}\right\rfloor+1, \ldots, d-1$ and we have an affine dependence among the $\gamma$-$h$-vectors, and hence the $\gamma$-vectors, of the simplices.
However, this is impossible since the $\gamma$-vectors of the simplices are affinely independent by the previous theorem. Therefore the $\gamma$-$h$-vectors and the $\alpha$-$f$-vectors of the $d$ constructed polytopes are affinely independent. 
\end{proof}

As an example, we can consider the set of simplicial polytopes in $\R^3$ whose $\alpha$-$f$-vectors span the affine space $A_S$.  This space is 2-dimensional and is spanned by the two simplices $P_{\delta}P^2$ and $P_NP^2$, one tall and one short tetrahedron, and by $T^3_1$, formed by stellar subdivision of a facet of a tetrahedron.  This last can also be thought of as a bipyramid over a triangle or two tetrahedra glued together along one face.  To pick a particular basis of $\alpha$-$f$-vectors, take the limiting cases of the simplices, $P_{0}P^{2}$ and $P_{\infty}P^{2}$, and the geometric realization of $T^3_1$ made by gluing together two regular tetrahedra.  This results in the following $\alpha$-$f$-vectors:
\begin{align*}
\alpha\text{-}f(P_{0}P^{2}) & = \left(\frac{1}{2}, \frac{3}{2}, 2, 1, 4,6,4,1\right)\\
\alpha\text{-}f(P_{\infty}P^{2}) & =\left(\frac{1}{4},\frac{5}{4},2,1,4,6,4,1\right)\\
\alpha\text{-}f(T^3_1) & = \left(\frac{6}{\pi}\arccos\left(\frac{1}{3}\right)-2,\frac{6}{\pi}\arccos\left(\frac{1}{3}\right),3,1,5,9,6,1\right)
\end{align*}

For the $\alpha$-$f$-vectors of simplicial polytopes, this shows that the dimensions beyond those determined combinatorially are found in variation of the angle sums of simplices.  This means that degrees of freedom in the geometry of simplicial polytopes beyond that of the simplex are purely combinatorial.

We can similarly build a set of polytopes whose $\alpha$-$f$-vectors span the space defined by the Gram and Euler relations. We will use a method similar to the proof of Theorem \ref{Angle sums of simplices}, but first we will prove the following lemma.
\begin{lemma}\label{inverse image of P on B*Q}
Let $Q$ be a $d$-polytope with $f$-vector $f = (f_0, f_1, \ldots, f_d)$ and let $\f = (1, \f_0, \f_1, \ldots, \f_{d-1})$ be the inverse image of $f$ under the pyramid transformation.  Also, let the $(d+1)$-polytope $B^{*}Q$ have $f$-vector $f^{*} = (f^{*}_0, f^{*}_1, \ldots, f^{*}_{d+1})$ and inverse image $\f^{*} = (1, \f^{*}_0, \f^{*}_1, \ldots, \f^{*}_{d})$ under the pyramid transformation.  Then \[\sum_{i=0}^{d} (-1)^i\f^{*}_i = \sum_{i=0}^{d-1} (-1)^i\f_i +1.\]
\end{lemma}
\begin{proof}
By \eqref{f and P}, if we extend the $f$-vector of a polytope $Q$ to $(1, f(Q))$,
\begin{equation}\label{Matrix B}
f(PQ) = \textbf{B}\begin{bmatrix} 1\\f(Q)\end{bmatrix},
\end{equation}
where
\begin{equation*}
\textbf{B} = \left[\begin{array}{ccccc}
 1&1&&&\\
 &1&1&&\\
 &&\ddots&\ddots&\\
 &&&1&1\\
 &&&&1
\end{array}\right],
\end{equation*}
a $(d+1)\times(d+1)$ matrix where all other entries are 0.  Therefore the matrix for the inverse transformation is
\begin{equation}\label{Binverse}
\textbf{B}^{-1} = \left[\begin{array}{rrrrr}
 1&-1&1&-1&\ldots\\
 &1 & -1 & 1 & \ldots\\
 &&\ddots&\ddots&\\
 &&&1&-1\\
 &&&&1
\end{array}\right],
\end{equation}
where all the entries below the main diagonal are 0.  Multiplication by $\textbf{B}^{-1}$ gives that 
\[\f_i = \sum_{j=i+1}^d (-1)^{i-j+1} f_j \quad \text{for } 0 \leq i \leq d-1\]
 and
 \begin{equation}\label{barf* to f*}
 \f^{*}_i = \sum_{j=i+1}^{d+1} (-1)^{i-j+1} f^{*}_j \quad \text{for } 0 \leq i \leq d.
\end{equation}
Since $f$ and $f^{*}$ are $f$-vectors of polytopes, we know the last entry is 1, i.e. $f_d = f^{*}_{d+1} = 1$.  We also know from \eqref{alpha and B*} that $f^{*}_i = 2f_i + f_{i-1}$ for $1 \leq i \leq d$.  Therefore, for $0 \leq i \leq d$, we can rewrite \eqref{barf* to f*} as follows:
\begin{equation*}
\begin{split}
\f^{*}_i &= 2\sum_{j=i+1}^{d} (-1)^{i-j+1}f_j + \sum_{j=i}^{d-1} (-1)^{i-j}f_{j}+(-1)^{i-d}f^{*}_{d+1}\\
& =f_i + \sum_{j=i+1}^{d} (-1)^{i-j+1}f_j \\
& = f_i + \f_i.
\end{split}
\end{equation*}
Now, taking the alternating sum we get
\begin{equation*}
\begin{split}
\sum_{i=0}^{d} (-1)^i\f^{*}_i & =\sum_{i=0}^{d} (-1)^i f_i + \sum_{i=0}^{d} (-1)^i \f_i\\
& = 1 + \sum_{i=0}^{d} (-1)^i \f_i,
\end{split}
\end{equation*}
where the last equality follows by the Euler relation on $Q$.
\end{proof}

In the proof of the next theorem, we will work with the $\gamma$-$f$-vector of a $d$-polytope $Q$: $(\gamma_0(Q), \ldots, \gamma_d(Q)| f_0(Q) \ldots, f_d(Q))$.
\begin{theorem}\label{Angle sums of general polytopes} The affine span of the $\alpha$-$f$-vectors of general $d$-polytopes has dimension $2d-3$ for $d \geq 2$.
\end{theorem}
\begin{proof}
The Euler and Gram equations provide two independent equations on the $\alpha$-$f$-vectors.  We also know that $\alpha_{d}(P)=1$, $f_{d}(P)=1$, and $\alpha_{d-1}(P)=\frac{1}{2}f_{d-1}(P)$ for all polytopes $P$.  As long as $d>1$, these equations are independent.  Therefore, the span of the $\alpha$-$f$-vectors is at most $2d+2-5=2d-3$ if $d\geq 2$. To show this whole space is spanned, we will again proceed inductively on $d$.

The statement is true in two dimensions, since the $\alpha$-$f$-vectors of the triangle and the square (denoted $P^{2}$ and $B^{*}P$, respectively) are $\left(\frac{1}{2},\frac{3}{2}, 1,3,3,1\right)$ and $(1,2,1,4,4,1)$. This gives an affine span of dimension 1.

Suppose the statement is true for dimension $d-1$.  That is, there are $2(d-1)-2=2d-4$ affinely independent $\alpha$-$f$-vectors of $(d-1)$-polytopes: $Q_{1}, Q_{2}, \ldots, Q_{2d-4}$. Then we claim that the $\alpha$-$f$-vectors of the polytopes 
\[P_{\infty}Q_{1}, P_{\infty}Q_{2}, \ldots, P_{\infty}Q_{2d-4},(B^{*})^{d-2}P^{2} \text{ and }(B^{*})^{d-1}P\] 
are affinely independent.  %We show this with an argument similar to that for Theorem \ref{Angle sums of simplices}.  

Since each $\alpha$-$f$-vector has $\alpha_d(\P) = f_d(\P) = 1$, affine independence of a set of $\alpha$-$f$-vectors is equivalent to their linear independence.  Also, since the linear transformation from the $\alpha$-vector to the $\gamma$-vector is invertible and independent of the $f$-vector, the linear independence of a set of $\alpha$-$f$-vectors is equivalent to the linear independence of the corresponding set of $\gamma$-$f$-vectors.  Therefore, we will show that the $\gamma$-$f$-vectors of $P_{\infty}Q_{1}$, $P_{\infty}Q_{2}$, \ldots, $P_{\infty}Q_{2d-4},(B^{*})^{d-2}P^{2}$ and $(B^{*})^{d-1}P$ are linearly independent. 

We will consider the effect of the $P_{\infty}$ construction on the $\gamma$-$f$-vector of a $(d-1)$-polytope $Q$.  To do this, we extend the $\gamma$-$f$-vector to the $2d$-vector $(\gamma(Q), 1, 1, f(Q))$, thinking of the additional entries as $\gamma_{d}(Q)$ and $f_{-1}(Q)$, respectively. Then we can write
\begin{equation}\label{Pinfty trans on gamma-f}
\gamma\text{-}f(P_{\infty}Q) = \textbf{C}\begin{bmatrix} \gamma(Q)\\1\\1\\f(Q)\end{bmatrix}
\end{equation}
where
\begin{equation*}
\textbf{C}=  \left[\begin{array}{c|c} \textbf{A}& 0 \\\hline 0 & \textbf{B}\end{array}\right],
\end{equation*}
a $(2d+2)\times(2d+2)$ matrix, with blocks $\textbf{A}$ \eqref{Pinfty trans on gamma} and $\textbf{B}$ \eqref{Matrix B}. Since this is an invertible matrix, the $\gamma$-$f$-vectors of $P_{\infty}Q_{i}$, $i=1,\ldots,2d-4$, are all linearly independent.  Therefore, we consider the $\gamma$-$f$-vectors of $(B^{*})^{d-2}P^2$ and $(B^{*})^{d-1}P$ in particular.

As in Theorem \ref{Angle sums of simplices}, we will consider the inverse images of the $\gamma$-$h$-vectors of these two polytopes in the $P_{\infty}$ transformation and show that they must increase the dimension by two. Let 
\[v_1:=\textbf{C}^{-1}\left(\gamma\text{-}f((B^{*})^{d-2}P^2)\right)^{T} \qquad \text{and} \qquad
v_2:=\textbf{C}^{-1}\left(\gamma\text{-}f((B^{*})^{d-1}P)\right)^{T} \]
and consider these vectors in relation to the span of the $(\gamma(Q_i), 1, 1, f(Q_i))$.

First we note that the alternating sum of the values of $f(Q)$ is 1 for any polytope $Q$ by the Euler relation.  Since this is the same value as $\alpha_d(Q)$ and $f_{-1}(Q)$, any polytope in the linear span of the $\gamma$-$f$-vectors of the $Q_i$ for $i=1,\ldots,2d-4$ must have the corresponding alternating sum equal to the $\alpha_d$ and $f_{-1}$ entries.  Since the $f_{-1}$ entry of $\textbf{C}^{-1}\left(f(Q)\right)$ for any $d$-polytope $Q$ is 1 by the Euler relation, the $f_{-1}$ entries of $v_1$ and $v_2$ are 1.  

However, by Lemma \ref{inverse image of P on B*Q}, the alternating sum of the entries $f_0, f_1, \ldots, f_{d-1}$ of $\textbf{C}^{-1}\left(f\left(B^{*}Q\right)\right)$ for a $(d-1)$-polytope $Q$ is one greater than the alternating sum of the entries $f_0, f_1, \ldots, f_{d-2}$ of $\textbf{C}^{-1}\left(f\left(Q\right)\right)$. If $Q$ is a pyramid over a polytope, this alternating sum is 1 by the Euler relation and therefore the alternating sum of the entries $f_0, f_1, \ldots, f_{d-1}$ of $v_1$ and $v_2$ must be at least 2 since $d \geq 3$.  Therefore we can see that neither can be a linear combination of the extended $\gamma$-$f$-vectors of $Q_{i}$ for $i=1,\ldots,2d-4$. Then if
\[k= \dim\left(\spanr\{(\gamma(Q_i), 1, 1, f(Q_i)),
v_1,v_2\}\right) - \dim\left(\spanr\{\gamma\text{-}h(Q_i)\}\right),\]
\[1\leq k \leq 2.\]

Considering the $\gamma$-vector entries will show that $k=2$.
By Proposition \ref{gamma and B*},  $\gamma((B^{*})^{d-2}P^2) = \left(0, \frac{1}{2}, 1, \ldots,1\right)$ and
$\gamma((B^{*})^{d-1}P) = \left(0,1,\ldots, 1\right).$
Therefore the $\gamma$-vector portions of $v_1$ and $v_2$ are
$(0,1, \ldots ,1)$ and $(0,2,0 ,2,\ldots)$, respectively.  Since $v_2$ has $\gamma_d \neq 1$, but $f_{-1}=1$, $\gamma_d \neq f_{-1}$ in $v_2$, even though $\gamma_d =f_{-1}$  for each of the vectors $v_1$ and $(\gamma(Q_i), 1, 1, f(Q_i))$ for $i = 1, \ldots, 2d-4$.  Therefore 
\[v_2 \notin \spanr\{(\gamma(Q_i), 1, 1, f(Q_i)), v_1\}\]
and the following set of $d$-polytopes inductively constructed above for dimension $d$ have affinely independent $\alpha$-$f$-vectors:
\begin{equation*}
\begin{split}
(P_{\infty})^{d-2}P^2, & (P_{\infty})^{d-2}B^{*}P, (P_{\infty})^{d-3}B^{*}P^2, (P_{\infty})^{d-3}(B^{*})^2P,\\
& \ldots,  P_{\infty}(B^{*})^{d-3}P^2, P_{\infty}(B^{*})^{d-2}P, (B^{*})^{d-2}P^2
, (B^{*})^{d-1}P.
\end{split}
\end{equation*}

Then by Lemma \ref{backing off limiting}, we know that we have a set of $d$-polytopes of size $2d-2$ with affinely independent $\alpha$-vectors. 
\end{proof}

The set of 3-dimensional polytopes given by the theorem that affinely span the space of all $\alpha$-$f$-vectors of polytopes has 4 elements: $(B^{*})^{2}P$, $B^{*}P^{2}$, $P_{N}P^2$, and $P_{N}B^{*}P$ for large enough $N$.
%Compare to Bayer & Billera set??  More?

We note that the set of polytopes which span the space of $\alpha$-$f$-vectors has significant duplication in the $\alpha$-vectors.  For instance, the polytopes $P_{\infty}(B^{*})^{k}P$ and $(B^{*})^{k}P^{2}$ have the same angle sums for all $k \geq 1$.

These results strengthen the correspondence between the geometric structure and the combinatorial structure of polytopes.  The Gram and Perles relations are close analogs of the Euler and Dehn-Sommerville relations.  In this chapter, we have shown that the affine dimensions closely correspond.  The affine span of the $\alpha$-vectors of $d$-simplices has the same dimension as the span of the $f$-vectors of simplicial $(d-1)$-polytopes.  Also, the affine span of the $\alpha$-$f$-vectors of simplicial $d$-polytopes has the same dimension as the span of the $f$-vectors of $d$-polytopes.  It would be interesting to speculate whether there is a deeper significance to this relationship. 

The use of the $\gamma$-vector also raises questions about the nature of this measure on angle sums.  For the $h$-vectors of simplicial polytopes, there are many results bounding the values.  The Upper Bound Theorem \cite{McM3} bounds the $h$-vector entries above by those of the cyclic polytope of the same dimension and the same number of vertices.  The Generalized Lower Bound Theorem \cite{McMWa} shows that the first $\left\lfloor \frac{d}{2} \right \rfloor$ entries of the $h$-vector are unimodal and the $g$-Theorem gives bounds on the differences between adjacent entries of the vector \cite{BL, Stanley}.  In the case of the $\gamma$-vector the bounds on its entries are unexplored.  Initial examples show that the $\gamma$-vector may be more tractable on non-simplicial polytopes than the $h$-vector; for example, the basis polytopes for the theorems in this chapter, many of which are not simplicial, all have non-decreasing $\gamma$-vectors.  This is not the case for all polytopes (for example, the bipyramid made by gluing two regular tetrahedra along a face), but unimodality may be true in general and monotonicity in specific cases such as for simplices.

\chapter{Angle Sums on Complexes}

%%Thesis Chapter on Complexes

In this chapter we will consider angle sums on polytopal complexes. We define the angle characteristic on the $\alpha$-vector in analogy with the Euler characteristic on the $f$-vector and study the effect of a few constructions upon Gram- and Perles-type relations. In particular, we will create a large set of odd-dimensional complexes whose angle characteristic is half the Euler characteristic and another large set of polytopes which satisfy the Perles relations. 

\section{Defining Angle Sums for Complexes}
%Defining angle sums for complexes- 
$\C$ is called a \emph{polytopal complex} if it is a cell complex where each cell is a convex polytope and the following two properties are satisfied:
\begin{itemize} 
\item If $\F \in \C$ and $\G$ is a face of $\F$, then $\G \in \C$.
\item If $\F, \G \in\C$, then $\F \cap \G$ is a face of both $\F$ and $\G$.
\end{itemize}

We will assume that a $d$-complex $\C$ is connected, pure, and embedded in $\R^d$. These requirements are chosen to simplify the complexes and make the measure of interior angles unambiguous. We require that the complex be pure because if a line segment met a triangle at a vertex $v$, the interior angle at $v$ might be measured into the segment, the triangle or some undefined composite thereof.  By assuming that $\C$ is $d$-dimensional and embedded in $\R^d$, we have guaranteed that the complex has a well-defined and geometrically fixed $d$-dimensional interior, $\int(\C)$, making it possible to define angles consistently since each angle has ambient space of the same dimension.  %If the complex consisted of the edges of a cube in $\R^4$, the angle at a vertex could be $3/2$ if we take the sum of angle into the edges or 0 if we consider the angles into the cube as the fraction of $\S^3$ in $\R^4$.  %Could add example of when it doesn't work...
%Is it semi-Eulerian or a homology manifold?
Since we have fixed the geometry of the complex so that each polytopal face is represented by a polytope in $\R^d$, we will let $|\C|$, the \emph{geometric realization} of $\C$, be the set of points in $\R^d$ which is the union of the sets of points in each of the polytopes which is a member of $\C$.  We will consider only the angles at boundary faces. %We could have just considered boundary complexes to begin with, but by considering the larger complex we are guaranteed nice properties...

Then the \emph{interior angle at a boundary face $\F$ of $\C$} is defined as 
\[\alpha(\F, \C) = \frac{vol\left(S_\e(x) \cap \int(\C)\right)}{vol\left(S_\e(x)\right)},\]
where $x$ is in the interior of $\F$ and $S_\e(x)$ is the $(d-1)$-sphere of radius $\e$ centered at $x$ for $\e$ sufficiently small.  Then the \emph{angle sums of $\C$} are defined as 
\[\alpha_{i}(\C)= \sum_{i-\text{faces } \F\subseteq \partial \C} \alpha(\F, \C) \quad \text{for }0\leq i\leq d-1.\]
Since we are only working with boundary faces, we will continue to think of $\alpha_d(\C)=1$, thinking of this as counting a connected interior.  Since disconnected complexes will have angle sums that are the sum of angle sums on the connected components, we will assume that we are working with connected polytopal complexes. We will call a polytopal complex that is connected, pure, and embedded in $\R^d$ a \emph{geometric polytopal complex}.
 
\begin{figure}
\begin{center}
\includegraphics[width=3.5in]{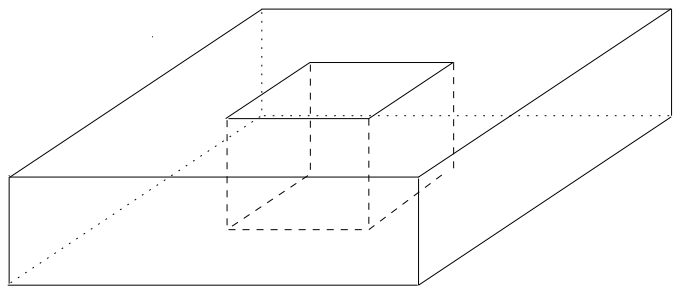}
\caption{The geometric torus $|T|$.}\label{geotorus}
\end{center}
\end{figure}

\begin{figure}
\begin{center}
\includegraphics[width=3.5in]{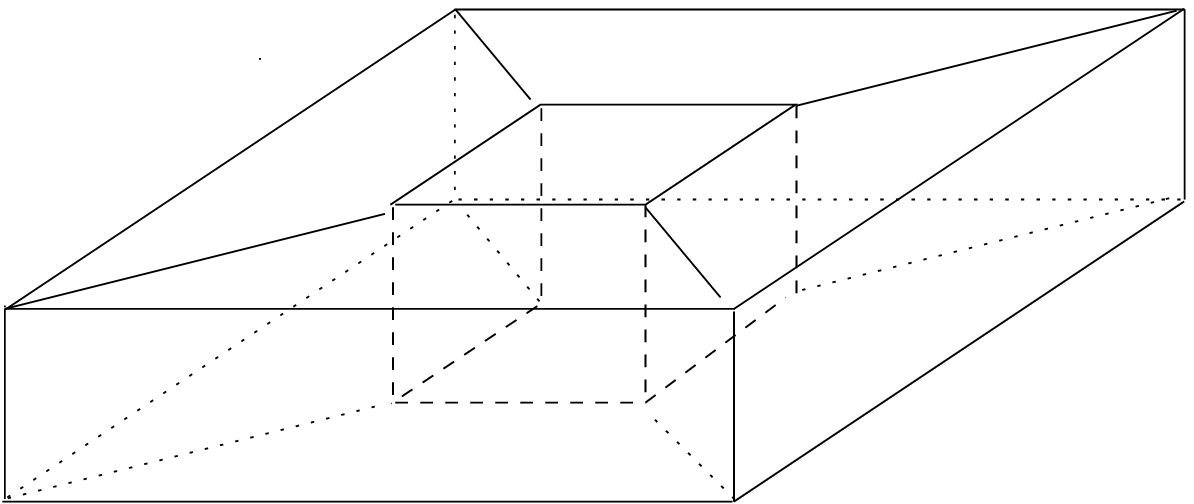}

\vspace{.3in}
\includegraphics[width=3.5in]{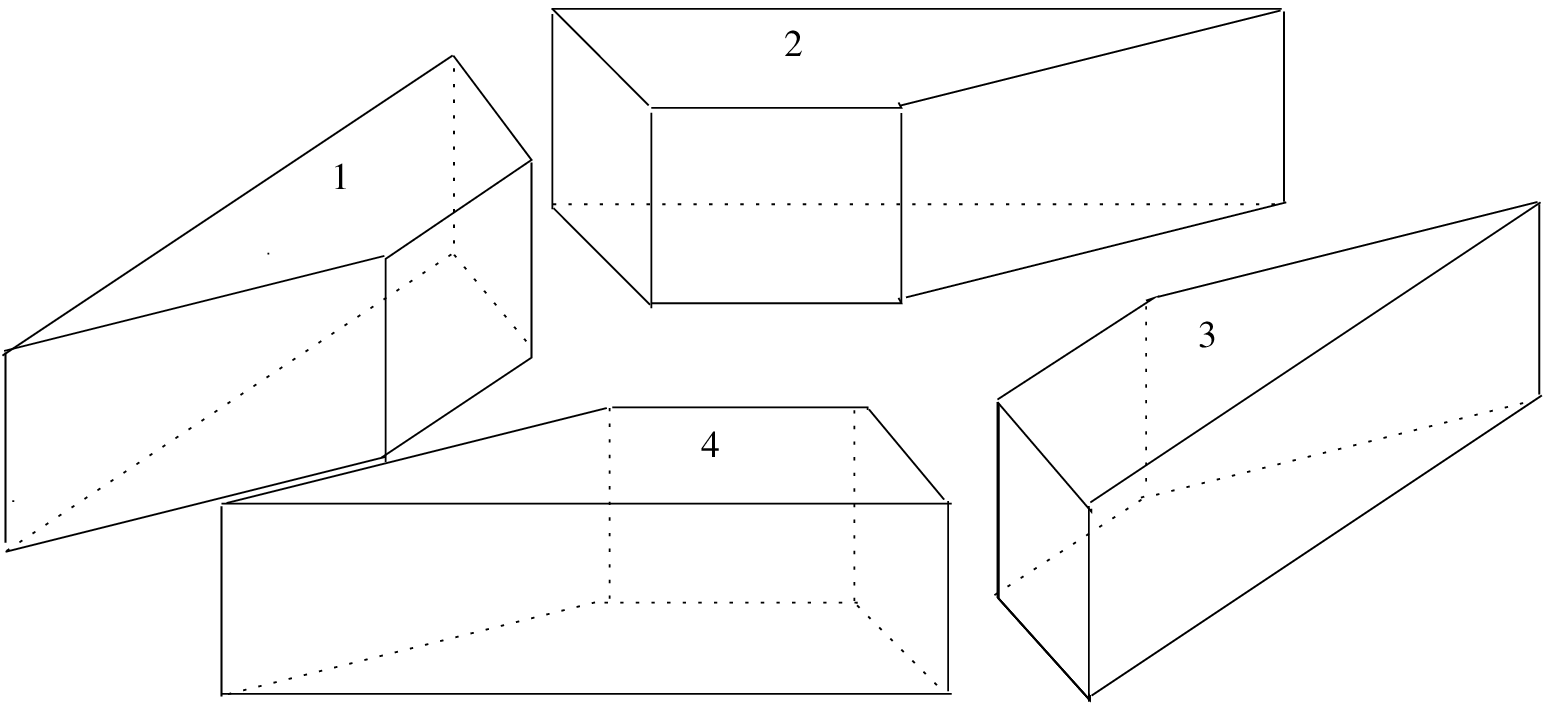}
\caption[A decomposition of a torus complex.]{At top, the subdivision of the boundary of the torus $T$, which is the boundary of the polytopal complex decomposed below.}\label{torusdecomp}
\end{center}
\end{figure}

For a complex, $\alpha(\F, \C)$ can be greater than $\frac{1}{2}$, in contrast to interior angles in convex polytopes.  For example, consider a complex whose geometric realization is a solid torus, $|T|$, where all angles are right angles, as in Figure \ref{geotorus}.
%EXAMPLE: TORUS
\singlespacing
\normalspacing
As this is, it is not the boundary of a polytopal complex, since some of the boundary faces, such as the top and bottom, are not homeomorphic to spheres.  Therefore we will subdivide these ring faces into four quadrilaterals to get the boundary complex of a polytopal complex made from four trapezoidal prisms.  The subdivision of the boundary and the polytopal decomposition of the solid torus are shown in Figure \ref{torusdecomp}.
\singlespacing
\normalspacing
For a vertex or edge on the outer ring of $T$, the angles are the same as the corresponding faces of a cube.  Similarly, the horizontal edges on the inner ring of $T$ have the same interior angle as edges of a cube.  However, for a vertex $v$ and a vertical edge $e$ on the inner ring of the torus, $\alpha(v, \C) = \frac{3}{8}$ and $\alpha(e, \C) = \frac{3}{4}$.  Also notice that an edge $e'$ that was inserted to subdivide the top and bottom of the torus has $\alpha(e', \C) = \frac{1}{2}$, the same as the interior angle at the facet it subdivides.  We then compute that 
\begin{equation*}
\begin{split}
\alpha_0(\C) & = 8\left(\frac{1}{8}\right) + 8\left(\frac{3}{8}\right) = 4, \\ 
\alpha_1(\C) & = 20\left(\frac{1}{4}\right) + 4\left(\frac{3}{4}\right) + 8\left(\frac{1}{2}\right)= 12, \\
\alpha_2(\C) & = 16\left(\frac{1}{2}\right) = 8, 
\end{split}
\end{equation*}
so \[\alpha_0(\C) - \alpha_1(\C) + \alpha_2(\C) = 0.\]
This is the same as the Euler characteristic of the torus.  If we consider the boundaries of odd-dimensional polytopes, the Gram relation gives us that the alternating sum of angle sums is 1, half the Euler characteristic, a relationship that would also apply to the results for the torus.

With these motivating examples, we will use our knowledge of the Gram and Perles relations for polytopes as a basis to determine relations on angle sums for a variety of complexes.  

%Angle characteristics
\section{The Angle Characteristic}
We define the following operator on a geometric polytopal complex $\C$ to study the patterns on the alternating sum of angle sums:
\[\chi_{\alpha}(\C) = \sum_{i=0}^{d-1} (-1)^{i} \alpha_{i} (\C).\]
Then for $d$-polytopes, $\P$, we can rewrite the Gram relation as $\chi_{\alpha}(\P) = (-1)^{d-1}$.  We will call $\chi_{\alpha}$ the \emph{angle characteristic} to parallel the Euler characteristic. 

In a first step to justifying the name of angle characteristic, we show that the angle characteristic is independent of the subdivision of the complex, just as the Euler characteristic is.

If a face in $\partial\C$ is the intersection of a supporting hyperplane $H$ with some polytope $\P$ in $\C$, the maximal connected part of $H \cap \partial\C$ which includes the face will be called a \emph{flat}. For the torus in Figure \ref{geotorus}, the horizontal flats would include ring-shaped 2-faces, so flats are not necessarily polytopal. We can decompose the set of faces of $\C$ into subcomplexes $\C^{\circ}(\F^{*})$, the complex of faces that are in the relative interior of some flat $\F^{*}$.  For every $\G \in \C^{\circ}(\F^{*})$, $\alpha(\G, \C)$ has the same value, and is equal to the interior angle from any interior point of the flat $\F^{*}$.  So we define $\alpha(\F^{*}, \C) = \alpha(\G, \C)$ for some $\G \in \C^{\circ}(\F^{*})$.  
\begin{lemma}\label{independence of subdivision}
Let $\C$ be a geometric polytopal complex and let $\C^{*}$ be the set of flats. Then
\[\chi_{\alpha}(\C) = \sum_{\F^{*} \in \C^{*}} (-1)^{\dim (\F^{*})}\alpha(\F^{*}, \C)\chi(\int(\F^{*})). \]
Therefore, if $\A$ and $\B$ are two $d$-complexes that have the same geometric realization, $|\A| = |\B|$, then $\chi_{\alpha}(\A) = \chi_{\alpha}(\B)$ and $\chi_{\alpha}(\C)$ is independent of the subdivision of $|\C|$.
\end{lemma}
\begin{proof}
We will show that the contribution of the faces of one flat, $\C^{\circ}(\F^{*})$, to $\chi_{\alpha}(\C)$ is dependent only on the topology of $\F^{*}$.

The contribution of $\C^{\circ}(\F^{*})$ to the angle characteristic of $\C$ is 
\begin{equation*}
\begin{split}
\sum_{\G \in C\left(\F^{*}\right)} (-1)^{\dim (\G)}\alpha\left(\G, \C\right)  & = \alpha\left(\F^{*}, \C\right)\sum_{\G \in C\left(\F^{*}\right)} (-1)^{\dim (\G)} \\
& = \alpha\left(\F^{*}, \C\right)\sum_{i=0}^{\dim(\F)} (-1)^{i}f_i\left(C\left(\F^{*}\right)\right) \\
& = \alpha\left(\F^{*}, \C\right)\chi\left(C\left(\F^{*}\right)\right).
\end{split}
\end{equation*}
However, $\chi(C^{\circ}(\F^{*}))$ is independent of the subdivision of $\left| \C^{\circ}(\F^{*})\right|$ since the Euler characteristic is a topological invariant and we have a fixed geometry for $C^{\circ}(\F^{*})$.

Since each face of $\C$ is contained in $C^{\circ}(\F^{*})$ for exactly one flat $\F^{*}$,  
\[\chi_{\alpha}(|\C|) = \sum_{\G \in \C} (-1)^{\dim (\G)} \alpha(\G, \C) = \sum_{\F^{*} \in \C^{*}} (-1)^{\dim (\F^{*})}\chi(\int(\F^{*}))\alpha(\F^{*}, \C) \]
and the angle characteristic of a complex is dependent only on its set of flats.

Flats are independent of subdivision and dependent only on the geometric realization, so if $|\A|= |\B|$, then 
\[\{\F_{\A}^{*}: \F_{\A}^{*}\text{ is a flat of }\A\} = \{\F_{\B}^{*}: \F_{\B}^{*}\text{ is a flat of }\B\}\]
and
\[\alpha(\F^{*}, \A) = \alpha(\F^{*}, \B).\]
 Therefore, $\chi_{\alpha}(|\A|) = \chi_{\alpha}(|\B|)$ .
\end{proof}

Based on this lemma, we will consider the angle characteristic as acting on the geometric realization rather than on a particular complex.  This will be implicit in the notation for the angle characteristic, where we will write $\chi_{\alpha}(\C)$ and $\chi_{\alpha}(|\C|)$ interchangeably.  The lemma also allows us to think of the angle characteristic as a weighted alternating sum on the geometric realization with weights given by the Euler characteristic.

%Building via connected sums & identifications
\section{Constructions of Complexes and their Angle Characteristics}
We will consider a few variants of a basic construction on complexes and consider their effect on both the angle characteristic and the Euler characteristic.

\begin{figure}
\begin{center}
\includegraphics[width=4in]{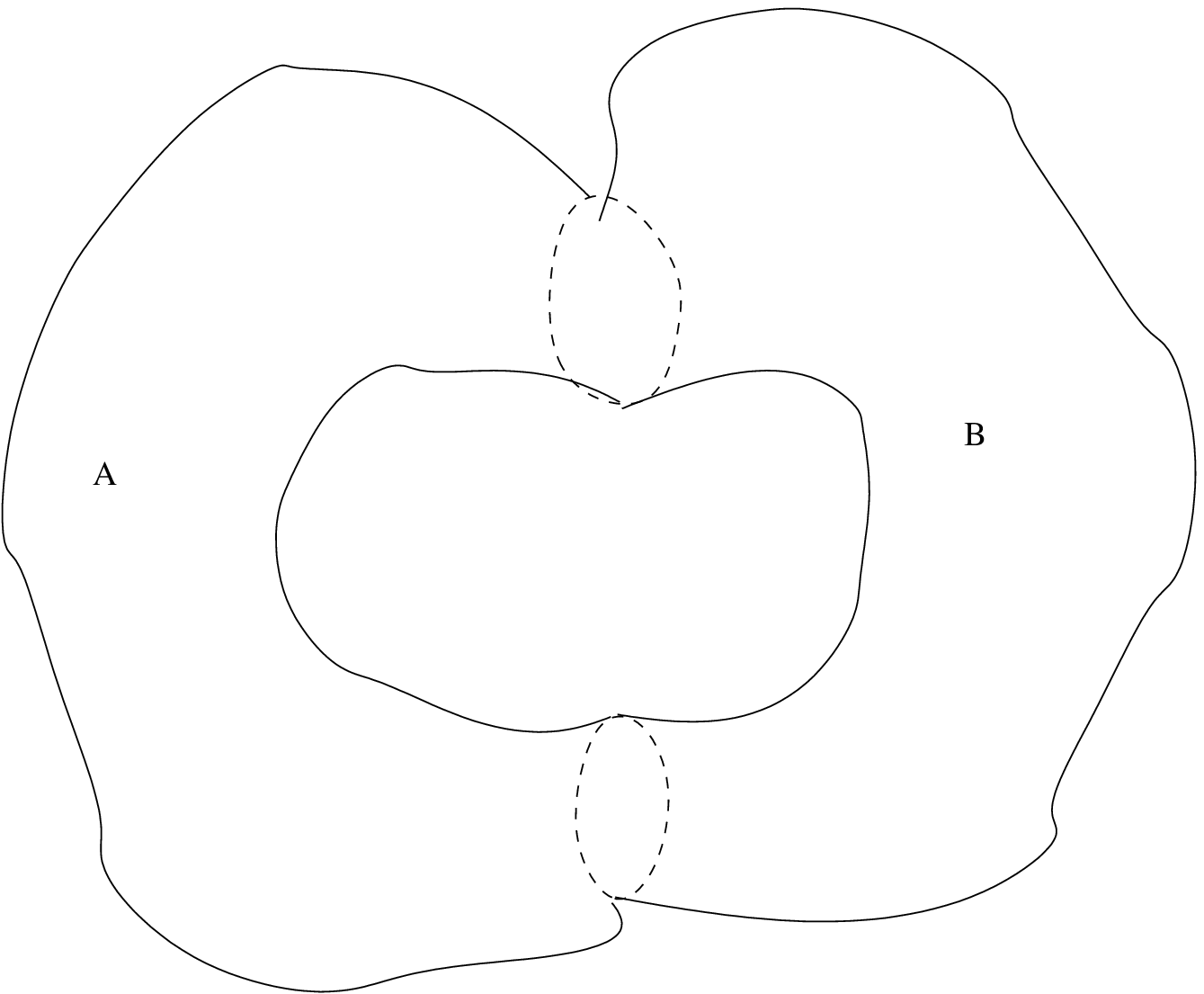}
\caption[A gluing, $\C = \A \glue \B$.]{Two complexes, $\A$ and $\B$, which meet along a sub-complex. The resulting complex is $\C = \A \glue \B$.}\label{gluing}
\end{center}
\end{figure}

If we have a $d$-complex $\C$ that can be decomposed into two disjoint complexes, $\A$ and $\B$, along a sub-complex $\C'=\A \cap \B  \subseteq \C$ which is contained in the boundary of $\A$ and $\B$, we will write $\C = \A \glue \B$ and say that $\C$ is the \emph{gluing} of $\A$ and $\B$.  The sub-complex $\A \cap \B$ does not need to be connected, as illustrated in Figure \ref{gluing}.  We will consider how the Euler and angle characteristics of $\C$ are related to those of $\A$ and $\B$.  Since we consider only the boundary faces when working with angle sums, we will also look at Euler characteristics of the boundary complex for comparison.
\singlespacing
\normalspacing

If $\A \cap \B$ is a $(d-1)$-polytope, $\P$, then the gluing is equivalent to a connected sum along $\P$.  If $\B$ is homeomorphic to a $d$-ball and $\A \cap \B$ consists of two $(d-1)$-polytopes, $\P$ and $\P'$ then $\A \glue \B$ is topologically the same as a handle addition to $\A$ by identifying $\P$ and $\P'$. However, since we have fixed the geometry of the complexes, we must add to the geometric realization of $\A$ in $\R^d$ to perform a handle addition rather than just identify $\P$ and $\P'$.  Therefore, we can think of the gluing construction as a generalization of the connected sum and handle addition constructions.

We will consider $\A \cap \B$ as a $(d-1)$-dimensional geometric polytopal complex.  That is, we will let $\int(\A \cap\B)$ be the set of points in $|\A \cap\B|$ that are contained in an open $(d-1)$-dimensional subset of $|\A \cap\B|$ and $\partial(\A \cap \B) = (\A \cap \B)\setminus\int(\A \cap\B)$.  Then if $\A \cap \B$ has dimension less than $d-1$, $\int(\A \cap\B) = \emptyset$ and $\partial(\A \cap \B)=\A \cap \B$.
%I don't like the fact that this doesn't look parallel - I don't want to have to write the boundary operator! 
\begin{lemma}\label{valuations}
If a $d$-complex $\C=\A \glue \B$, then the following relations hold on the Euler and angle characteristics:
\[\chi\left(\partial \C\right) = \chi\left(\partial \A\right) +\chi\left(\partial \B\right) - 2\chi\left(\int\left(\A \cap \B\right)\right) - \chi\left(\partial\left(\A \cap \B\right)\right)\] 
\[\chi_{\alpha}\left(\C\right) = \chi_{\alpha}\left(\A\right)+\chi_{\alpha}\left(\B\right) - \chi\left(\int\left(\A \cap \B\right)\right) \]
%Should include the valuations on face numbers and angle sums in the statement.
\end{lemma}
\begin{proof}
We consider what happens to the $f$-vectors and $\alpha$-vectors of $\A$, $\B$ and $\C$. Let $\F$ be an $i$-face.  If $\F \in \partial \A$, it is either (i) in $\int(\A \cap \B)$ or (ii) in $\partial \C$.  In the first case, $\F$ contributes nothing to either $f_i(\C)$ or $\alpha_{i}(\C)$, but it does contribute to $\alpha_i$ and $f_i$ for both $\A$ and $\B$. In the second case, if $\F$ is not in $\A \cap \B$, then the contribution to $f_{i}(\partial \C)$ or $\alpha_{i}(\C)$ is the same as it is to $f_{i}(\partial \A)$ or $\alpha_{i}(\A)$ and it makes no contributions to $\B$.  If, on the other hand, $\F \in \partial(\A \cap \B)$, $\F$ contributes to both $f_{i}(\A)$ and $f_{i}(\B)$, double-counting the contribution to $f_{i}(\C)$.  However, there is no over-counting in the angle sum since the angles at $\F$ in $\A$ and $\B$ are concatenated to make the angle in $\C$.  Therefore we have the following relations on face numbers and angle sums:
\begin{equation}\label{f-valuations}
f_i\left(\partial \C\right)  = f_i\left(\partial \A\right) + f_i\left(\partial \B\right) - 2f_i\left(\int\left(\A \cap \B\right)\right) - f_i\left(\partial\left(\A \cap \B\right)\right)
\end{equation}
and
\begin{align}\label{a-valuations}
\alpha_i\left(\C\right) & = \alpha_i\left(\A\right) + \alpha_i\left(\B\right) - \sum_{\F \subseteq \int\left(\A \cap\B\right)} \big[\alpha\left(\F, \A\right) + \alpha\left(\F, \B\right)\big]\\
& = \alpha_i\left(\A\right) + \alpha_i\left(\B\right) - f_i\left(\int\left(\A \cap\B\right)\right). \nonumber
\end{align}
The last equality follows from the fact that for $\F\subseteq \int(\A \cap\B)$, $\F$ is an interior face of $\C$, so the angles at $\F$ in $\A$ and $\B$ must add to a full angle, giving a total contribution of 1 for each face in $\int(\A \cap\B)$.
Then the following is immediate:  
\begin{equation*}
\begin{split}
\chi_{\alpha}\left(\C\right) & = \sum_{i=0}^{d-1} (-1)^{i} \alpha_{i} \left(\C\right) \\ 
& = \sum_{i=0}^{d-1} (-1)^{i} \left[\alpha_{i}\left(\A\right) + \alpha_{i}\left(\B\right) - f_{i}\left(\int\left(\A \cap\B\right)\right)\right] \\
& = \chi_{\alpha}\left(\A\right) + \chi_{\alpha}\left(\B\right) - \chi\left(\int\left(\A \cap\B\right)\right)\\
\end{split}
\end{equation*} 
and
\begin{equation*}
\begin{split}
\chi\left(\partial \C\right) & = \sum_{i=0}^{d-1} (-1)^{i} f_{i} \left(\partial\left(\A \glue \B\right)\right) \\ 
& = \sum_{i=0}^{d-1} (-1)^{i} \left[f_i\left(\partial \A\right) + f_i\left(\partial \B\right) - 2f_i\left(\int\left(\A \cap \B\right)\right) - f_i\left(\partial\left(\A \cap \B\right)\right)\right] \\
& = \chi\left(\partial\A\right) + \chi\left(\partial\B\right) - 2\chi\left(\int\left(\A \cap \B\right)\right) - \chi\left(\partial\left(\A \cap \B\right)\right).
\end{split}
\end{equation*} 
\end{proof}

This allows us to compute the angle characteristic of complexes made by gluings and shows how changes in the angle characteristic compare to changes in the Euler characteristic.  We will consider a number of different possibilities for $\A \cap \B$ and the effect on the angle and Euler characteristics.
\begin{theorem}\label{Gluing along (d-1)-balls}
If $\C = \A \glue \B$ and $\A \cap \B$ is a union of $m$ disjoint complexes homeomorphic to $(d-1)$-balls, then
\begin{equation}
\chi(\partial \C) = \chi(\partial \A) + \chi(\partial \B) - m\left(1+(-1)^{d-1}\right)
\end{equation}
and
\begin{equation}
\chi_{\alpha}(\C) = \chi_{\alpha}(\A)+\chi_{\alpha}(\B) - m\left((-1)^{d-1}\right).
\end{equation}
In particular, if $d$ is odd, the difference between the Euler characteristic of $\C$ and the sum of the Euler characteristics of $\A$ and $\B$ is twice the corresponding difference on the angle characteristics.  On the other hand, in even dimensions $\chi_{\alpha}(\C)$ is greater than the sum of the angle characteristics of $\A$ and $\B$ while $\chi(\partial\C)$ is the sum of the Euler characteristics of $\A$ and $\B$.
\end{theorem}
\begin{proof}
We know that if $K$ is a $(d-1)$-ball, then $\chi(\int (K))=(-1)^{d-1}$ and that $\partial K$ is a $(d-2)$-sphere, so $\chi(\partial K) = 1+(-1)^{d-2}$. Also, the Euler characteristic of a disconnected complex is the sum of the Euler characteristics of the connected components.

Then, from Lemma \ref{valuations}, 
\begin{align*}
\chi(\partial \C) & = \chi(\partial \A) +\chi(\partial \B) - 2\chi(\int(\A \cap \B)) - \chi(\partial(\A \cap \B))\\
&= \chi(\partial \A) + \chi(\partial \B) - 2m((-1)^{d-1}) - m(1+(-1)^d)\\
& = \chi(\partial \A) + \chi(\partial \B) - m(1+(-1)^{d-1})
\end{align*}
and
\begin{align*}
\chi_{\alpha}(\C) &= \chi_{\alpha}(\A)+\chi_{\alpha}(\B) - \chi(\int(\A \cap \B)) \\
& = \chi_{\alpha}(\A)+\chi_{\alpha}(\B) - m((-1)^{d-1}).
\end{align*} 
\end{proof}

Now, if both $\A$ and $\B$ satisfy the Gram relation, the above simplifies to \[\chi(\partial\C) = (1+(-1)^{d-1})(2-m)\]
and 
\[\chi_{\alpha}(\C) =(-1)^{d-1}(2-m).\]
Therefore, if we begin with odd-dimensional polytopes and glue along sets of disjoint polytopes to get an odd-dimensional complex $\C$, $\chi(\partial\C) = 2\chi_{\alpha}(\C)$. We can then see that Gram's relation will hold for any $d$-complex which can be built from polytopes by gluing each new polytope along a complex homeomorphic to a single ball of dimension $d-1$.
\begin{corollary}\label{Gram's Relation for special spherical polytopal complexes}
If a pure polytopal $d$-complex $\C$ can be built from polytopes by gluing one polytope at a time to the complex along a single complex homeomorphic to a $(d-1)$-ball, then $\C$ satisfies the Gram relation:
\[\chi_{\alpha}(\C) = \sum_{i=0}^{d-1} (-1)^{i}\alpha_{i}(\C) = (-1)^{d-1}.\]
In particular, any complex whose boundary is a stacked sphere satisfies the Gram relation.
\end{corollary}
\begin{proof}
At each stage as we add a polytope, the previous theorem implies that the new angle characteristic is $(-1)^{d-1}(2-1) = (-1)^{d-1}$, so $\C$ satisfies the Gram relation.
\end{proof}

This corollary gives more justification for calling $\chi_{\alpha}(\C)$ the angle characteristic of the complex.  In this case, any complex built in this manner will have a spherical boundary and have the same angle characteristic as a sphere. These complexes are not necessarily convex, thereby extending the set of complexes to which the Gram relation applies.  

Corollary \ref{Gram's Relation for special spherical polytopal complexes} also gives a variant of the Gr\"unbaum proof of the Gram relation.  We can use complexes to accomplish the second and third step of Gr\"unbaum's proof of the Gram relation more directly. We will start with the Gram relation for simplices and build the polytope via the corollary.  Each polytope can be decomposed into simplices via a barycentric subdivision.  The barycentric subdivision of the boundary of a polytope is shellable and the complex formed by taking the pyramid over each boundary face from an interior point is also shellable with the same shelling order as the boundary.  Therefore, these simplices can be ordered so that any intersection between one simplex and all the previous simplices is homeomorphic to a $(d-1)$-ball.  Therefore, we can iteratively add a simplex sharing at least one facet with the previously constructed complex and get the whole polytope.  This construction is shown by the corollary to preserve the angle characteristic, so if we have the Gram relation on simplices, the Gram relation on polytopes results. 

\begin{figure}
\begin{center}
\includegraphics[width=4in]{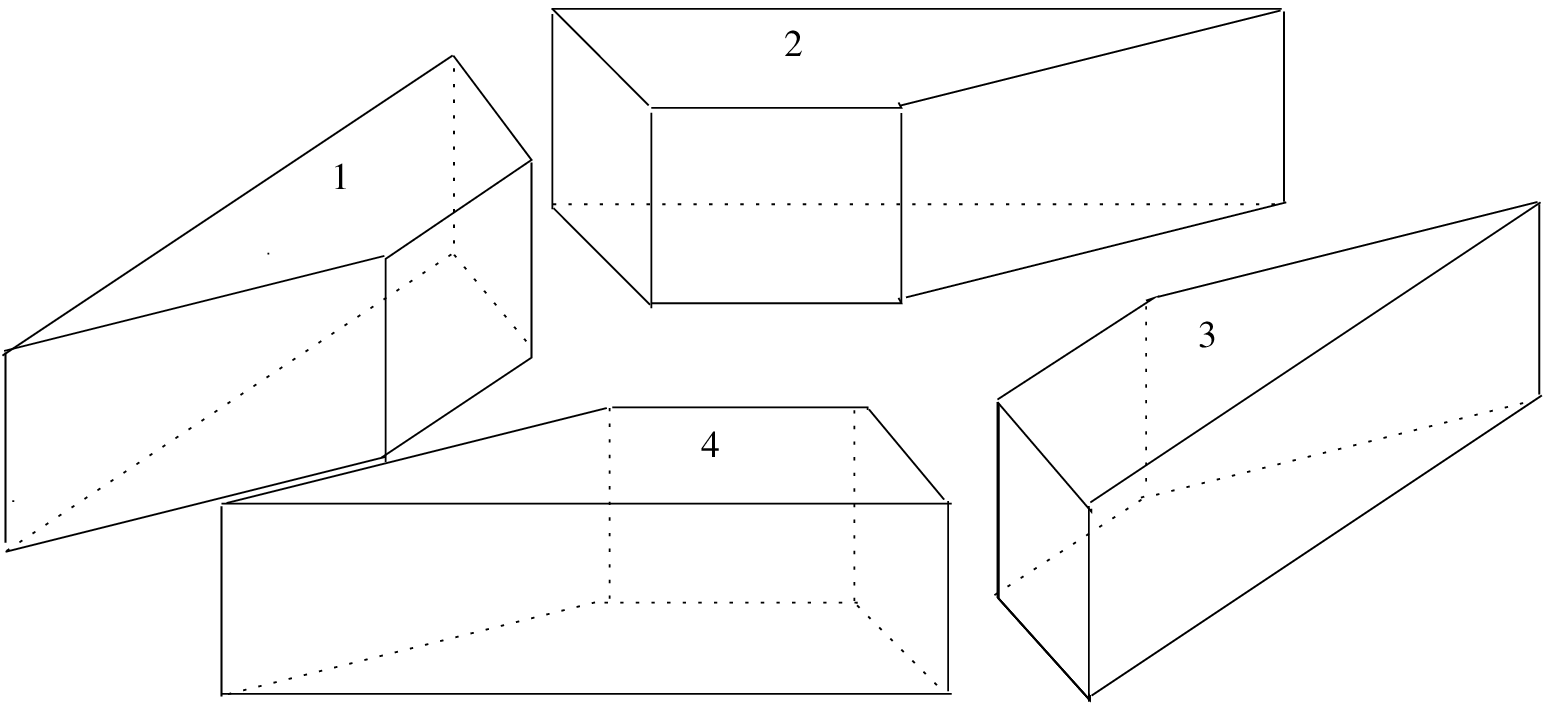}
\caption{A decomposition of the solid torus into four polytopes.}\label{brokentorus}
\end{center}
\end{figure}

Theorem \ref{Gluing along (d-1)-balls} also suggests that changes in the angle characteristic occur when complexes are glued along more than one $(d-1)$-ball, changing the topology.  For example, we can consider a decomposition of the torus complex into polytopes, as in Figure \ref{brokentorus}.
\singlespacing
\normalspacing
We will call the four polytopes $\P_1, \P_2, \P_3$ and $\P_4$, as numbered in the figure.  As $\P_2$ and $\P_3$ are added sequentially to $\P_1$ they meet the previous complex along one 2-face.  When $\P_4$ is added it meets the previous complex along two disjoint 2-faces.  Since each of the polytopes originally follows the Gram relation and the Euler characteristic of a $(d-1)$-ball is $(-1)^{d-1}$, we can use the theorem to see that 
\[\chi_{\alpha}(T) = \chi_{\alpha}(\P_1)+\chi_{\alpha}(\P_2) + \chi_{\alpha}(\P_3)+\chi_{\alpha}(\P_4)- 4((-1)^{3-1}) = 0.\]

The addition of $\P_2, \P_3$ and $\P_4$ to $\P_1$ is topologically equivalent to handle addition. However, the geometry is fixed by the embedding, so we will call it a \emph{geometric handle addition}.  If $\C'$ is obtained by geometric handle addition to $\C$, a gluing of a complex homeomorphic to a ball with $\C$ along two $(d-1)$-balls, Theorem \ref{Gluing along (d-1)-balls} says that  
\[\chi_{\alpha}(\C') = \chi_{\alpha}(\C) - (-1)^{d-1}.\]  
If we build a polytopal complex by $g$ geometric handle additions to a $3$-polytope, the boundary complex is a surface of genus $g$.  Therefore we can determine the angle characteristic of these complexes constructed by geometric handle additions. 
\begin{corollary} \label{Gram's Relation for surface of genus g}
If $\C$ is a polytopal 3-complex that can be built via geometric handle additions to a $3$-polytope so that the boundary is homeomorphic to a surface of genus $g$, then 
\[\chi_{\alpha}(\C) = 1-g.\]
\end{corollary}
\begin{proof}
Theorem \ref{Gluing along (d-1)-balls} shows that $\chi_{\alpha}(\C)$ changes as $\chi(\partial \C)$ does.  Since a 3-polytope has angle characteristic of 1 and Euler characteristic of 2 on its boundary, the ratio of 1:2 is maintained through construction by geometric handle additions and the end result is independent of the particular construction process since the Euler characteristic is.  Since $\chi(\C) = 2-2g$, $\chi_{\alpha}(\C) = 1-g.$ 
\end{proof}

\begin{figure}
\begin{center}
\includegraphics[width=3in]{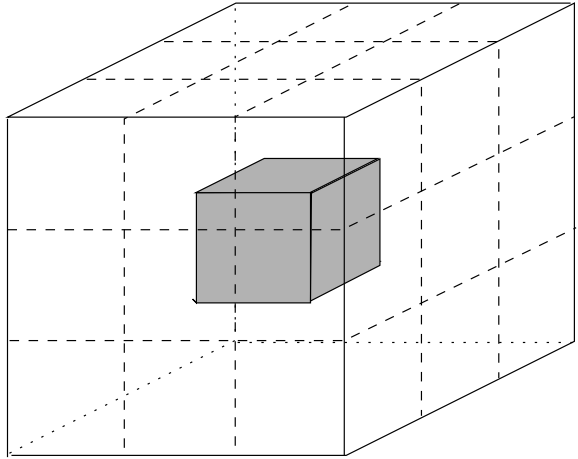}
\caption[The complex $\Gamma$.]{The complex $\Gamma$, a cube of side length 3 with the center unit cube removed.}\label{Cube gamma}
\end{center}
\end{figure}

\begin{figure}
\begin{center}
\includegraphics[width=3in]{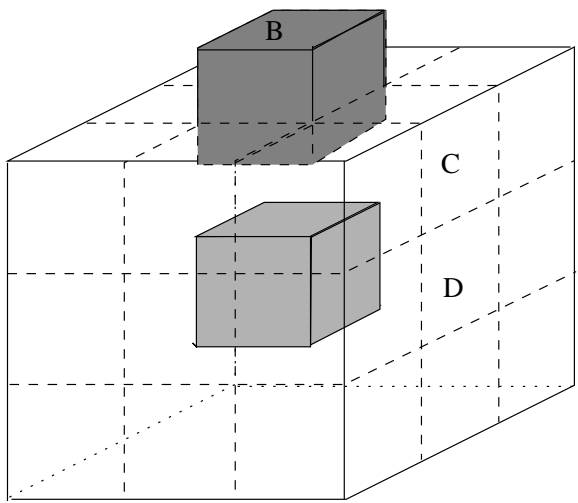}
\caption[An attempted decomposition of $\Gamma$.]{An attempted decomposition of $\Gamma$, trying to order the addition of unit cubes $B$, $C$ and $D$.}\label{cube decomp}
%2nd picture of cube\end{center}
\end{center}
\end{figure}

However, gluing along disjoint complexes homeomorphic to open $(d-1)$-balls will not make all geometric polytopal complexes.  Consider $\Gamma$, a cube of side length 3 built with unit cubes where the center cube is removed, as shown in Figure \ref{Cube gamma}.
\singlespacing
\normalspacing
Suppose these cubes can be ordered so that as each is added it only meets the others along disjoint 2-balls. Consider the unit cubes that are in the center of the faces of the larger cube.  One of these, $B$, must be the last of the center cubes in the order.  Orient the cube so that $B$ is in the middle of the top face. When $B$ is added to the complex, one of the cubes, $C$, that shares a face with $B$ must not yet have been added or else $B$ would meet the previous complex along an annulus. The cube below $C$ in the complex, $D$, cannot be present when $B$ is added or they will meet along an edge, which is not allowed since all maximal faces of the intersection must be two-dimensional.  This is illustrated in Figure \ref{cube decomp}.  
However, this contradicts the choice of $B$ since the cube below $C$ is also the center cube.  
%picture(s) of decomposed cube

\singlespacing
\normalspacing
\singlespacing
\begin{figure}
\begin{center}
\includegraphics[width=3in]{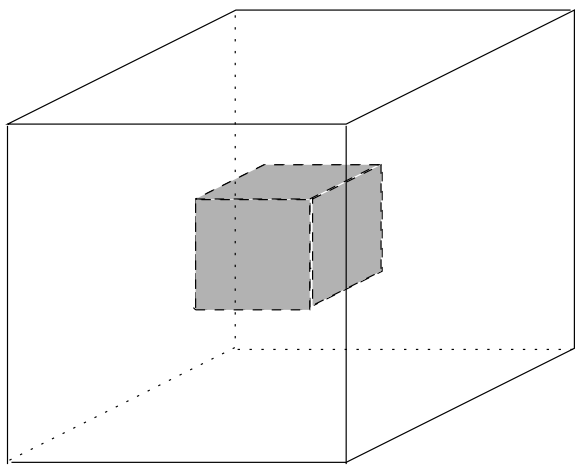}
\caption{$|\Gamma|$, the geometric realization of $\Gamma$.}\label{geo3x3cube}
\end{center}
\end{figure} 
\normalspacing

However, we can compute the Euler and angle characteristics for this complex.  By Lemma \ref{independence of subdivision}, we consider $|\Gamma|$ as shown in Figure \ref{geo3x3cube} and compute
$\chi(\Gamma) = 16-24+12 = 4$ and 
\begin{align*}
\alpha_0\left(\Gamma\right) & = 8\left(\frac{1}{8}\right) + 8\left(\frac{7}{8}\right) = 8,\\
\alpha_1\left(\Gamma\right) & = 12\left(\frac{1}{4}\right) + 12\left(\frac{3}{4}\right) = 12,\\
\alpha_2\left(\Gamma\right) & = 6\left(\frac{1}{2}\right) + 6\left(\frac{1}{2}\right) = 6.\\
\end{align*}
Therefore $\chi_{\alpha}\left(\Gamma\right) = 8-12+6 = 2.$ and we can see that $\chi(\Gamma) = 2 \chi_{\alpha}(\Gamma)$, just as for all odd-dimensional complexes built by gluing polytopes along disjoint sets of balls.  Therefore, we will consider other constructions besides gluing along disjoint $(d-1)$-balls.
\begin{theorem}\label{gluing along an annulus}
If $\C = \A \glue \B$ and $\A \cap \B$ is a union of $m$ disjoint $(d-1)$-dimensional annuli (i.e. closed $(d-1)$-balls with an open $(d-1)$-ball removed from the interior or complexes homeomorphic to $S^{d-2}\times [0,1]$), then 
\[  \chi(\partial \C)= \chi(\partial \A) + \chi(\partial \B) \] 
and
\[\chi_{\alpha} (\C) = \chi_{\alpha}(\A)+\chi_{\alpha}(\B) +  m \left(1 +(-1)^{d} \right) \]
In particular, if $d$ is odd, the Euler characteristic and the angle characteristic of $\C$ is the sum of the respective characteristics of $\A$ and $\B$.  On the other hand, in even dimensions the angle characteristic is greater than the sum of the angle characteristics of $\A$ and $\B$ while the Euler characteristic equals the sum of the Euler characteristics of $\A$ and $\B$.
\end{theorem}
\begin{proof}
For a $(d-1)$-annulus $K$, $\chi(\int (K)) = -1+(-1)^{d-1}$, found by subtracting the Euler characteristics of a closed $(d-1)$-ball and a $(d-2)$-sphere from that of a closed $(d-1)$-ball.  Also, $\chi(\partial K) = 2(1+(-1)^{d-2})$, since the boundary of the annulus is two $(d-2)$-spheres.  Then by Lemma \ref{valuations}
\begin{equation*}
\begin{split}
\chi\left(\partial \C\right) & = \chi\left(\partial \A\right) +\chi\left(\partial \B\right) - 2\chi\left(\int\left(\A \cap \B\right)\right) - \chi\left(\partial\left(\A \cap \B\right)\right)\\
&= \chi\left(\partial \A\right) + \chi\left(\partial \B\right) - 2m\left(-1+(-1)^{d-1}\right) - m\left(2\left(1+(-1)^d\right)\right)\\
& = \chi\left(\partial \A\right) + \chi\left(\partial \B\right)\\
\end{split}
\end{equation*} 
and
\begin{equation*}
\begin{split}
\chi_{\alpha}\left(\C\right) &= \chi_{\alpha}\left(\A\right)+\chi_{\alpha}\left(\B\right) - \chi\left(\int\left(\A \cap \B\right)\right) \\
& = \chi_{\alpha}\left(\A\right)+\chi_{\alpha}\left(\B\right) - m\left(-1+(-1)^{d-1}\right)\\
& = \chi_{\alpha}\left(\A\right)+\chi_{\alpha}\left(\B\right) + m \left(1 +(-1)^{d}\right).
\end{split}
\end{equation*} 
\end{proof}

Since $\Gamma$ can be constructed from polytopes by gluings along balls and annuli, the last two theorems explain why the ratio of 1:2 between the angle and Euler characteristics on polytopes is preserved for $\Gamma$.  We could also construct $\Gamma$ by removing a cube from the center of the solid cube, a construction we consider in the next theorem.

\begin{theorem}\label{gluing along a sphere}
If $\C = \A \glue \B$ and $\A \cap \B$ is a disjoint union of closed geometric polytopal $(d-1)$-complexes without boundary, the difference between the Euler characteristic of $\C$ and the sum of the Euler characteristics of $\A$ and $\B$ is twice the corresponding difference on the angle characteristics. In particular, if $\A \cap \B$ is a union of $m$ disjoint $(d-1)$-spheres, then 
\[  \chi(\partial \C)=\chi(\partial \A) + \chi(\partial \B) - 2m\left(1+(-1)^{d-1}\right) \]
and 
\[\chi_{\alpha} (\C) = \chi_{\alpha}(\A)+\chi_{\alpha}(\B) - m\left(1+(-1)^{d-1}\right). \]
\end{theorem}
\begin{proof}
Since $\A \cap \B$ is a union of closed geometric polytopal complexes, \newline $\int(\A \cap \B) = \A \cap \B$ and $\partial(\A \cap \B) = \emptyset$. Then Lemma \ref{valuations} gives that 
\[\chi(\partial \C) = \chi(\partial \A) +\chi(\partial \B) - 2\chi\left(\int(\A \cap \B)\right)\]
and
\[\chi_{\alpha}(\C) = \chi_{\alpha}(\A)+\chi_{\alpha}(\B) - \chi\left(\int(\A \cap \B)\right).\]
Therefore, the difference between the Euler characteristics is double the corresponding difference between the angle characteristics.

If $K$ is a $(d-1)$-sphere, then $\chi(\int (K)) = \chi(K) = 1+(-1)^{d-1}$, so   
\[\chi(\partial \C)= \chi(\partial \A) + \chi(\partial \B) - 2m\left(1+(-1)^{d-1}\right)\]
and
\[\chi_{\alpha}(\C) = \chi_{\alpha}(\A)+\chi_{\alpha}(\B) - m\left(1+(-1)^{d-1}\right).\]
\end{proof}

We define the set of semi-constructible complexes as the set of complexes that can be made from polytopes by iteratively gluing along a disjoint union of $(d-1)$-balls, $(d-1)$-annuli, or $(d-1)$-complexes without boundary.  Therefore, the next theorem follows immediately from the previous theorems.
\begin{theorem}\label{angle char half of Euler char}
If $\C$ is a $d$-dimensional \compname  geometric polytopal complex where $d$ is odd, then $\chi_{\alpha}(\C) = \frac{1}{2}\chi(\partial\C)$.
\end{theorem}

We conclude our set of constructions with one last case of gluings, that where $\A \cap \B$ is of lower dimension.
\begin{theorem}\label{gluing along a lower dimensional complex}
Let $\C = \A \glue \B$ where $\A \cap \B$ has dimension less than $d-1$.  Then:
\[ \chi(\partial \C) = \chi(\partial \A) + \chi(\partial \B)  - \chi(\partial(\A \cap \B))\] 
and
\[\chi_{\alpha} (\C) = \chi_{\alpha}(\A)+\chi_{\alpha}(\B).\]
\end{theorem}
\begin{proof}
Since $\A \cap \B$ has dimension less than $d-1$, $\int\left(\A \cap \B\right)=\emptyset$ and \newline $\partial(\A \cap \B) = \A \cap \B$.  Then Lemma \ref{valuations} gives that 
\begin{align*}
\chi(\partial \C) & = \chi(\partial \A) +\chi(\partial \B) - 2\chi\left(\int(\A \cap \B)\right) - \chi\left(\partial(\A \cap \B)\right)\\
&= \chi(\partial \A) + \chi(\partial \B)  - \chi\left(\partial(\A \cap \B)\right) 
\end{align*}
and
\begin{align*}
\chi_{\alpha}(\C) &= \chi_{\alpha}(\A)+\chi_{\alpha}(\B) - \chi\left(\int(\A \cap \B)\right) \\
& = \chi_{\alpha}(\A)+\chi_{\alpha}(\B).
\end{align*}
\end{proof}

Therefore, whenever the angle and Euler characteristics or $\A$ and $\B$ have a 1:2 ratio and $\chi\left(\partial(\A \cap \B)\right) = 0$, $\chi_{\alpha}(\partial\C)$ and $\chi(\C)$ will maintain this ratio. In particular, if $\partial(\A \cap \B)$ is an odd-dimensional semi-Eulerian manifold, we have that the ratio is maintained since $\chi\left(\partial(\A \cap \B)\right)= 0$ as a result of Theorem \ref{Klee DS}.

\section{Complexes that are \compname}
%Furch's knotted hole ball example
For all odd-dimensional \compname  complexes, we can see that there is a 1:2 ratio between the angle characteristic and the Euler characteristic.  We can ask whether all polytopal complexes can be constructed in this manner, and if not, whether these other polytopal complexes still maintain the 1:2 ratio between the angle characteristic and the Euler characteristic.

%To intro??
First, we compare to a few classes of complexes already in the literature. In all cases we will focus on those $d$-complexes which can be PL-embedded in $\R^d$ and are \emph{pseudomanifolds}, that is, $d$-dimensional pure simplicial complexes in which each $(d-1)$-dimensional face belongs to at most two facets \cite{Ha}.  By the way we have constructed our geometric polytopal complexes, each is a pseudomanifold.

 A pure $d$-dimensional complex is \emph{shellable} if its facets can be ordered $F_1, \ldots, F_t$ so that $\displaystyle{\left(\bigcup_{i=1}^{j-1} F_i\right) \cap F_j}$ is a pure $(d-1)$-complex which has a shelling that extends to all of $\partial F_j$ for $2 \leq j \leq t$.  A shellable pseudomanifold is also \compname since as each facet is added it is added along a $(d-1)$-ball or sphere, a subset of the boundary of the polytope.  Bruggesser and Mani showed that all boundary complexes of polytopes are shellable \cite{BM}.  However, there are balls and spheres for all $d\geq3$ that are not shellable \cite{Ha}.

\emph{Constructible} complexes satisfy the following recursive definition:
\begin{itemize}
\item Every simplex is constructible.
\item A $d$-complex which is not a simplex is constructible if and only if it can be written as $\C=\A \cup \B$ where $\A$ and $\B$ are constructible $d$-complexes and $\A \cap \B$ is a constructible $(d-1)$-complex.
\end{itemize}
This construction is a gluing along a $(d-1)$-ball. Therefore, since simplices are semi-constructible, all constructible pseudomanifolds are \compname complexes.  However, for a constructible complex $\A \cap \B$ must be connected, which is not the case when constructing a \compname  complex. As a result, the set of \compname  complexes is larger than the set of embedded constructible complexes.  Even though every shellable complex must also be constructible, there are balls and spheres that are not constructible.

Shellable and constructible pseudomanifolds are all homeomorphic to balls or spheres \cite{Ha}, and our constructions clearly make boundary complexes with more complex topology.  Since semi-Eulerian complexes have a much wider variety of topology, we also compare the sets of semi-Eulerian complexes and \compname  complexes.  Any semi-Eulerian complex that cannot be PL-embedded in $\R^d$ cannot be a \compname complex.  There are many such complexes, such as the projective plane.  On the other hand, the boundaries of polytopes are semi-Eulerian manifolds and and gluings along subcomplexes of semi-Eulerian manifolds homeomorphic to $(d-1)$-balls and $(d-1)$-annuli will maintain the proper Euler characteristic on the link of each face, so the result is also a semi-Eulerian manifold.  The same will be true if a complex is formed by a gluing along a semi-Eulerian manifold, but not necessarily along a closed manifold.  This is true even though the operation of geometric handle addition seems to be less restrictive than that of handle addition. The former can be accomplished between two facets with different combinatorics, and vertices of the facets which would be identified in a handle addition can be connected by an edge; the geometry guarantees that it is still a geometric polytopal complex since we are adding extra edges along the handle.  Also, we can do geometric handle addition along complexes and not just along a single facet.  In this case, the geometry may limit facets along which handles can be added.  Still, as long as we are gluing along $(d-1)$-balls, $(d-1)$-annuli and semi-Eulerian $d$-complexes, the result is a semi-Eulerian complex.

\begin{figure}
\begin{center}
\includegraphics[width=3.5in, height=4in, angle=270]{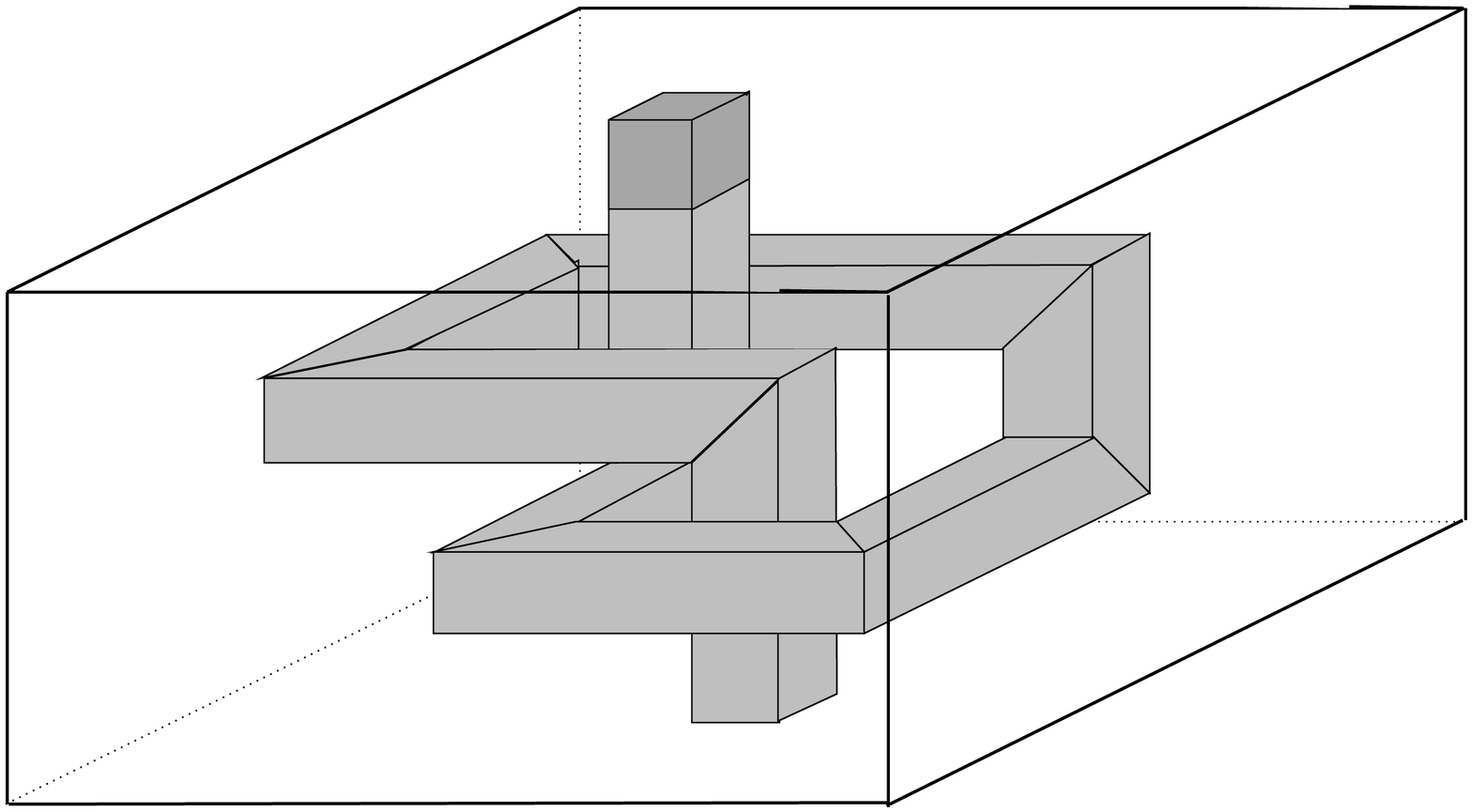}
\caption{Furch's knotted hole ball.}\label{Furch}
\end{center}
\end{figure}

Even though the class of \compname  complexes is larger than the classes of constructible or shellable pseudomanifolds, we can consider complexes which are not constructible as candidates for a non-\compname  complex.  One of the common ways to create a non-constructible, and hence non-shellable, complex is to introduce a knot into the complex.  Therefore, we consider Furch's knotted hole ball as shown in Figure \ref{Furch}.
This complex is created by starting with a cube and hollowing out a knotted tunnel through the cube.  This tunnel is then ``plugged'' with the darker cube in the figure. This complex is homeomorphic to a ball, and is neither shellable or constructible \cite{Ha}.   The fact that the complex is not constructible says that it cannot be constructed by gluings, each of which is along a single 2-ball.  Therefore Corollary \ref{Gram's Relation for special spherical  polytopal complexes} does not apply to $\FF$ to give us the anticipated relationship between the Euler and angle characteristics.%We can in fact show a stronger statement.
\singlespacing
\normalspacing

\begin{figure} \begin{center} 
\includegraphics[width=3.5in, height=4in, 
angle=270]{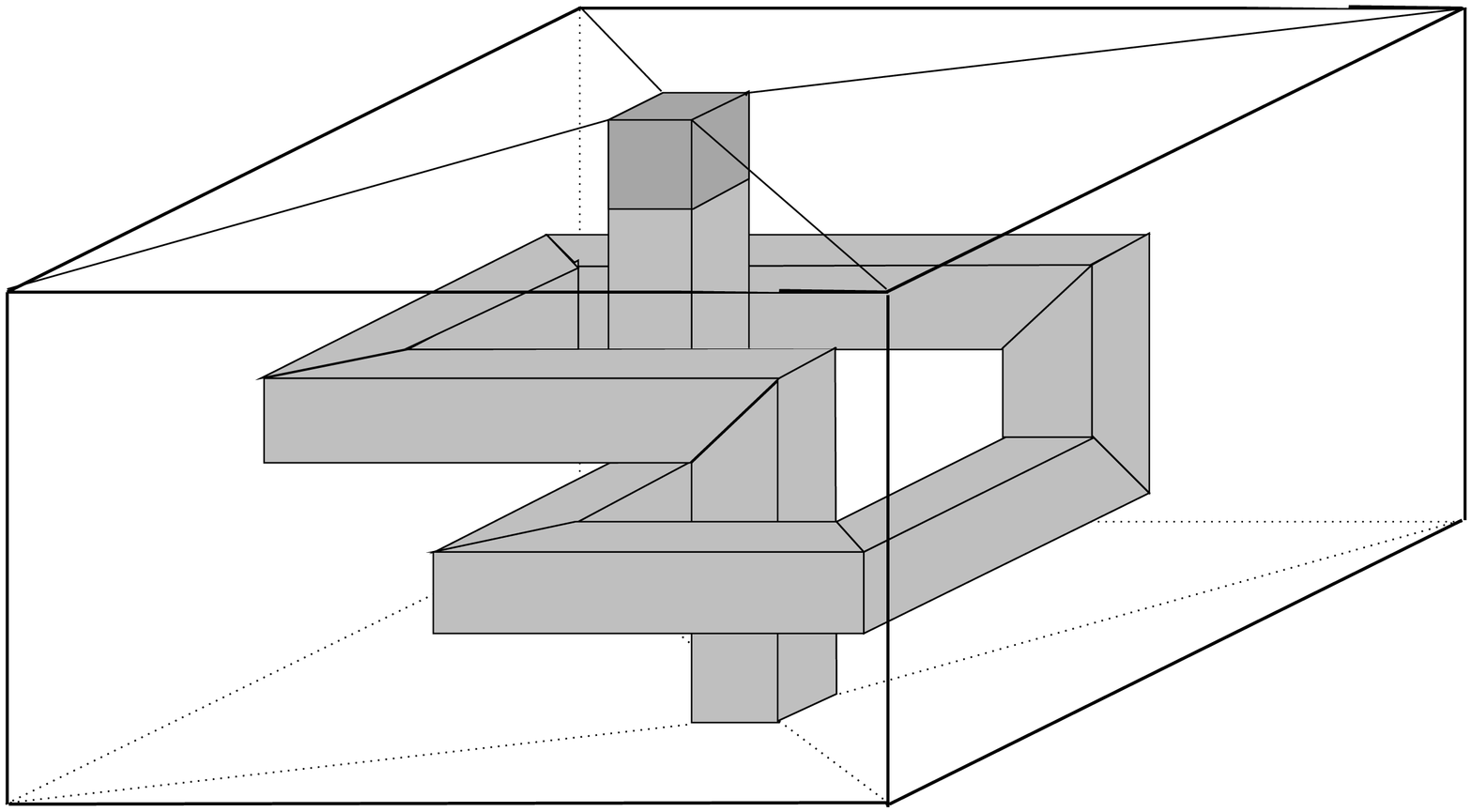} \caption{Furch's knotted hole 
ball, subdivided.}\label{Furchsubdiv} \end{center} \end{figure}

However, we can compute that Furch's knotted hole ball, which we denote 
$\FF$ still satisfies the Gram relation.  We work with the subdivision 
shown in Figure \ref{Furchsubdiv} where two of the outer faces are 
subdivided. \singlespacing  
\normalspacing Then the angle sums are: \begin{equation*} \begin{split} 
\alpha_{0}(\FF) & = 8\left(\frac{1}{8}\right) + 4\left(\frac{1}{2}\right) 
+ 4\left(\frac{3}{8}\right)+ 8\left[2\left(\frac{7}{8}\right)+ 
2\left(\frac{5}{8}\right)\right] + 4\left(\frac{7}{8}\right) = 32,\\ 
\alpha_{1}(\FF) & = 
\left[12\left(\frac{1}{4}\right)+8\left(\frac{1}{2}\right)\right] + 
4\left(\frac{1}{2}\right) + 4\left(\frac{1}{4}\right)\\ & \qquad \qquad + 
8\left[4\left(\frac{3}{4}\right)+ 2\left(\frac{1}{2}\right)+ \frac{3}{4} + 
\frac{1}{4}\right] + 
\left[4\left(\frac{3}{4}\right)+4\left(\frac{3}{4}\right)\right] = 56, \\
\alpha_{2}(\FF)  & =
12\left(\frac{1}{2}\right) + \frac{1}{2} + 0 +
8\left[4\left(\frac{1}{2}\right)\right] + 5\left(\frac{1}{2}\right) =
25,
\end{split} \end{equation*} 
where the first term of each refers to faces on the exterior of the 
cube, the second to the top of the plug cube, the third to the entrance 
into the knot, the fourth to each bend in the knot (including edges on the 
side furthest from the entrance to the knot, but not those closest to the 
entrance), and finally to the last section of the knot ending with the 
plug cube. Then we can see that \[\chi_{\alpha}(\FF) = 32-56+25 = 1.\] 
This is the same as the angle characteristic of a polytope, so this 
satisfies the Gram relation.  Since $\FF$ is a 3-ball, $\chi(\partial\FF) 
= 2$ and $\chi(\partial\FF) = 2\chi_{\alpha}(\FF)$.

This can also be shown by using Theorem \ref{Gluing along (d-1)-balls} to compute the angle characteristic of one of the component complexes in a gluing.  If we consider the complex $\FF'$ formed by placing a cube $C$ into the knot, pushing it in until it shares a face with the plug cube and starts to plug the hole, we can see that 
\[\chi_{\alpha}(\FF') = \chi_{\alpha}(\FF) + \chi_{\alpha}(C) - \left((-1)^{3-1}\right) = \chi_{\alpha}(\FF)\]
since the cube meets $\FF$ along five of its faces, homeomorphic to a 2-ball, and $\chi_{\alpha}(C) = 1$.  In the same way, we can fill in the hole with cubes until we have a standard cube, which the computation shows to have the same angle characteristic as Furch's knotted hole ball.  

Therefore, although we cannot build this ball via gluings along one 2-ball at a time, %$\FF$ is not \compname, 
the theorems in the previous section still allow us to conclude that its angle characteristic is half the Euler characteristic of its boundary. All other examples of geometric polytopal complexes we know which are counterexamples to constructibility are either \compname  or are component complexes in a gluing along $(d-1)$-balls, $(d-1)$-annuli, or closed $(d-1)$-complexes that gives a polytope as above.  %DOUBLE CHECK! 

%D-S and Perles variants
\section{Perles Relations on Complexes}
A similar extension of the Perles equations can be determined for \emph{geometric simplicial complexes}, geometric polytopal complexes whose boundary consists of simplices.  Interior polytopes of such a complex can be subdivided so that all the maximal polytopes are simplicial, so we assume our complex is composed of $d$-simplices.    We will consider how different gluings affect the Dehn-Sommerville and Perles relations.  Therefore, we define the following operators on a $d$-dimensional geometric polytopal complex $\C$:
\begin{equation}
DS_{k}(\C) \equiv \sum_{j=k}^{d} (-1)^i {j+1\choose k+1} f_{j}(\C)
\end{equation}
and
\begin{equation}
Pe_{k}(\C) \equiv \sum_{j=k}^{d-1} (-1)^i {j+1\choose k+1} \alpha_{j}(\C).
\end{equation}
Then, for a $d$-polytope $\P$, we can rewrite the Perles relations as  \[Pe_{k}(\P) = (-1)^d\left[\alpha_{k}(\P) - f_{k}(\P)\right] \quad \text{for } 0 \leq k \leq d-2\] and the Dehn-Sommerville relations as \[DS_{k}(\partial\P) = (-1)^{d-1}f_{k}(\P) \quad \text{for } 0 \leq k \leq d-2.\] $Pe_{k}$ will be called the Perles operator and $DS_{k}$ the Dehn-Sommerville operator.  Usually, we will use $DS_{k}$ to act on the boundary of a $d$-complex $\C$ so that we sum over the same faces that we do to determine the angle sums. That is, \[DS_{k}(\partial\C)  = \sum_{j=k}^{d-1} (-1)^i {j+1\choose k+1} f_{j}(\partial\C).\]
If $\C$ is of dimension $l<d$, 
\[DS_{k}(\partial\C) \equiv \sum_{j=k}^{l-1} (-1)^i {j+1\choose k+1}f_{j}(\partial\C) = \sum_{j=k}^{d-1} (-1)^i {j+1\choose k+1} f_{j}(\partial\C)\]
since $f_{j}(\partial\C)=0$ if $j\geq l$.  Therefore we will leave the dimension implicit in naming $DS_k$, even though it has some impact.

Theorem \ref{Klee DS} tells us that for $0 \leq k \leq d-2$ the 
Dehn-Sommerville relations hold on all semi-Eulerian $(d-1)$-complexes.  
However, we have seen that the set of \compname complexes is not 
easily compared to the set of semi-Eulerian complexes.  Therefore, we will 
still consider the effects on the Dehn-Sommerville and Perles 
operators in parallel, occasionally using the Dehn-Sommerville relations 
for semi-Eulerian complexes.  We start by considering a generic gluing of 
two $d$-complexes, $\A$ and $\B$, where all the $d$-polytopes are 
simplicial.  \begin{lemma}\label{simplicial valuations} Let $\C = \A \glue 
\B$ where $\A$ and $\B$ are pure simplicial polytopal $d$-complexes in 
$\R^d$.  Then for $0 \leq k \leq d-1$ \[DS_{k}(\partial\C) = 
DS_{k}(\partial \A) + DS_{k}(\partial \B) - 2 DS_{k}\left(\int(\A \cap 
\B)\right) - DS_{k}\left(\partial(\A \cap \B)\right)\] and \[Pe_{k}(\C) = 
Pe_{k}(\A) + Pe_{k}(\B) - DS_{k}\left(\int(\A \cap \B)\right)\] 
\end{lemma} \begin{proof} Using \eqref{f-valuations} and 
\eqref{a-valuations}, for $0 
\leq i \leq d-1$ \begin{align*} DS_{k}(\partial\C) & \equiv 
\sum_{j=k}^{d-1} (-1)^j {j+1 \choose k+1} f_{j}\left(\partial\C\right)\\ 
&= \sum_{j=k}^{d-1} (-1)^j {j+1 \choose k+1} \left[f_j(\partial\A) + 
f_j(\partial \B)\right. \\ & \qquad \qquad \left.- 2f_j\left(\int(\A \cap 
\B)\right) - f_j\left(\partial(\A \cap \B)\right)\right] \\ & = 
DS_{k}(\partial \A) + DS_{k}(\partial \B) - 2 DS_{k}\left(\int(\A \cap 
\B)\right) - DS_{k}\left(\partial(\A \cap \B)\right) \end{align*} and 
\begin{align*} Pe_{k}(\C) & \equiv \sum_{j=k}^{d-1} (-1)^j {j+1 \choose 
k+1} \alpha_{j}(\A \glue \B)\\ & = \sum_{j=k}^{d-1} (-1)^j {j+1 \choose 
k+1} \left[\alpha_i(\A) + \alpha_i(\B) - f_i\left(\int(\A 
\cap\B)\right)\right]\\ & = Pe_{k}(\A) + Pe_{k}(\A) - DS_{k}\left(\int(\A 
\cap \B)\right). \end{align*} \end{proof}

\begin{figure}
\begin{center}
\includegraphics[width=4.5in]{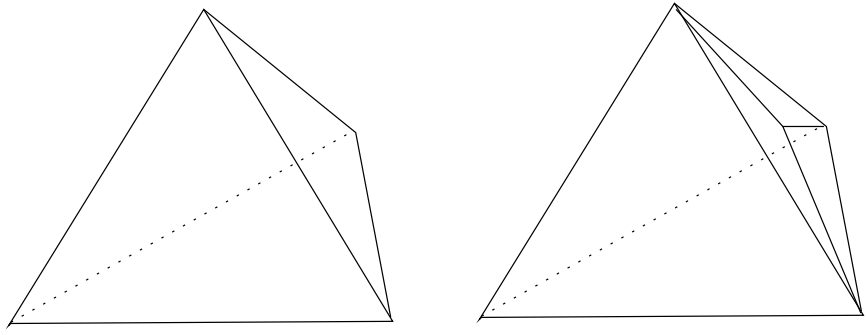}
\caption[A tetrahedron and stellar subdivision of one of its faces.]{The regular tetrahedron $\D$ and $\D'$, the stellar subdivision of one face of the tetrahedron.}\label{D and D'}
\end{center}
\end{figure}

Since the simplicial subdivision of a simplicial polytope also results in a simplicial polytope, we know that the Dehn-Sommerville and Perles relations will still hold under simplicial subdivision of a polytope. Unlike the Euler and angle characteristics, we would expect the values of the operators to change even though the relations still hold since these relations include the number of faces in a given dimension. For example, we can consider the regular tetrahedron $\D$ and the stellar subdivision of one face of the tetrahedron $\D'$, as shown in Figure \ref{D and D'}.  We have that 
\[\alpha\text{-}f(\D) = \left(\frac{3}{\pi}\arccos\left(\frac{1}{3}\right)-1, \frac{3}{\pi}\arccos\left(\frac{1}{3}\right),2,1, 4, 6, 4, 1\right)\]
and
\[\alpha\text{-}f(\D') = \left(\frac{3}{\pi}\arccos\left(\frac{1}{3}\right)-\frac{1}{2}, \frac{3}{\pi}\arccos\left(\frac{1}{3}\right)+\frac{3}{2},3,1, 5, 9, 6, 1\right).\]
Then we can compute 
$DS_{1}(\D) =  6$, $DS_{1}(\D') = 9$, 
$Pe_{1}(\D) = 6-\frac{3}{\pi}\arccos\left(\frac{1}{3}\right)$, and $Pe_{1}(\D') = 9 - \left( \frac{3}{\pi}\arccos\left(\frac{1}{3}\right)+\frac{3}{2}\right)$, giving different answers on different subdivisions as expected.
\singlespacing
\normalspacing

The following lemmas will be useful in computing the Dehn-Sommerville and Perles operators for certain classes of gluings.
%MCMULLEN & WALKUP - relating f(?)-vector of whole to that of boundary

\begin{lemma}\label{DS on open balls}
Let $K$ be a simplicial complex homeomorphic to a $(d-1)$-ball.  Then   \[DS_{k}\left(\int(K)\right) = (-1)^{d-1}f_k(K)\] and \[DS_{k}(K) = (-1)^{d-1}f_k\left(\int (K)\right).\]
\end{lemma}
\begin{proof}
Let $K^{*}$ be the complex formed by identifying two copies of $K$ along corresponding faces of the boundary.  Then $K^{*}$ is a $(d-1)$-sphere and the faces in $K^{*}$ are the disjoint union of faces of two copies of $\int(K)$ and one copy of $\partial K$, which is homeomorphic to a $(d-2)$-sphere.  So $f_i(K^{*}) = 2f_i\left(\int(K)\right) + f_i(\partial K)$. Then 
\begin{equation*}
\begin{split}
(-1)^{d-1}f_k\left(K^{*}\right) & = DS_{k}\left(K^{*}\right) \\
&= \sum_{j=k}^{d-1} (-1)^j {j+1 \choose k+1} [2f_{j}(\int(K))+f_{j}(\partial K)]\\
&=  2DS_{k}(\int(K))+ DS_{k}(\partial K)\\
&= 2DS_{k}(\int(K)) + (-1)^{d-2}f_k(\partial K),
\end{split}
\end{equation*}
where the first and last equalities follow from the Dehn-Sommerville relations on spheres. %(Note that there are no $(d-1)$-faces in $\partial K$.} 

\noindent We rearrange to get:
\begin{equation*}
\begin{split}
DS_{k}\left(\int(K)\right) & = \frac{1}{2}(-1)^{d-1}\left[f_k\left(K^{*}\right)  + f_k(\partial K)\right]\\
& = \frac{1}{2}(-1)^{d-1}\left[\left(2f_k(\int(K)) + f_k(\partial K)\right)+ f_k(\partial K)\right] \\
& = (-1)^{d-1}\left(f_k(\int(K)) + f_k(\partial K)\right)\\
& = (-1)^{d-1}f_k(K).
\end{split}
\end{equation*}
This gives the relation for the open ball $\int(K)$.

Similarly, we can get $DS_{k}(K)$ by writing 
\begin{equation*}
\begin{split}
DS_{k}(K) & = DS_{k}\left(K^{*}\right) - DS_{k}\left(\int(K)\right)\\
&  = (-1)^{d-1}f_k\left(K^{*}\right) - (-1)^{d-1}f_k(K) \\
& = (-1)^{d-1}f_k\left(\int (K)\right).
\end{split}
\end{equation*}
\end{proof}

\begin{lemma}\label{DS on annuli}
Let $K$ be a simplicial complex homeomorphic to a $(d-1)$-annulus.  Then   \[DS_{k}(\int(K)) = (-1)^{d-1}f_k(K)\] and \[DS_{k}(K) = (-1)^{d-1}f_k\left(\int K\right).\]
In addition, if $K^{*}$ is the complex formed by two copies of the annulus $K$ identified along corresponding faces of $\partial K$,
\[DS_{k}(K^{*}) = (-1)^{d-1}f_k(K^{*}).\]
\end{lemma}
\begin{proof}
If $K$ is a $(d-1)$-annulus, we can extend it to a $(d-1)$-ball by adding a vertex $v$ and faces which are the convex hulls of $v$ and the $(d-2)$-faces on the inner boundary of the annulus.  Denote the new faces by $K'$, and the complex induced by these faces  by $\cl(K')$, the closure of $K'$.  Then $K'$ is homeomorphic to an open $(d-1)$-ball and $\cl(K')$ to a closed $(d-1)$-ball.  Then, using Lemma \ref{DS on open balls}, we can write 
\begin{equation*}
\begin{split}
DS_{k}(K) & = DS_{k}\left(K \cup K'\right) - DS_{k}(K') \\
&= (-1)^{d-1}f_{k}\left(\int\left(K \cup K'\right)\right)- (-1)^{d-1}f_{k}\left(\cl(K')\right)\\
& = (-1)^{d-1}f_{k}\left(\int(K)\right).
\end{split}
\end{equation*}
Likewise, since the boundary of an annulus is two $(d-2)$-spheres, \begin{equation*}
\begin{split}
DS_{k}\left(\int(K)\right) & = DS_{k}(K) - DS_{k}(\partial K) \\
& = (-1)^{d-1}f_{k}\left(\int(K)\right) - (-1)^{d-2}f_{k}(\partial K) \\
& = (-1)^{d-1}f_{k}(K).
\end{split}
\end{equation*}

When we consider the complex $K^{*}$, formed by two copies of the annulus $K$ identified along their boundary we see that
\begin{equation*}
\begin{split}
DS_{k}\left(K^{*}\right) & = DS_{k}\left(\int(K)\right) + DS_{k}(K) \\
& = (-1)^{d-1}f_{k}(K) + (-1)^{d-1}f_{k}\left(\int(K)\right)\\
&  = (-1)^{d-1}f_{k}\left(K^{*}\right).
\end{split}
\end{equation*}
\end{proof}

Now we use Lemmas \ref{DS on open balls} and \ref{DS on annuli} to consider $DS_k(\C)$ and $Pe_k(\C)$ on complexes constructed by gluing along $m$ disjoint $(d-1)$-balls and $(d-1)$-annuli.  
\begin{theorem}\label{DS on ball and annuli gluings}
Let $\C = \A \glue \B$, where $\A \cap \B$ is the union of $m$ disjoint $(d-1)$-balls or $(d-1)$-annuli, $K_1, K_2, \ldots, K_m$.  For $i=1,2,\ldots,m$, define $K_{i}^{*}$ as the complex made by identifying two copies of $K_{i}$ along the corresponding faces of the boundary.  Then
\[DS_{k}(\partial(\C)) = DS_{k}(\partial\A) + DS_{k}(\partial\B) - \sum_{i=0}^{m} (-1)^{d-1}f_{k}\left(K_{i}^{*}\right).\]
As a result, the Dehn-Sommerville and Perles relations for \mbox{$0 \leq k \leq d-2$} hold on $\C$ if the relations hold on $\A$ and $\B$.  %That is 
%\[DS_{k}(\partial\C) = (-1)^{d-1}f_k(\C)\]
%and 
%\[Pe_{k} (\C) = (-1)^{d}\left[\alpha_{k}(\C) - f_k(\C)\right].\]
\end{theorem}
\begin{proof}
The terms $ 2 DS_{k}\left(\int(\A \cap \B)\right) + DS_{k}\left(\partial(\A \cap \B)\right)$ from Lemma \ref{simplicial valuations} count the contribution of the Dehn-Sommerville operator on two copies of the interior of the intersection and one copy of the boundary.  Then if $K^{*}$ is formed by two copies of $K$ identified along corresponding parts of their boundary, these terms give $\displaystyle{\sum_{i=0}^{m} DS_{k}\left(K_{i}^{*}\right)}$.  Therefore we use Lemmas \ref{DS on open balls} and \ref{DS on annuli} and the fact that if $K$ is a $d$-ball, $K^{*}$ is a $d$-sphere to get
\begin{equation*}
\begin{split}
DS_{k}(\partial\C) & = DS_{k}(\partial\A) + DS_{k}(\partial\B) - \sum_{i=0}^{m} DS_{k}\left(K_{i}^{*}\right)\\
& = DS_{k}(\partial\A) + DS_{k}(\partial\B) - \sum_{i=0}^{m} (-1)^{d-1}f_{k}(K_{i}^{*}).
\end{split}
\end{equation*}
If the Dehn-Sommerville relations hold for $\A$ and $\B$, then the above simplifies:
\begin{equation*}
\begin{split}
DS_{k}(\partial\A) + DS_{k}&(\partial\B)  - \sum_{i=0}^{m} (-1)^{d-1}f_{k}\left(K_{i}^{*}\right)\\
& = (-1)^{d-1}f_k(\partial\A) +(-1)^{d-1}f_k(\partial\B) - \sum_{i=0}^{m}(-1)^{d-1}f_k\left(K_{i}^{*}\right) \\
& = (-1)^{d-1}f_k(\partial\C).
\end{split}
\end{equation*}

Similarly, if we assume $\A$ and $\B$ satisfy the Perles relations we have the following:
\begin{align*}
Pe_{k} (\C) & = Pe_{k}(\A) + Pe_{k}(\B) - \sum_{i=1}^{m} DS_{k}\left(\int(K_i)\right)\\
&= (-1)^d\left[\alpha_{k}(\A) - f_k(\partial\A)\right] +(-1)^d\left[\alpha_{k}(\B) - f_k(\partial\B)\right]  \\
&\qquad \qquad -\sum_{i=1}^{m} (-1)^{d-1}f_{k}(K_{i}) \quad \text{by Lemmas } \ref{DS on open balls} \text{ and } \ref{DS on annuli}\\
&= (-1)^d\left[\left(\alpha_{k}(\A) + \alpha_{k}(\B) - \sum_{l=1}^{m} f_{k}\left(\int(K_{i})\right)\right)  \right.\\
& \qquad \qquad \qquad \qquad \left. 
- \left(f_k(\A)+f_k(\B)- \sum_{l=1}^{m}f_k\left(K_{i}^{*}\right)\right) \right] \\
&=(-1)^{d}\left[\alpha_{k}(\C) - f_k(\partial\C)\right],
\end{align*}
where the last equality follows from \eqref{f-valuations} and \eqref{a-valuations}.
\end{proof}

It is not surprising that the Dehn-Sommerville relations hold on these complexes, because all the complexes made are semi-Eulerian, so Theorem \ref{Klee DS} guarantees the Dehn-Sommerville relations.  

As with the angle characteristic, we can also consider the result if complexes are glued along $(d-1)$-complexes without boundary.
%DS and P for sphere identifications
\begin{theorem}\label{DS on sphere-gluings}
Let $\C = \A \glue \B$, where $\A \cap \B$ is the union of $m$ disjoint semi-Eulerian $(d-1)$-complexes, $K_1, K_2, \ldots K_m$. Then the Dehn-Sommerville and Perles relations on $\C$ for $0 \leq k \leq d-2$ hold if the relations hold on $\A$ and $\B$.  That is 
\[DS_{k}(\partial\C) = (-1)^{d-1}f_k(\partial\C)\]
and 
\[Pe_{k} (\C) = (-1)^{d}(\alpha_{k}(\C) - f_k(\partial\C)).\]
\end{theorem}
\begin{proof}
Since $\int(\A \cap \B)) = \A \cap \B$ and $\partial(\A \cap \B)=\emptyset$, Lemma \ref{simplicial valuations} gives that
\[DS_{k}(\partial\C) =  DS_{k}(\partial \A) + DS_{k}(\partial \B) - 2 DS_{k}(\A \cap \B)\]
and
\[Pe_{k}(\C) = Pe_{k}(\A) + Pe_{k}(\B) - DS_{k}(\A \cap \B).\]
$\A \cap \B$ is the union of $m$ semi-Eulerian $(d-1)$-complexes, so $DS_{k}(\A \cap \B) = (-1)^{d-1}f_k(\A \cap \B)$ by Theorem \ref{Klee DS}.
Therefore,
\[DS_{k}(\partial\C) = DS_{k}(\partial \A) + DS_{k}(\partial \B) - 2(-1)^{d-1}f_k(\A \cap \B)\]
and
\[Pe_{k}(\C) = Pe_{k}(\A) + Pe_{k}(\B) - (-1)^{d-1}f_k(\A \cap \B).\]

If the Dehn-Sommerville relations hold for $\A$ and $\B$, then the above simplifies to
\begin{align*}
DS_{k}(\partial\C) & = (-1)^{d-1}f_k(\partial\A) +(-1)^{d-1}f_k(\partial\B) - 2(-1)^{d-1}f_k(\A \cap \B) \\
& = (-1)^{d-1}f_k(\C).
\end{align*}
Similarly, if we assume $\A$ and $\B$ satisfy the Perles relations we have that:
\begin{align*}
Pe_{k} (\C) & = Pe_{k}(\A) + Pe_{k}(\B) - (-1)^{d-1}f_k(\A \cap \B)\\
&= (-1)^d[\alpha_{k}(\A) - f_k(\partial\A)] +(-1)^d[\alpha_{k}(\B) - f_k(\partial\B)]  +(-1)^{d}f_k(\A \cap \B) \\
&= (-1)^d\left[\left(\alpha_{k}(\A) + \alpha_{k}(\B) - (-1)^{d}f_k(\A \cap \B)\right)\right.  \\
& \qquad \qquad - \left.\left(f_k(\partial\A)+f_k(\partial\B)- 2(-1)^{d}f_k(\A \cap \B)\right) \right] \\
&=(-1)^{d}[\alpha_{k}(\C) - f_k(\partial\C)],
\end{align*}
where the last equality follows from \eqref{f-valuations} and \eqref{a-valuations}.
\end{proof}

The proof of this last theorem shows that even if the Dehn-Sommerville and Perles relations do not hold on $\A$ and $\B$, the difference between the Dehn-Sommerville operator on $\C$ and the sum of the Dehn-Sommerville operators on $\A$ and $\B$ is double the corresponding difference of Perles operators if $\A\cap \B$ is a disjoint union of semi-Eulerian $(d-1)$-complexes.  This is reminiscent of what we were finding with the Euler and angle characteristics.

The last proof would also apply for intersections which were semi-Eulerian manifolds of dimension $l$, where $l$ has the same parity as $d-1$.  This follows since the only effect this would have is in using the Dehn-Sommerville relations, which give the same result as long as the dimensions have the same parity. Of course, if $l < k$, $f_k(\A \cap \B) = 0$.  We now consider other intersections along lower dimensional complexes.
%DS and P for lower dimensional intersections
\begin{theorem}\label{DS on lower dimensional intersections}
Let $\C = \A \glue \B$, where $\A \cap \B$ is of dimension $l \leq d-2$ and suppose the Dehn-Sommerville relations hold for $\A$ and $\B$. Then the following relations hold on $\A$ and $\B$:  
\begin{equation*}
\begin{split}
DS_{k}(\partial\C) = (-1)^{d-1}f_k(\partial\A)& +(-1)^{d-1}f_k(\partial\B) - (-1)^{k}f_k(\A \cap\B)\\ - &\sum_{j=k+1}^{d-1} (-1)^j {j+1 \choose k+1} f_{j}(\A \cap \B)
\end{split}
\end{equation*}
and
\[Pe_{k} (\C) =(-1)^{d}[\alpha_{k}(\C) - (f_k(\partial\C)+f_k(\A \cap \B))].\]
In particular, the Dehn-Sommerville and Perles relations hold for $l < k \leq d-2$ and also for $l = k$ if $k$ is the same parity as $d-1$.
\end{theorem}
\begin{proof}
Since $\A \cap \B$ is of dimension $l \leq d-2$, $\int(\A \cap \B)= \emptyset$ and $\partial(\A \cap \B) = \A \cap \B$.  Therefore Lemma \ref{simplicial valuations} gives us that \[DS_{k}(\partial\C) =  DS_{k}(\partial \A) + DS_{k}(\partial \B) - DS_{k}(\A \cap \B)\]
and
\[Pe_{k}(\C) = Pe_{k}(\A) + Pe_{k}(\B).\]
If the Dehn-Sommerville relations hold for $\A$ and $\B$, then
\begin{align*}
DS_{k}(\partial\C) & = DS_{k}(\partial\A) + DS_{k}(\partial\B) - DS_{k}(\A \cap \B)\\
& = (-1)^{d-1}f_k(\partial\A) +(-1)^{d-1}f_k(\partial\B) - \sum_{j=k}^{d-1} (-1)^j {j+1 \choose k+1} f_{j}(\A \cap \B)\\
& = (-1)^{d-1}f_k(\partial\A) +(-1)^{d-1}f_k(\partial\B) - (-1)^kf_k(\A \cap\B) \\
& \qquad \qquad - \sum_{j=k+1}^{d-1} (-1)^j {j+1 \choose k+1} f_{j}(\A \cap \B).
\end{align*}
If $k$ is the same parity as $(d-1)$, this simplifies to
\[(-1)^{d-1}f_k(\partial\C) - \sum_{j=k+1}^{d-1}(-1)^j {j+1 \choose k+1} f_{j}(\A \cap \B).\]
If $k$ is different parity than $(d-1)$, we simplify to get 
\[(-1)^{d-1}(f_k(\partial\C) + 2f_k(\A \cap\B)) - \sum_{j=k+1}^{d-1} (-1)^j {j+1 \choose k+1} f_{j}(\A \cap \B).\]
In either case, if $k > l$, $f_k(\A \cap\B)) = 0$ and $\displaystyle{ \sum_{j=k+1}^{d-1} (-1)^j {j+1 \choose k+1} f_{j}(\A \cap \B)= 0}$.  Therefore, the Dehn-Sommerville relations hold for $l < k \leq d-2$.  If $k$ is the same parity as $(d-1)$, we only need $k \geq l$ since no $f_k(\A \cap \B)$ terms are involved.

If we assume $\A$ and $\B$ satisfy the Perles relations for $0 \leq k \leq d-1$ we have the following:
\begin{align*}
Pe_{k} (\C) & = Pe_{k}(\A) + Pe_{k}(\B)\\
&= (-1)^d[\alpha_{k}(\A) - f_k(\A)] +(-1)^d[\alpha_{k}(\B) - f_k(\B)] \\
&= (-1)^d\big[(\alpha_{k}(\A) + \alpha_{k}(\B)) - (f_k(\A)+f_k(\B))\big] \\
&=(-1)^{d}[\alpha_{k}(\C) - (f_k(\partial\C)+f_k(\A \cap \B))],
\end{align*}
where the last equality follows since $\int(\A \cap \B)= \emptyset$.  Therefore the $k^{\text{th}}$ Perles relation will hold on $\C$ if and only if $f_k(\A \cap \B) = 0$.  This will happen whenever $l<k\leq d-1$.
\end{proof}

\section{Conjectures}
This chapter shows that many gluings of geometric polytopal complexes or geometric simplicial complexes still satisfy the Gram, Dehn-Sommerville, and Perles relations.  With these results in mind, we make a few conjectures.  

Since gluings along lower dimensional complexes do not preserve the ratio of the Euler and angle characteristics, we consider a subclass of complexes which will include \compname  complexes.  A $d$-complex is said to be \emph{strongly connected} if any two $d$-faces can be connected by a path through faces of dimension $d$ or $d-1$.  
Then we make the following conjecture:
\begin{conjecture}
If $\C$ is a pure, strongly connected, $d$-dimensional complex embedded in $R^d$, where $d$ is odd, then $\chi_{\alpha}(\C) = \frac{1}{2}\chi(\C)$.
\end{conjecture}
This would be a further step toward showing that the angle characteristic is invariant under PL-homeomorphisms. 

Since the Perles equations also seem to be satisfied where the Dehn-Sommerville relations are satisfied, we make the following conjecture.
\begin{conjecture}
If $\C$ is a semi-Eulerian $(d-1)$-complex embedded in $R^d$, then $\C$ satisfies the Perles relations.
\end{conjecture}
The conjecture encompasses any embedded semi-Eulerian complex, not just those that are the boundary of a geometric simplicial $d$-complex. This may or may not be more general, but in either case, angle sums are still well-defined since $\C$ will be closed and have a defined $d$-dimensional interior.

%Exactly which complexes can we make this way? TIE IN CONSTRUCTIBLE AND SHELLABLE

%Thesis Chapter on Angles sums in different geometries
\chapter{Angle Sums in Other Geometries}

Thus far, all the results have concerned Euclidean polytopes and polytopal complexes.  In this chapter we will consider whether and how these results generalize to spherical and hyperbolic polytopes and polytopal complexes.

\section{Spherical and Hyperbolic Polytopes}

A \emph{pointed cone} is the intersection of linear hyperplanes, where the intersection of all the hyperplanes is exactly the origin.  A \emph{spherical $d$-polytope} is the intersection of a pointed cone in $\R^{d+1}$ with $\S^d$.  We could also think of this as the intersection of at least $d+1$ hemispheres in $\S^d$ which share a common point. If we consider the great spheres on $\S^d$ as hyperplanes and the hemispheres as the half-spaces they define, this definition is analogous to the intersection of half-spaces definition given for Euclidean polytopes.  If we choose an affine hyperplane parallel to a supporting hyperplane of a pointed cone at the origin, its intersection with the cone is a Euclidean $d$-polytope that has the same combinatorial structure as the spherical polytope defined by the cone.  Therefore it is clear that the Euler Relation holds for spherical polytopes and the Dehn-Sommerville Relations hold for simplicial spherical polytopes.  

A \emph{hyperbolic $d$-polytope} is formed by the intersection of hyperplanes in $\H^d$, just as Euclidean polytopes are.  In the Klein model of $\H^d$, hyperbolic space is viewed as a $d$-dimensional ball - as in the more familiar Poincar\'e model -  but where the hyperplanes are modeled by Euclidean hyperplanes. %PICTURE?
Therefore there is a one-to-one correspondence between hyperbolic polytopes and Euclidean polytopes where corresponding polytopes share the same combinatorial structure. As a result, the Euler relation holds for hyperbolic polytopes and the Dehn-Sommerville relations hold for simplicial hyperbolic polytopes.  

Therefore, we can summarize the relations on the $f$-vectors of polytopes in all the standard geometries:
\begin{theorem}\label{Euler & DS for all polytopes}
Let $\P$ be a $d$-polytope in $\E^d$, $\S^d$ or $\H^d$ and let $f_{i}(\P)$ be the number of faces of $\P$ of dimension $i$.  Then $\P$ satisfies the Euler relation:
\[\sum_{i=0}^{d-1} (-1)^{i}f_{i}(\P) = 1 + (-1)^{d+1}.\]
If, in addition, $\P$ is simplicial, then $\P$ satisfies the Dehn-Sommerville relations:
\[\sum_{j=k}^{d-1} (-1)^{j} {j+1\choose k+1} f_{j}(\P) =
(-1)^{d-1}f_{k}(\P) \quad \text{for }-1 \leq k \leq d-1.\]
\hfill $\Box$
\end{theorem} 

\section{Gram Relations on Spherical and Hyperbolic Polytopes}

We can define angle sums for spherical and hyperbolic polytopes in the same way they are defined for Euclidean polytopes.  For each face $\F$ of $\P$, we pick an interior point of the face and center a small sphere at that point. The sphere is not necessarily a Euclidean sphere, but rather intrinsic to each geometry, consisting of a set of points equidistant from the point.  Then the interior angle at $\F$ is the fraction of the sphere which is contained within the polytope. Since each of the geometries has constant curvature, the fraction of the sphere is independent of the size of the sphere since the edges and faces follow hyperplanes.  For example, consider the spherical triangle $T$ on $\S^2$ with right angles at all vertices so it is one-eighth of the sphere.  Then each angle at a vertex is $\frac{1}{4}$, and each angle at an edge is $\frac{1}{2}$.  Therefore the angle sums for this triangle are $\alpha_{0}(T) = \frac{3}{4}$ and $\alpha_{1}(T) = \frac{3}{2}$.   

\begin{figure}
\begin{center}
\includegraphics[width=4in]{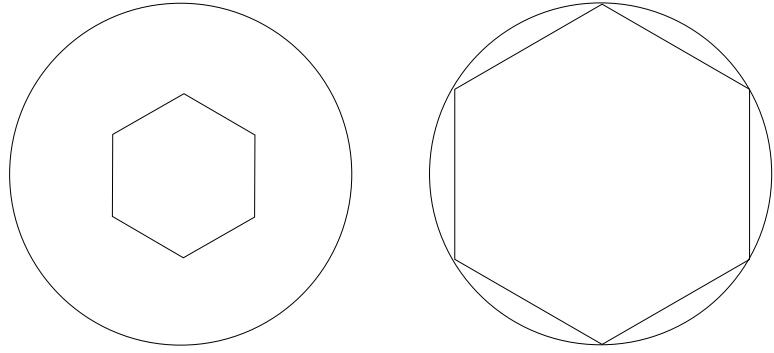}
\caption[Hexagons in the Klein model.]{Two hyperbolic hexagons in the Klein model, the first with non-zero angles and the second with vertex angles of measure 0 since the vertices are points at infinity.  The corresponding Euclidean hexagons have the same angles.}\label{Kleinpolygons}
\end{center}
\end{figure}

%PICTURE??
However, the case for the Gram and Perles relations is not so straightforward as that for the Euler and Dehn-Sommerville relations.  The angles of the corresponding Euclidean and spherical polytopes do not correspond as the combinatorial structure did; we can choose affine hyperplanes whose intersection with the cone would give us many different Euclidean polytopes with the same combinatorics but varying angles, and these angles will not agree with those of the spherical polytope. Similarly, Euclidean polytopes with constant angles and combinatorics can be made via dilation, but since the Klein model is not conformal, the corresponding hyperbolic polytopes will have differing angles even though the combinatorics is constant, as illutrated in Figure \ref{Kleinpolygons}. In fact, we can dilate the Euclidean polytope far enough that it will not fit in the Klein model that supported the original correspondence!

\singlespacing
\normalspacing

Despite these differences, there are similar results.  The major difference is a term which involves the volume of the polytope.  The introduction of this term is not surprising since the area of polygons in $\H^2$ and $\S^2$ is determined by the sum of the vertex angles.  Therefore, for both spherical and hyperbolic $\P$, we define a normalized volume,
\[\alpha_{-1}(\P) = \frac{\vol(\P)}{\vol(\S^d)},\]
where for spherical $\P$ we divide by the volume of the underlying sphere, therefore viewing the polytope as a certain fraction of the whole sphere, and for hyperbolic $\P$ we divide by the volume of a Euclidean sphere of the same radius as the ambient hyperbolic space.  In either case, the normalization adjusts for different curvature, so we can assume unit curvature.  Writing this normalized volume as $\alpha_{-1}(\P)$ is consistent with our earlier definition of $\alpha_{-1}(\P)=0$ for Euclidean $\P$.  This could be thought of as normalized volume by viewing $\E^d$ as a $d$-sphere of infinite radius. Another interpretation of $\alpha_{-1}(\P)$ for $\P$ spherical or Euclidean is as the angle at the empty face which is placed at the center of the sphere.  In either case, the angle subtended is the normalized volume, which is negligible in the Euclidean case.  For this chapter, we will extend the $\alpha$-vector to include $\alpha_{-1}(\P)$: $\alpha(\P) = \left(\alpha_{-1}(\P), \alpha_{0}(\P), \ldots, \alpha_{d}(\P)\right)$.  With this definition, there are Gram relations for hyperbolic and spherical polytopes:
\begin{theorem}[Sommerville]\label{spherical Gram}
Let $\P$ be a spherical $d$-polytope.
Then $\P$ satisfies the following relation on its angle sums:
\[\sum_{i=0}^{d}(-1)^{i} \alpha_{i}(\P)= \big(1+(-1)^d\big)\alpha_{-1}(\P).\]
\end{theorem}
This relation is usually called the \emph{Sommerville relation} \cite{So, McM2}.  It is frequently stated as a theorem on the pointed polyhedral cone in $\E^{d+1}$ used to define the spherical polytope and usually proved using these cones.

\begin{theorem}[Heckman]\label{hyperbolic Gram}
Let $\P$ be a hyperbolic $d$-polytope.
Then $\P$ satisfies the following relation on its angle sums:
\[\sum_{i=0}^{d}(-1)^{i} \alpha_{i}(\P)= (-1)^{\frac{d}{2}}\big(1+(-1)^d\big)\alpha_{-1}(\P).\]
\end{theorem}
The term on the right hand side is zero when $d$ is odd, so the difficulties of $-1$ raised to a fractional power are averted.  Heckman cites Hopf \cite{Hopf} as the first to extend the spherical result for simplices to a hyperbolic simplex.

Since all of the Gram theorems are very similar, we will rewrite this as one theorem, following Heckman \cite{He}.
\begin{theorem}[Gr\"unbaum, Sommerville, Heckman]\label{all-geometries Gram}
Let $\P$ be a Euclidean, spherical or hyperbolic $d$-polytope.  Define a curvature indicator 
\[\e = \begin{cases}
        1 \qquad \text{if } \P \text{ is spherical}\\
        0 \qquad \text{if } \P \text{ is Euclidean}\\
        -1 \qquad \text{if } \P \text{ is hyperbolic}.
\end{cases}\]
Then $\P$ satisfies the following relation on its angle sums:
\begin{equation}\label{GRAM}
\sum_{i=0}^{d}(-1)^{i} \alpha_{i}(\P)= (\e)^{\frac{d}{2}}\big(1+(-1)^d\big)\alpha_{-1}(\P).
\end{equation}
\end{theorem}

We will usually refer to all three of these relations in the different geomtries as the Gram relation for convenience, relying on context for the particular form.  If it is necessary to differ between these relations and the basic Gram relation on Euclidean polytopes, we will call these the generalized Gram relations.  Following Heckman \cite{He}, we will prove this theorem using a normalized form of the Schl\"afli Differential formula on simplices.  

\begin{theorem}[Heckman \cite{He}]\textbf{Normalized Schl\"afli Differential Formula} \label{Schlafli Differential formula}
Let $\D$ be a spherical or hyperbolic $d$-simplex.  Let $\F$ be a codimension 2 face of $\D$.  By convention, let $\alpha_{-1}(\cdot) = 1$. Then the partial derivative of the normalized volume of $\D$ with respect to the dihedral angle at $\F$ is a multiple of the $(d-2)$-dimensional normalized volume of $\F$, i.e.,
\[\frac{\partial \alpha_{-1}(\D)}{\partial \alpha(\F, \D)} = \e \alpha_{-1}(\F) \]
where $\e$ is 1 if $\D$ is spherical and $-1$ is $\D$ is hyperbolic.  \end{theorem}
Notice that $\alpha_{-1}(\D)$ is measuring a $d$-dimensional normalized volume while $\alpha_{-1}(\F)$ is a $(d-2)$-dimensional normalized volume.  To differentiate between these throughout this section, we will write $\displaystyle{\alpha_{-1}^d(\D) = \frac{\vol_d(\D)}{\vol_d(S^d)}}$ where the $d$ emphasizes the dimension of $\D$. Then the theorem can be rewritten as $\displaystyle{\frac{\partial \alpha^d_{-1}(\D)}{\partial \alpha(\F, \D)} = \e \alpha_{-1}^{d-2}(\F) }$.  Since $\E^d$ can be viewed as a $d$-sphere of infinite radius, any Euclidean polytope has normalized volume 0, so the Normalized Schl\"afli Differential Formula also holds for Euclidean polytopes with $\e=0$.  

In order to consider other angle sums, we remember that interior angles of any simplex are normalized volumes of spherical simplices.  That is, 
\[\alpha(\G, \D) = \frac{\vol_{d-1}(S^{d-1} \cap \D)}{\vol_{d-1}(S^{d-1})}= \alpha_{-1}^{d-1}(S^{d-1}\cap\D).\] 
If $\G$ is a $(d-j)$-face of $\D$, $\G$ is the intersection of $\D$ and $j$ hyperplanes, each of which is the affine span of a facet.  Let $J$ be the set of these hyperplanes. Then $J$ defines a cone $\D^J$ with apex $\G$ which contains $\D$ and the angle at $\D$ is 
the normalized volume of the spherical simplex determined by $\D^J$. In fact, we can write
\[\alpha(\G, \D) = \frac{\vol_{|J|-1}(S^{|J|-1} \cap \D)}{\vol_{|J|-1}(S^{|J|-1})} = \alpha_{-1}^{|J|-1}\left(S^{|J|-1} \cap \D^J\right),\]
where the sphere is in the normal space to the face.  This equivalence also holds for $\alpha_{-1}^{-1}(\cdot) = \alpha(\D, \D) = 1$. This allows us to translate between angle sums and normalized volumes of spherical simplices.
 \begin{lemma}\label{derivative of angle sum}
Let $\D$ be a $d$-simplex and $\F$ a codimension 2 face of $\D$.  Then \begin{align*}\frac{\partial \alpha_k(\D)}{\partial \alpha(\F, \D)} = \begin{cases}
\e\alpha_{-1}(\F) \qquad &\text{if } k=-1\\
\alpha_k(\F) \qquad &\text{if } 0\leq k \leq d-2\\
0 \qquad &\text{if } k > d-2
\end{cases}
\end{align*}
where $\e$ is 1 if $\D$ is spherical, 0 if $\D$ is Euclidean, and $-1$ if $\D$ is hyperbolic.  
\end{lemma}
\begin{proof}
We use the Normalized Schl\"afli Differential Formula.
When $k = -1$, this follows directly from the formula.

If $S$ is the set of supporting hyperplanes $H_i$, $i=1,2,\ldots d+1$, which define the facets of $\D$, then $\emptyset \neq J \subsetneq S$ defines a proper face $\G=\left(\bigcap_{i \in J}H_i\right)\cap \D$ where $\dim \G = d-|J|$. The set of hyperplanes in the affine span of $\G$ which define the facets of $\G$ is $\{H_i \cap \left(\bigcap_{i \in J}H_i\right): H_i \in S \setminus J\}$.  Therefore, the set of faces of $\G$ can also be identified with subsets of $S\setminus J$, and we will write $\G^{K}$ for the face of $\G$ defined by $K \subseteq S\setminus J$.  Therefore, if $k \geq 0$, we can write
\begin{equation*}
\begin{split}
\frac{\partial \alpha_{i}(\D)}{\partial \alpha(\F, \D)}
& = \sum_{i\text{-faces }\F^i} \frac{\partial \alpha(\F^i, \D)}{\partial \alpha(\F, \D)}\\
& = \sum_{I \subseteq S, |I| = d-i} \frac{\alpha_{-1}^{|I|-1}(S^{|I|-1} \cap \D^I)}{\partial \alpha(\F, \D)}\\
& = \sum_{J \subset I \subseteq S, |I| = d-i} \frac{\alpha_{-1}^{|I|-1}(S^{|I|-1} \cap \D^I)}{\partial \alpha(\F, \D)} \\ & \qquad \qquad \qquad \text{since other terms reduce to 0}\\
& = \sum_{J \subset I \subseteq S, |I| = d-i} \frac{\alpha_{-1}^{|I|-1}(S^{|I|-1} \cap \D^I)}{\partial \alpha(\F, \D)}  \\
& = \sum_{K \subsetneq S \setminus J, |K| = d-i-2}\alpha_{-1}^{|K|-1}(S^{|K|-1} \cap \F^K),
\end{split}
\end{equation*}
where the last step follows from the Normalized Schl\"afli Differential Formula where $\e = 1$ for a spherical simplex.  $S^{|K|-1} \cap \F^K$ is the face of $S^{|I|-1} \cap \D^I$ determined by $J$. Therefore, this last term is exactly $\alpha_i(\F)$.

Now, if $d-2<i = \dim (\F)$, $\alpha_{i}(\F) = 0$. This also follows since the sets $I$ have cardinality $|I| < d-(d-2) = 2$, so the set of $I$ containing $J$ is empty. 
\end{proof}

\begin{proof}[Proof of Theorem \ref{all-geometries  Gram}] 
This will be proved in a manner analogous to Gr\"unbaum's proof for the Euclidean Gram relation:  first the theorem is proved for simplices; second we show that if we decompose the polytope into pyramids, each the convex hull of a facet and an interior point where each pyramid follows the Gram relation, so does the whole;  third, we show that every $d$-pyramid follows the Gram relation.  Taken together, this shows that every spherical or hyperbolic polytopes satisfies the corresponding theorem.

The proof given in Chapter 1 of the second and third steps is sufficient since the proof is dependent only on the additive nature of interior angles and the Euler relation, which holds in all constant curvature geometries. Therefore, we only need to prove the theorem for simplices.

We will show the result for simplices by induction on $d$.  For $d=1$, the statement is obvious for simplices.  For $d=2$, we recall that if we have a hyperbolic triangle $T$ in $\S^2$ or $\H^2$ (of unit positive or negative curvature respectively) with angles $\alpha$, $\beta$ and $\gamma$, then $\vol_2(T) = \e((\alpha+\beta+\gamma) - \pi)$ \cite{Hend}. Dividing both sides by $2\pi$ and using that $\vol_2(S^2) = 4\pi$,  this becomes $2\alpha_{-1}(T) = \e\left(\alpha_0(T) - \frac{1}{2}\right)$.  This agrees with the Gram relation since $\alpha_1(T) = \frac{3}{2}$ and $\alpha_2(T) = 1$.  Likewise, for a triangle in $\E^2$, $\alpha_0(T)=\frac{1}{2}$, $\alpha_1(T) = \frac{3}{2}$ and $\alpha_2(T) = 1$, so the relation is satisfied.

Now suppose $d \geq 3$. Let $\F$ be a codimension 2 face of $\D$, defined by a set of hyperplanes $J \subsetneq S$.  We will check that the derivatives of both sides of the relation with respect to the dihedral angle, $\alpha(\F, \D)$, are equal.  This implies that the formula is correct up to an additive constant.  By Lemma \ref{derivative of angle sum}  the derivative of the right hand side is
\begin{equation*}
\begin{split}
\e^{\frac{d}{2}}\left(1+(-1)^d\right) \frac{\partial \alpha_{-1}(\P)}{\partial \alpha(\F, \D)}
& = \e^{\frac{d-2}{2}}\left(1+(-1)^d\right)\alpha_{-1}(\F)
\end{split}
\end{equation*}
and that of the left hand side is
\begin{equation*}
\begin{split}
\sum_{i=0}^{d}(-1)^{i} \frac{\partial \alpha_{i}(\D)}{\partial \alpha(\F, \D)}
& = \sum_{i=0}^{d}(-1)^{i} \alpha_{i}(\F).
\end{split}
\end{equation*}

Therefore the derivatives of the two sides give the Gram relation in dimension $d-2$ and we need only check to see if there is an additive constant in the relation in dimension $d$.  In the spherical case, we can take a simplex with all dihedral angles equal to $\frac{\pi}{2}$.  Then $\alpha_{-1}^d(\D) = 2^{-d-1}$ and the formula gives the correct identity \[\left(1+(-1)^d\right)2^{-d-1} = \sum_{k=0}^d{n+1 \choose k}\left(-\frac{1}{2}\right)^k.\] Therefore the theorem holds for spherical simplices.  We can take a set of spherical simplices of decreasing size (or, equivalently, spherical simplices of the same volume on increasingly larger spheres) to show the angle sum on the left hand side vanishes for a Euclidean simplex $\D$ since the normalized volume tends to 0.  Likewise, a limiting argument for hyperbolic simplices of decreasing size shows the additive constant is 0.
\end{proof}

\section{Perles Relations on Spherical and Hyperbolic Polytopes} 
We can also consider whether there are relations on simplicial spherical and hyperbolic polytopes analogous to the Perles relations.  In the case of simplicial spherical polytopes, we do have such relations, shown by Perles and Shephard \cite{PS}. %CHECK THIS!

\begin{theorem}\label{Spherical Perles Relations}
Let $\P$ be a simplicial spherical polytope.  Then 
\[Pe_k(\P) \equiv \sum_{j=k}^{d-1} (-1)^j {j+1 \choose k+1} \alpha_j(\P) = (-1)^{d} \left[\alpha_k(\P) - f_k(\P) \right]\]
for $-1 \leq k \leq d-1$, where $k=-1$ gives Sommerville's relation on spherical polytopes.
\end{theorem}
These relations are sometimes referred to as the Perles-Shephard relations, but since they have the same form in $\S^d$ as in $\E^d$, we will call both the Perles relations or generalized Perles relations for simplicity. Sometimes the Perles relations above are written as
\[\left((-1)^k - (-1)^d\right)\alpha_k (\P) + \sum_{j=k+1}^{d-1} (-1)^j {j+1 \choose k+1} \alpha_j(\P) = (-1)^{d-1} f_k(\P).\]

It is not known whether simplicial hyperbolic polytopes satisfy a similar relation.  The relations in the spherical and Euclidean cases are usually proved using cones in Euclidean space \cite{PS}, the symmetry of the sphere \cite{KS} or by the projection argument in Chapter 1.  None of these methods are usable to prove relations on simplicial hyperbolic polytopes. %CHECK REFERENCES

We will consider a few examples to determine what form such equations might take. Since the spherical Perles relations have the same form as the Euclidean Perles relations, it is not surprising that the examples of simplicial hyperbolic polytopes below satisfy a relation of the same form. In general, the $k=-1$ case of the Perles relations reduces to the Gram relations.  Therefore, we will consider only cases where $0 \leq k \leq d-1$. 

For $d=1$, the only hyperbolic polytope is a segment $\P$. Then $\alpha_0(\P) = 1$, counting both the vertices with an angle of $\frac{1}{2}$.  We also know $\alpha_1(\P) = 1$.  Then for $k=0$ then 
\[Pe_0(\P) \equiv \sum_{j=0}^{0} (-1)^j {j+1 \choose 1} \alpha_j(\P) = 1 =   (-1)^1\left(\alpha_0(\P) - f_0(\P)\right). \]

For $d=2$, we consider hyperbolic polygons.  Let $\P_n$ be a polygon with $n$ sides. Then $\alpha_1(\P_n) = \frac{n}{2}$.  Then for $k=0$, 
\[Pe_0(\P_n) \equiv \sum_{j=0}^{1} (-1)^j {j+1 \choose 1} \alpha_j(\P_n) = \alpha_0(\P_n) - n = (-1)^2\left(\alpha_0(\P_n) - f_0(\P_n)\right).\]
For $k=1$, 
\[Pe_1(\P_n) \equiv \sum_{j=1}^{1} (-1)^j {j+1 \choose 2} \alpha_j(\P_n) = -\frac{n}{2} = (-1)^2\left(\alpha_1(\P_n) - f_1(\P_n)\right).\]

For $d=3$, we can see the relations also hold for a simplex $\D$ for $k=0,1$.
For $k=0$, 
\begin{equation*}
\begin{split}
Pe_0(\D) & \equiv \sum_{j=0}^{2} (-1)^j {j+1 \choose 1} \alpha_j(\D)\\
& = \alpha_0(\D) - 2\alpha_1(\D) + 3\alpha_2(\D) \\
& = 2(\alpha_0(\D) - \alpha_1(\D) + \alpha_2(\D)) -\alpha_0(\D)+ \alpha_2(\D)\\
& = 2\alpha_3(\D) -\alpha_0(\D) + \alpha_2(\D) \quad \text{by the Gram relation} \\
& = - \alpha_0(\D)+4 \\
& = (-1)^3\left(\alpha_0(\D) - f_0(\D)\right).
\end{split}
\end{equation*}
For $k=1$, 
\begin{equation*}
\begin{split}
Pe_1(\D) & \equiv \sum_{j=1}^{2} (-1)^j {j+1 \choose 2} \alpha_j(\D) = -\alpha_1(\D)+ 3 \alpha_2(\D) \\
& = -\alpha_1(\D)+ 6 = (-1)^3\left(\alpha_1(\D) - f_1(\D)\right).
\end{split}
\end{equation*}

For any hyperbolic $d$-polytope $\P$ we can see that for $k=d-1$, 
\[Pe_{d-1}(\P) \equiv (-1)^{d-1}\alpha_{d-1}(\P) = (-1)^{d-1}\left(\alpha_{d-1}(\P) - f_{d-1}(\P)\right)\] since $\alpha_{d-1}(\P) = \frac{1}{2}f_{d-1}(\P)$.
%Therefore the Perles relations hold for $k=d-1$.
For $k=d-2$,
\begin{equation*}
\begin{split}
Pe_{d-2}(\P) & \equiv (-1)^{d-2}\alpha_{d-2}(\P) + d(-1)^{d-1}\alpha_{d-1}(\P) \\
& = (-1)^{d-2}\left(\alpha_{d-2}(\P) - f_{d-2}(\P)\right)
\end{split}
\end{equation*}
since $f_{d-2}(\P) = \frac{d}{2} f_{d-1}(\P) = d\alpha_{d-1}(\P)$ because $\P$ is simplicial. In all of these examples, the Perles relations hold, so we make the following conjecture:
\begin{conjecture}\label{hyperbolic Perles}
Let $\P$ be a simplicial hyperbolic polytope of dimension $d$.  Then, for $0 \leq k \leq d-1$,
\[Pe_k(\P) \equiv \sum_{j=k}^{d-1} (-1)^j {j+1 \choose k+1} \alpha_j(\P) = (-1)^{d} \left(\alpha_k(\P) - f_k(\P) \right).\]
\end{conjecture}

We will call these conjectured relations the hyperbolic Perles relations. Taking $k=-1$ in this relation does not give the Gram relation for hyperbolic polytopes since the curvature constant is missing.  However, if we define a variant on the $\alpha$-vector which takes the curvature into account, we can rewrite the Gram and Perles relations for the different geometries in a more homogenous way.  Define the $\alphah$-vector as 
\[\left(\alphah_{-1}(P), \alphah_{0}(P),  \ldots, \alphah_{d}(P)\right) =\left((\e)^{\frac{d}{2}}\alpha_{-1}(P), \alpha_{0}(P),  \ldots, \alpha_{d}(P)\right). \]
We can either think of the $(\e)^{\frac{d}{2}}$ term abstractly and only evaluate it in the context of the relations or, in the hyperbolic case, we can use $(\e)^{\frac{d}{2}} = \cos\left(\frac{d\pi}{2}\right)$ to avoid raising $-1$ to a fractional power.

Then the Gram relations for all the geometries become
\[\sum_{i=0}^{d}(-1)^{i} \alphah_{i}(\P)= \big(1+(-1)^d\big)\alphah_{-1}(\P).\]
The hyperbolic Perles relations are equivalent with $\alpha_i(\P)$ replaced with $\alphah_i(\P)$ and Conjecture \ref{hyperbolic Perles} may be extended to include the $k=-1$ case, which would agree with the Gram relation.  Then we can write the following conjecture which would include the Gram and Perles relations in all three geometries:
\begin{conjecture}\label{all-geometries Perles}
Let $\P$ be a Euclidean, spherical or hyperbolic $d$-polytope.  
Then for $-1 \leq k \leq d-1$
\[\sum_{j=k}^{d-1} (-1)^j {j+1 \choose k+1} \alphah_j(\P) = (-1)^{d} \left(\alphah_k(\P) - f_k(\P) \right).\]
\end{conjecture}

The examples above show that the conjecture is always true for  3-simplices and for $k=d-1,d-2$.   From this last, we can show the conjecture is true for $d=3$ using the following argument which is similar to the second step of Gr\"unbaum's proof of the Gram relation.  The same argument would prove the conjecture for simplicial $d$-polytopes of any geometry, assuming that the Perles relations were true for simplices of dimension $d$.

Let $\P$ be a simplicial polytope.  Then we choose a point $0$ in the interior of $\P$ and decompose $\P$ into simplices $\P_1, \P_2, \ldots, \P_m$ where $m = f_{d-1}(\P)$ and each $\P_i$ is the pyramid with one of the facets as a base and apex $0$.  Assuming the Perles relations on the simplices $\P_i$, $i=1\ldots m$, we note the following relations:
\begin{equation}\label{Fact 1}
\alpha_j(\P) = \sum_{i} \alpha_j(\P_i) - f_{j-1}(\P)
\end{equation}
and
\begin{equation}\label{Fact 2} 
\sum_i f_j(\P_i) = {d \choose j} f_{d-1}(\P)  +  {d \choose j+1} f_{d-1}(\P) \quad 0 \leq i \leq d.
\end{equation}

The angles at the $j$-faces of $\P$ will be included in the angles at the $j$-faces of the $\P_i$, but the $\P_i$ also include new $j$-faces. If $\F$ is an interior $j$-face, the interior angles from $\F$ into the $\P_i$ will sum to 1 since angles in all directions are counted.  Therefore, the excess counted in all the angles at the $j$-faces of the $\P_i$ is exactly the number of new interior $j$-faces created in the decomposition.  This number is $f_{j-1}(\P)$ since each interior $j$-face is the convex hull of $0$ and a $(j-1)$-face of $\P$. This gives \eqref{Fact 1}.

To get \eqref{Fact 2} we count the $j$-faces of the $\P_i$, which can occur in two ways: either as a $j$-face of $\P$ or as the pyramid over a $(j-1)$-face of $\P$ with apex $0$.  For each facet $\F_i$, there are ${d \choose j+1}$ $j$-faces and ${d \choose j}$ $(j-1)$-faces of $\P$ counted as $j$-faces in $\P_i$.  Therefore, the sum counts the number of $j$-faces in all the $\P_i$.

Therefore, using \eqref{Fact 1}:
\begin{equation*}
\begin{split}
Pe_k(\P) & \equiv \sum_{j=k}^{d-1} (-1)^j {j+1 \choose k+1} \alpha_j(\P) \\
& = \sum_i \sum_{j=k}^{d-1} (-1)^j {j+1 \choose k+1} \alpha_j(\P_i) - \sum_{j=k}^{d-1} (-1)^j {j+1 \choose k+1} f_{j-1}(\P).
\end{split}
\end{equation*} 
Then, assuming the Perles relations on the $\P_i$ we get
\begin{equation*}
\begin{split}
Pe_k(\P) & = \sum_i(-1)^d\left(\alpha_k(\P_i) - f_k(\P_i)\right) - \sum_{j=k}^{d-1} (-1)^j {j+1 \choose k+1} f_{j-1}(\P)\\
& = \sum_i(-1)^d\left(\alpha_k(\P_i) - f_k(\P_i)\right) \\
& \qquad \qquad - \Bigg[\sum_{j=k}^{d-1} (-1)^j {j \choose k} f_{j-1}(\P) +\sum_{j=k}^{d-1} (-1)^j {j \choose k+1} f_{j-1}(\P) \Bigg]\\
& = \sum_i(-1)^d\alpha_k(\P_i) - \sum_i(-1)^d f_k(\P_i)  \\
& \qquad \qquad + \sum_{m=k-1}^{d-2} (-1)^m {m+1 \choose k} f_{m}(\P) +\sum_{n=k-1}^{d-2} (-1)^n {n+1 \choose k+1} f_{n}(\P) \\
& = (-1)^d\left(\alpha_k(\P) + f_{k-1}(\P)\right)- \sum_i(-1)^d f_k(\P_i) \\
& \qquad \qquad + \sum_{m=k-1}^{d-1} (-1)^m {m+1 \choose k} f_{m}(\P) -(-1)^{d-1}{d \choose k} f_{d-1}(\P)  \\
& \qquad \qquad +\sum_{n=k}^{d-1} (-1)^n {n+1 \choose k+1} f_{n}(\P)  - (-1)^{d-1} {d \choose k+1} f_{d-1}(\P)
\end{split}
\end{equation*} 
 by \eqref{Fact 1}.  Then we can use the Dehn-Sommerville relations followed by \eqref{Fact 2} to simplify:
\begin{equation*}
\begin{split}
Pe_k(\P) & = (-1)^d\alpha_k(\P) + (-1)^d f_{k-1}(\P) - \sum_i(-1)^d f_k(\P_i)  \\
& \qquad \qquad + (-1)^{d-1}f_{k-1}(\P)  -(-1)^{d-1}{d \choose k} f_{d-1}(\P)  \\
& \qquad \qquad +(-1)^{d-1}f_k(\P)  - (-1)^{d-1} {d \choose k+1} f_{d-1}(\P) \\
& = (-1)^d\alpha_k(\P)+(-1)^{d-1}f_k(\P)  \\
& \qquad \qquad + (-1)^{d-1}\left[\sum_{i}f_k(\P_i)-{d \choose k} f_{d-1}(\P) -{d \choose k+1} f_{d-1}(\P) \right]\\
& = (-1)^{d}\left[\alpha_k(\P) - f_k(\P)\right]. 
\end{split}
\end{equation*} 
Therefore, to prove Conjecture \ref{all-geometries Perles} we need only prove the Perles relation for hyperbolic simplices.  

As noted above, the standard methods of proof for spherical and Euclidean Perles relations do not apply to hyperbolic simplices.  Trying to apply Heckman's method of proof does not work, since $f_k(\D)$ is a constant and therefore has derivative 0 rather than $f_k(\F)$.  This is true even though for a codimension 2 face $\F$ of $\D$, Lemma \ref{derivative of angle sum} shows that the derivative of $Pe_k(\D)$ with respect to $\alpha(\F, \D)$ is $Pe_k(\F)$ and the left sides match.  
     
\section{Hyperbolic and Spherical Polytopal Complexes}
Using the generalized Gram relations, we can consider analogous relations on spherical and hyperbolic polytopal complexes.  As earlier, we will define a \emph{geometric polytopal complex} as a connected, pure complex with polytopes of maximal dimension $d$ which is embedded in $d$-space, either $\E^d$, $\S^d$ or $\H^d$.  We will assume a complex is composed homogeneously of Euclidean, spherical or hyperbolic polytopes.  We will also define \emph{gluings} as before, allowing a gluing $\C = \A \glue \B$ if $\A$ and $\B$ are both complexes in the same underlying space, $\E^d$, $\S^d$ or $\H^d$.

Since the Euler relation and Euler characteristic are independent of the underlying geometry, Lemma \ref{independence of subdivision}  applies and the angle characteristic is still based on the geometric realization of a complex rather than the particular subdivision of the flats.  Similarly, Lemma \ref{valuations} still holds since it is only dependent on the topology of a gluing, and therefore 
\[\chi\left(\partial \C\right) = \chi\left(\partial \A\right) +\chi\left(\partial \B\right) - 2\chi\left(\int\left(\A \cap \B\right)\right) - \chi\left(\partial\left(\A \cap \B\right)\right)\] 
and
\[\chi_{\alpha}\left(\C\right) = \chi_{\alpha}\left(\A\right)+\chi_{\alpha}\left(\B\right) - \chi\left(\int\left(\A \cap \B\right)\right). \]

If a complex $\C$ satisfies the Gram relation, we can rewrite this as \begin{equation}
\chi_{\alpha}(\C) = (-1)^{d-1} + (\e)^{\frac{d}{2}}\big(1+(-1)^d\big)\alpha_{-1}(\C).
\end{equation}
We will consider what happens when gluing together complexes that satisfy the Gram relation, using the fact that the volume on the complexes is additive:
\begin{align*}
\chi_{\alpha}(\C) & = \chi_{\alpha}(\A)+\chi_{\alpha}(\B) - \chi(\int(\A \cap \B))\\
& = (-1)^{d-1} + (\e)^{\frac{d}{2}}\big(1+(-1)^d\big)\alpha_{-1}(\A) + (-1)^{d-1} \\
& \qquad \qquad + (\e)^{\frac{d}{2}}\big(1+(-1)^d\big)\alpha_{-1}(\B) - \chi(\int(\A \cap \B)) \\
& = 2(-1)^{d-1} + (\e)^{\frac{d}{2}}\big(1+(-1)^d\big)\alpha_{-1}(\C)- \chi(\int(\A \cap \B)).
\end{align*}
This proves the following:
\begin{lemma}\label{all-geometries relation}
Let $\C = \A \glue \B$, where $\A$ and $\B$ are both geometric polytopal complexes in either $\E^d$, $\S^d$ or $\H^d$.  Then if $\A$ and $\B$ satisfy the Gram relation, 
\[ \chi_{\alpha}(\C) = 2(-1)^{d-1} + (\e)^{\frac{d}{2}}\big(1+(-1)^d\big)\alpha_{-1}(\C)- \chi(\int(\A \cap \B)).\]
\end{lemma}

We now consider specific intersections, $\A \cap \B$. 
\begin{proposition}\label{all-geometries gluing along balls}
Let $\C = \A \glue \B$, where $\A$ and $\B$ are both geometric polytopal complexes in either $\E^d$, $\S^d$ or $\H^d$ and $\A \cap \B$ is the union of $m$ $(d-1)$-balls.  Then if $\A$ and $\B$ satisfy the Gram relation,
\[ \chi_{\alpha}(\C) = (-1)^{d-1} + (\e)^{\frac{d}{2}}\big(1+(-1)^d\big)\alpha_{-1}(\C) - (m-1)(-1)^{d-1}.\]
In particular, if $m=1$, $\C$ also satisfies the Gram relation.
\end{proposition}
\begin{proof}
If $K$ is a $(d-1)$-ball, $\chi(\int (K)) = (-1)^{d-1}$. Therefore, by Lemma \ref{all-geometries relation},
\begin{align*}
\chi_{\alpha}(\C) & = 2(-1)^{d-1} + (\e)^{\frac{d}{2}}\left(1+(-1)^d\right)\alpha_{-1}(\C)- \chi(\int(\A \cap \B)) \\
& = 2(-1)^{d-1} + (\e)^{\frac{d}{2}}\left(1+(-1)^d\right)\alpha_{-1}(\C) - m(-1)^{d-1} \\ 
& = (-1)^{d-1} + (\e)^{\frac{d}{2}}\left(1+(-1)^d\right)\alpha_{-1}(\C) - (m-1)(-1)^{d-1}.
\end{align*}
\end{proof}
In particular, this proposition states that stacked polytopes in any geometry satisfy the Gram relation, including those whose embedding is not convex.

\begin{proposition}\label{S/H gluing along an annulus}
Let $\C = \A \glue \B$, where $\A$ and $\B$ are both geometric polytopal complexes in either $\E^d$, $\S^d$ or $\H^d$ and $\A \cap \B$ is a union of $m$ disjoint $(d-1)$-dimensional annuli (i.e. closed $(d-1)$-balls with an open $(d-1)$-ball removed from the interior). Then if $\A$ and $\B$ both satisfy the Gram relation, 
\[ \chi_{\alpha} (\C)=  (2-m)(-1)^{d-1} + (\e)^{\frac{d}{2}}\left(1+(-1)^d\right)\alpha_{-1}(\C)+  m.\]
\end{proposition}
\begin{proof}
Since $\chi(\int(\A \cap \B)) = m(-1 + (-1)^{d-1})$, we can use Lemma \ref{all-geometries relation} to write
\begin{equation*}
\begin{split}
 \chi_{\alpha} (\C) & = 2(-1)^{d-1} + (\e)^{\frac{d}{2}}\left(1+(-1)^d\right)\alpha_{-1}(\C)- \chi(\int(\A \cap \B)) \\
& = 2(-1)^{d-1} + (\e)^{\frac{d}{2}}\left(1+(-1)^d\right)\alpha_{-1}(\C)+  m (1 +(-1)^{d}) \\
& = (2-m)(-1)^{d-1} + (\e)^{\frac{d}{2}}\left(1+(-1)^d\right)\alpha_{-1}(\C)+  m.
\end{split}
\end{equation*}
\end{proof}

\begin{proposition}\label{S/H gluing along a sphere}
Let $\C = \A \glue \B$, where $\A$ and $\B$ are both geometric polytopal complexes in either $\E^d$, $\S^d$ or $\H^d$ and $\A \cap \B$ is a union of $m$ disjoint $(d-1)$-dimensional spheres. Then if $\A$ and $\B$ both satisfy the Gram relation,
\[ \chi_{\alpha} (\C)= (2-m)(-1)^{d-1} + (\e)^{\frac{d}{2}}\big(1+(-1)^d\big)\alpha_{-1}(\C) - m .\]
\end{proposition}
\begin{proof}
If $K$ is a $(d-1)$-sphere, $\chi(\int (K)) = \chi(K) = 1+(-1)^{d-1}$. Therefore, by Lemma \ref{all-geometries relation},
\begin{align*}
\chi_{\alpha}(\C) & = 2(-1)^{d-1} + (\e)^{\frac{d}{2}}\big(1+(-1)^d\big)\alpha_{-1}(\C)- \chi(\int(\A \cap \B)) \\
& = 2(-1)^{d-1} + (\e)^{\frac{d}{2}}\big(1+(-1)^d\big)\alpha_{-1}(\C) - m(1+(-1)^{d-1}) \\ 
& = (2-m)(-1)^{d-1} + (\e)^{\frac{d}{2}}\big(1+(-1)^d\big)\alpha_{-1}(\C) - m.
\end{align*}
\end{proof}

In all of these constructions, since the relations between $\chi(\C)$, $\chi(\A)$ and $\chi(\B)$ and between $\chi_{\alpha}(\C)$, $\chi_{\alpha}(\A)$ and $\chi_{\alpha}(\B)$ are dependent only on the topology of $\A \cap \B$, the difference between the angle characteristic of $\C$ and the angle characteristic on $\A$ and $\B$ is half the corresponding difference in the Euler characteristics.  Unlike Euclidean complexes, however, the resulting complexes will not have a ratio of 1:2 between the angle characteristic and the Euler characteristic, since polytopes in $\S^d$ and $\H^d$ do not have this ratio. 

\begin{proposition}\label{S/H gluing along lower dimensional}
Let $\C = \A \glue \B$, where $\A$ and $\B$ are both geometric polytopal complexes in either $\E^d$, $\S^d$ or $\H^d$ and $\A \cap \B$ is of dimension $l \leq d-2$. Then if $\A$ and $\B$ both satisfy the Gram relation,   
\[ \chi_{\alpha} (\C)= 2(-1)^{d-1} + (\e)^{\frac{d}{2}}\big(1+(-1)^d\big)\alpha_{-1}(\C).\]
\end{proposition}
\begin{proof}
If $\A \cap \B$ has dimension less than $d-1$, $\int(\A \cap \B) = \emptyset$. Therefore, by Lemma \ref{all-geometries relation},
\begin{align*}
\chi_{\alpha}(\C) & = 2(-1)^{d-1} + (\e)^{\frac{d}{2}}\big(1+(-1)^d\big)\alpha_{-1}(\C)- \chi(\int(\A \cap \B)) \\
& = 2(-1)^{d-1} + (\e)^{\frac{d}{2}}\big(1+(-1)^d\big)\alpha_{-1}(\C). \end{align*}
\end{proof}

%Perles-type results on complexes

There are also analogs of the Perles relation on simplicial geometric polytopal complexes in other geometries.  As the topology of the complexes rather than the geometry determines the relationship between the Dehn-Sommerville and Perles operators on $\A$, $\B$, and $\C = \A \glue \B$, the relations in Lemma \ref{simplicial valuations} still hold.  In the same way, the earlier results on the Dehn-Sommerville operator, $DS_k(\partial\C)$, hold in all the geometries.  In particular, Klee's Theorem \ref{Klee DS} on semi-Eulerian complexes applies to semi-Eulerian spherical and hyperbolic geometric polytopal complexes since it is based on the Euler characteristics of links only, and hence on the topology of the complexes rather than the underlying geometry. 

Since the Perles equations are identical on both Euclidean and spherical simplicial complexes and are conjectured to be the same for hyperbolic simplicial complexes, we can state the following, which has the same proof as Theorems \ref{DS on ball and annuli gluings} and \ref{DS on sphere-gluings}.

\begin{theorem}
Let $\A$ and $\B$ be geometric simplicial complexes in either $\E^d$, $\S^d$, or $\H^d$, $\C = \A \glue \B$, and suppose $\A$ and $\B$ satisfy the Perles relations \ref{all-geometries Perles}.  If $\A \cap \B$ is a disjoint union of $(d-1)$-balls, $(d-1)$-annuli, and semi-Eulerian $(d-1)$-complexes then $\C$ satisfies the Perles relations.
\[Pe_{k} (\C) = (-1)^{d}\left[\alpha_{k}(\C) - f_k(\C)\right].\]\hfill $\Box$
\end{theorem}

\end{document}